\documentclass[12pt]{article}
\usepackage{amsmath,amssymb}
\usepackage[dvipdfmx]{graphicx}
\usepackage{bm}
\newtheorem{tm}{Theorem}[section]
\newtheorem{lm}[tm]{Lemma}
\newtheorem{co}[tm]{Corollary}
\newtheorem{re}[tm]{Remark}

\newtheorem{pr}[tm]{Proposition}

\newtheorem{ex}[tm]{Example}
\makeatletter
\newcommand{\subscripts}[3]{%
  \@mathmeasure\z@\displaystyle{#2}%
  \global\setbox\@ne\vbox to\ht\z@{}\dp\@ne\dp\z@
  \setbox\tw@\box\@ne
  \@mathmeasure4\displaystyle{\copy\tw@_{#1}}%
  \@mathmeasure6\displaystyle{{#2}_{#3}}%
  \dimen@-\wd6 \advance\dimen@\wd4 \advance\dimen@\wd\z@
  \hbox to\dimen@{}\mathop{\kern-\dimen@\box4\box6}%
}
\makeatother
\newcommand{\qed}{~~\hbox{\rule{4pt}{8pt}}}

\def\Re{\mathop{\rm Re}\nolimits}
\def\Im{\mathop{\rm Im}\nolimits}
\makeatletter
 
 \@addtoreset{equation}{section}
\makeatother
\begin{document}
\setlength{\baselineskip}
{15.5pt}
\title{
Long time asymptotics of non-symmetric \\
random walks 
on crystal lattices}
\author{\Large
{Satoshi Ishiwata
\footnote{Department of Mathematical Sciences, Faculty of Science, Yamagata University,
1-4-12, Kojirakawa, Yamagata 990-8560, Japan
(e-mail: {\tt ishiwata@sci.kj.yamagata-u.ac.jp})}
\footnote{Partially supported by Grant-in-Aid for Young Scientists (B)(21740034, 25800034), JSPS}, 
Hiroshi Kawabi 
\footnote{Department of Mathematics, Faculty of Science, Okayama University,
3-1-1, Tsushima-Naka, Kita-ku, Okayama 700-8530, Japan (e-mail: {\tt{kawabi@math.okayama-u.ac.jp}})}
\footnote{Partially supported by Grant-in-Aid for Young Scientists (B)(23740107), Grant-in-Aid for
Scientific Research (C)(26400134), JSPS}}~~and
{\Large{Motoko Kotani
\footnote{Advanced Institute for Material Research and Mathematical Institute, Tohoku University, 2-1-1, Katahira, Aoba-ku, 
Sendai 980-8577, Japan (e-mail: {\tt{m-kotani@m.tohoku.ac.jp}})}
\footnote{Partially supported by Grant-in-Aid for
Scientific Research (A)(24244004), JSPS}}}
}
\date{   }
\maketitle 
%
%
\begin{abstract}
In the present paper, we study long time asymptotics of
non-symmetric random walks on crystal lattices from 
a view point of discrete geometric analysis due to Kotani and Sunada \cite{KS00, S}.
We observe that the Euclidean metric associated with the standard realization
of the crystal lattice, called the Albanese metric, naturally appears in
the asymptotics.
In the former half of the present paper, we establish two kinds of 
(functional) central limit theorems for random walks.
We first show that the Brownian motion on the Euclidean space with the
Albanese metric appears as the scaling limit of the
usual central limit theorem for the random walk. 
Next we introduce a family of random walks which interpolates between the
original non-symmetric random walk and the symmetrized one. We then 
capture the Brownian motion with a constant drift of the asymptotic direction
on the Euclidean space with the Albanese metric associated with the
symmetrized random walk through another kind of central limit theorem 
for the family of random walks.
In the latter half of the present paper, we give a spectral geometric proof of
the asymptotic expansion of the $n$-step transition probability 
for the non-symmetric random walk. This asymptotic
expansion is a refinement of the local central limit theorem 
obtained by Sunada \cite{Sunada-Lecture} and is a generalization of 
the result in \cite{KS00} for symmetric random walks on crystal
lattices to non-symmetric cases.
\vspace{2mm} \\
{\bf{2010 AMS Classification Numbers:}}~60J10, 60F05, 60G50, 60B10.
\vspace{2mm} \\
{\bf{Keywords:}}~Crystal lattice, non-symmetric random walk, central limit theorem,
asymptotic expansion, Albanese metric, modified standard realization.
\end{abstract}
\section{Introduction}\label{Introduction}
Let $X=(V,E)$ be an oriented, locally finite connected graph (which may have multiple edges and loops).
For an oriented edge $e \in E$, the origin and the terminus of $e$ are denoted by 
$o(e)$ and $t(e)$, respectively. The inverse edge of $e\in E$ is denoted by $\overline{e}$.
Let $E_{x}=\left\{ e \in E \vert~ o(e)=x \right\}$ be the set of edges with $o(e)=x\in V$. 
A path $c$ of $X$ of length $n$ is a sequence $c=(e_{1},\ldots, e_{n})$ of oriented edges 
$e_{i}$ with $o(e_{i+1})=t(e_{i})$ ($i=1,\ldots, n-1$). We denote by 
${\Omega}_{x,n}(X)$ ($x\in V, n\in \mathbb N \cup \{\infty \}$)
the set of all paths of length $n$ for which origin $o(c)=x$. For simplicity, we also write
${\Omega}_{x}(X):={\Omega}_{x,\infty}(X)$.

A random walk on $X$ is a stochastic process with values in $X$ 
characterized effectively by a transition probability,
a non-negative function $p:E\to \mathbb{R}$ satisfying
\begin{equation}
\sum_{e \in E_{x}} p(e)=1 \quad (x\in V),  \qquad  
p(e)+p({\overline e})>0 \quad (e \in E),
\label{p-def-intro}
\end{equation}
where $p(e)$ stands for the probability that a particle at $o(e)$ moves to 
$t(e)$ along the edge $e$ in one unit time. 
The transition operator $L$ on $X$ 
associated with the random walk 
is defined by
\begin{equation*}
Lf(x):=\sum_{e \in E_x} p(e) f(t(e))
\qquad (x\in V).
\end{equation*}
The $n$-step transition 
probability $p(n,x,y)$  ($n\in \mathbb N$, $x,y\in V$) 
is defined by
\begin{equation}
p(n,x,y):=\sum_{c=(e_{1},\ldots, e_{n})} p(e_{1})\cdots p(e_{n}), 
\label{nstep-p-intro}
\end{equation}
where the sum is taken over all paths 
$c =(e_{1},\ldots, e_{n})$ of length $n$ with the origin $o(c)=x$ and
the terminus $t(c)=y$.
We mention
\begin{equation*}
L^{n}f(x)=\sum_{y \in V} p(n,x,y) f(y) \qquad (x\in V).
\end{equation*}
In a natural manner, the transition probability $p$ induces
the probability measure ${\mathbb P}_{x}$ on the set 
${\Omega}_{x}(X)$.
The random walk associated with $p$ is 
the time homogeneous 
Markov chain $({\Omega}_{x}(X), {\mathbb P}_{x}, \{w_{n}\}_{n=0}^{\infty} )$
with values in $X$ defined by
$$w_{n}(c):=o(c(n+1)) \qquad (n=0,1,2,\ldots,~ c\in {\Omega}_{x}(X)), $$
where $c(n)$ is the $n$th edge of 
the infinite path 
$c\in {\Omega}_{x}(X)$.
If, in addition, there exists a positive function $m:V\rightarrow \mathbb{R}$ 
such that 
\begin{equation*}
p(e)m(o(e))=p(\overline{e})m(t(e)) \qquad (e \in E),
\end{equation*}
then the random walk is said to be  {\it symmetric} (or {\it reversible}), and the function $m$ is called
a {\it{reversible measure}} for the random walk.
Note that $m$ is uniquely determined up to a constant multiple.
The most canonical symmetric random walk is the simple random walk with
the transition probability given by
$p(e)=({\rm deg} \hspace{0.5mm} o(e))^{-1}$ ($e\in E$).

\begin{figure}[htbp]
\begin{center}
\scalebox{1.0}{\input{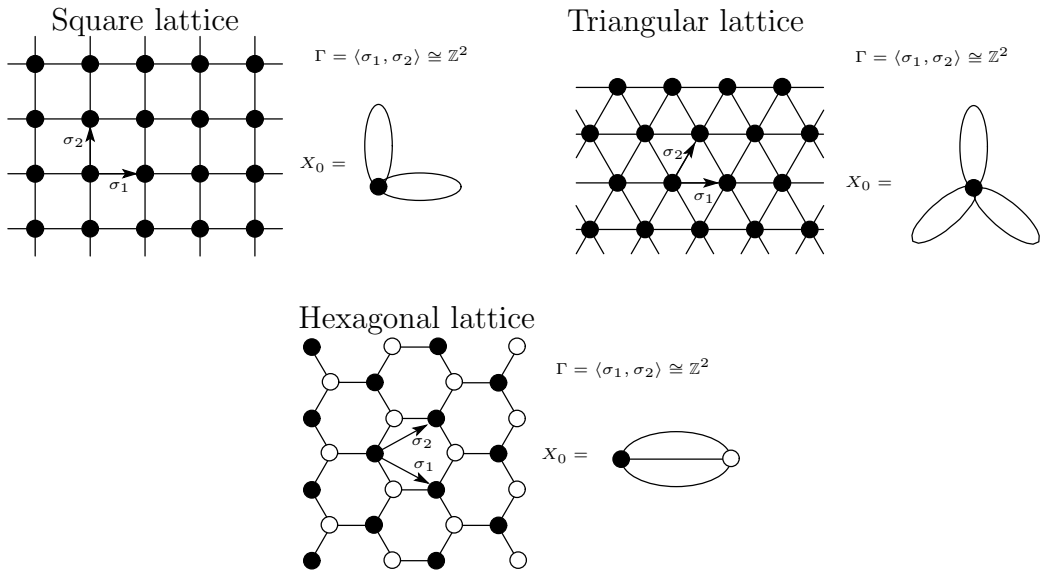}}
\end{center}
\caption{Crystal lattices}
\label{figure:example}
\end{figure}

Studying the long time asymptotics for random walks is 
a central theme in probability theory. In particular,
the central limit theorem (CLT), a generalization of the Laplace--de Moivre
theorem, has been studied by many authors in various settings. 
For basic results, see Spitzer \cite{Spitzer}, Woess \cite{Woess}, 
Lawler \cite{Lawler} and literatures therein.
As mentioned in Spitzer \cite{Spitzer}, 
the periodicity of the graph plays a crucial role 
to obtain such asymptotics. 
From this viewpoint, Kotani, Shirai and Sunada \cite{KSS}
applied {\it discrete geometric analysis} to study
the long time asymptotics of symmetric random walks 
on {\it crystal lattices}. 
We also refer to Guivar'ch \cite{Gui} and 
Kramli--Szasz \cite{Kramli} for related early works.
Here $X=(V,E)$ is called a ($\Gamma$-)crystal lattice if there exists an abelian group $\Gamma$ 
acting on $X$ freely and its quotient $X_0=(V_0, E_0)=\Gamma \backslash X$ is a finite graph. 
In other words, $X$ is an abelian covering graph of a finite graph $X_0$ for which
covering transformation group is $\Gamma$. Examples we have in mind are the square lattice 
$\mathbb Z^{d}$, the triangular lattice and the hexagonal lattice (see Figure \ref{figure:example}).
For simplicity, we assume $\Gamma$ is torsion free,
therefore isomorphic to $\mathbb{Z}^d$.
If both $p$ and $m$ are $\Gamma$-invariant, then the symmetric random walk
$\{ w_{n} \}_{n=0}^{\infty}$ induces a symmetric random walk on the quotient graph $X_0$ through
the covering map $\pi: X \to X_{0}$, and vice versa.

Later in \cite{KS00}, Kotani and Sunada studied
a numerical estimate of the $n$-step transition probability 
$p(n,x,y)$ for fixed $x,y\in V$ as $n\to \infty$,
called the local central limit theorem (LCLT), 
for the symmetric random walk on the crystal lattice $X$
by placing a special emphasis on the geometric feature.
Moreover, they also established 
the asymptotic expansion 
\begin{eqnarray}
p(n,x,y)m(y)^{-1} 
&\sim & a_{0}(X) n^{-d/2}
\exp \big( -\frac{d_\Gamma(x,y)^2}{2n} \big) 
\nonumber \\
&\mbox{  }& 
\times
 (1+a_1(x,y)n^{-1}+a_2(x,y) n^{-2}+ \cdots) 
 \quad \mbox{as  }n\to \infty,
 \label{AE-intro}
\end{eqnarray}
where $d_\Gamma(x,y):=\vert \Phi_{0}(y)-\Phi_{0}(x) \vert_{\Gamma \otimes \mathbb R}$ 
is a Euclidean pseudo-distance appearing through 
the {\it standard realization} $\Phi_0: X \rightarrow 
\Gamma \otimes \mathbb{R} \cong
\mathbb{R}^d$ (cf. \cite{KS00-2}).
The asymptotic expansion (\ref{AE-intro}) yields the CLT
\begin{eqnarray}
L^{[nt]}P_{n^{-1/2}} f(x_n) &:= &
\sum_{y \in V} p([nt], x_n, y) f \big(
n^{-1/2} \Phi_0 (y) \big)  \nonumber \\
& \rightarrow &
e^{-t\Delta/2} f(\mathbf{x})
\qquad (\mathbf{x} \in \Gamma \otimes \mathbb{R}) \quad \mbox{as  }n\to \infty,
\label{rough-CLT}
\end{eqnarray}
where $\{ x_n \}_{n=1}^\infty$ is a sequence in $V$ such that 
$\lim_{n\rightarrow \infty}n^{-1/2} \Phi_0(x_n)=\mathbf{x}$ and 
$f$ is a continuous function on $\Gamma \otimes \mathbb{R}$ 
vanishing at infinity.
Here $e^{-t\Delta}$ is the heat semigroup generated by the (positive) 
{\it Albanese Laplacian} $\Delta$ on $\Gamma \otimes \mathbb{R}$.
We mention that (\ref{rough-CLT}) is also obtained as a special case of 
\cite{Ishiwata, Kotani} by applying Trotter's approximation theory \cite{Trotter}.
In terms of probability theory,  (\ref{rough-CLT}) means that, as $n\to \infty$,
a sequence of $\Gamma \otimes \mathbb R$-valued random variables
$\{ n^{-1/2} \Phi_0(w_{[nt]}) \}_{n=1}^{\infty}$ $(t\geq 0)$ converges 
to $B_{t}$ as $n\to \infty$ in law, where $(B_{t})_{t\geq 0}$ is 
the standard Brownian motion on $\Gamma \otimes \mathbb{R}$ 
starting from ${\bf x}$.
Nevertheless, this kind of long time asymptotics
for random walks on crystal lattices which are not necessarily symmetric 
has not been studied 
satisfactorily although a large deviation principle (LDP)
is obtained in \cite{Kotani contemp, KS06}.

The main purpose of the present paper is to discuss the
long time asymptotics for non-symmetric random walks on
crystal lattices. In particular, 
we establish two kinds of (functional) CLTs (Theorems \ref{CLT-1}--\ref{FCLT-2})
and extend the asymptotic expansion formula (\ref{AE-intro})
to non-symmetric cases (Theorem \ref{AE-LCLT}). 
These CLTs are extensions of (\ref{rough-CLT}) to 
non-symmetric cases and Theorem \ref{AE-LCLT} is
a refinement of the LCLT for non-symmetric random walks on crystal lattices
presented by Sunada \cite{Sunada-Lecture}.

The rest of the present paper is organized as follows: In Section 2, 
we formulate our problem briefly and state the main results. 
In Section 3, we make a preparation from the discrete geometric analysis, 
some ergodic theorems and Trotter's approximation theory.
In Section 4, we prove the CLT of the first kind (Theorems \ref{CLT-1} and \ref{FCLT-1}). 
In Section 5, we prove the CLT of the second kind (Theorems \ref{CLT-2} and \ref{FCLT-2}). 
In Section 6, we concentrate on giving a spectral geometric
proof of Theorem \ref{AE-LCLT}. In the proof, perturbation arguments on
eigenvalues and eigenfunctions of the twisted transition operators $L_{\chi}$ and its
transposed operator $^t L_{\chi}$
play crucial roles, and we also make careful use of the classical
Laplace method to get the desired asymptotic expansion formula.
Finally, in Section \ref{section Example}, we present several concrete examples 
of the modified standard realization of crystal lattices associated 
with non-symmetric random walks. 

Throughout the present paper, $C$ denotes a positive constant that may change
at every occurrence, and $O(\cdot)$ stands for the Landau symbol. When 
the dependence of the $O(\cdot)$ term is significant, we denote it for example as 
$O_{N}(\cdot)$.
\section{Statement of the main results}
In this section, we state the main results. 
(For more details on the discrete geometric analysis, see Section 3.)
Let $p:E \to \mathbb R$ be a $\Gamma$-invariant 
non-negative function satisfying (\ref{p-def-intro}).
The random walk on the crystal lattice $X$ associated with
the transition probability $p$ induces a random walk on the quotient graph $X_{0}$
through the covering map $\pi: X \to X_{0}$, and 
the $n$-step transition probability $p(n,x,y)$ ($n\in \mathbb N, x,y\in V_{0}$) 
is also defined by (\ref{nstep-p-intro}).
Throughout the present paper, we assume 
that the random walk on $X_{0}$ is {\it{irreducible}}, that is, 
for every $x,y\in V_{0}$, there exists some $n=n(x,y)\in \mathbb N$ such that $p(n,x,y)>0$.
Then by applying Perron--Frobenius theorem (cf. Parry--Pollicott \cite[Theorem 2.2]{PP}),
we find a unique positive function $m:V_{0} \to \mathbb R$,
called the {\it{invariant probability measure}}, satisfying
\begin{equation}
m(x)=\sum_{e\in (E_{0})_{x}} p({\overline{e}}) m(t(e)) \quad (x\in V_{0}), 
\qquad
\sum_{x\in V_{0}} m(x)=1.
\label{m-PF-intro}
\end{equation}
We also write $m: V \to \mathbb R$ for  the ($\Gamma$-invariant)
lift of the invariant measure $m$.

For a topological space ${\cal T}$, let $C_\infty ({\cal T})$ 
denote the space of continuous functions on ${\cal T}$ vanishing at infinity. 
This space is endowed with usual uniform topology $\Vert \cdot \Vert_{\infty}$.
We note that a graph has the discrete topology defined by the graph distance.
Let $\mathrm{H}_1(X_0,\mathbb{R})$ and $\mathrm{H}^1(X_0, \mathbb{R})$ 
 be the first homology group and the first cohomology group on $X_0$, respectively.
Let $\rho_\mathbb{R} $ be the canonical 
surjective linear map
from $\mathrm{H}_1(X_0, \mathbb{R})$ to 
$\Gamma \otimes \mathbb{R}$. 
We define  
the {\it homological direction} $\gamma_{p}$
by
$$\gamma_{p}:=\sum_{e\in E_{0}} p(e)m(o(e))e
\in \mathrm{H}_1(X_0, \mathbb{R}),$$
and call $\rho_{\mathbb R}(\gamma_{p})(\in \Gamma \otimes \mathbb R)$ 
the {\it{asymptotic direction}}.
It should be noted that $\gamma_p=0$ if and only if the random walk on $X_{0}$
is symmetric, and 
$\rho_{\mathbb R}(\gamma_{p})={\bf{0}}$ does not
always imply the symmetry of the random walk on $X$.
We also introduce the {\it{modified harmonic realization}}
$\Phi_0 :X\rightarrow \Gamma\otimes \mathbb{R}$ by
\begin{equation}
\sum_{e\in E_x}p(e) \big (
\Phi_0(t(e))-\Phi_0(o(e)) \big)=
\rho_\mathbb{R}(\gamma_p) \qquad (x\in V).
\label{modified harmonic realization}
\end{equation} 
Note that $\Phi_0$ is uniquely determined 
up to translation (cf. \cite[page 854]{KS06}).
%

To present the CLT of the first kind, we need to introduce
the {\it transition-shift operator} $\mathcal{L}_{\gamma_p}$ acting on 
$C_\infty (X\times \mathrm{H}_1(X_0, \mathbb{R}))$ by
\begin{equation*}
\mathcal{L}_{\gamma_p} f( x, {\mathbf{z}}):= \sum_{e \in E_x} p(e)f(t(e), \mathbf{z}+\gamma_p)
\qquad (x\in V,  \hspace{1mm} {\bf{z}}\in {\mathrm{H}}_1(X_0, \mathbb{R})),
\end{equation*}
and the {\it{approximation operator}} 
${\cal P}_\varepsilon :C_\infty (\Gamma\otimes \mathbb{R}) \rightarrow
C_\infty (X \times \mathrm{H}_1(X_0, \mathbb{R}))$ ($0 \leq \varepsilon \leq 1$) by
\begin{equation*}
{\cal P}_\varepsilon f(x,{\bf{z}}):= f(\varepsilon(\Phi_0(x)-\rho_\mathbb{R}({\bf{z}}))).
\end{equation*}

Then the CLT of the first kind is stated as follows:
\begin{tm}\label{CLT-1}
For $0\leq s <t$ and $f \in C_\infty (\Gamma \otimes \mathbb{R})$,
\begin{equation}
\lim_{ n\rightarrow \infty} \left\| \mathcal{L}_{\gamma_p}^{[nt]-[ns]} {\cal P}_{n^{-1/2}} f 
- {\cal P}_{n^{-1/2}} 
e^{ -\frac{t-s}{2}\Delta} f \right\|_\infty =0,
\label{Trotter-1}
\end{equation}
where $\Delta$ is the {\rm{(}}positive{\rm{)}} Laplacian associated with the Albanese metric $g_{0}$ on 
$\Gamma \otimes \mathbb{R}$. {\rm{(}}See
Section {\rm{3}} for the definition of the Albanese metric $g_{0}$.{\rm{)}}

In particular, for any sequence $\{ (x_n, \mathbf{z}_n) \}_{n=1}^{\infty}$ in 
$V \times \mathrm{H}_1(X_0,\mathbb{R})$ with 
\begin{equation*}
\lim_{n\rightarrow \infty } 
n^{-1/2} \big( \Phi_{0}(x_n)-\rho_\mathbb{R}(\mathbf{z}_n) \big) =\mathbf{x}
 \in \Gamma \otimes \mathbb{R},
\end{equation*}
and for any $f \in C_\infty (\Gamma \otimes \mathbb{R})$, we have
\begin{equation}
\lim_{n\rightarrow \infty }\mathcal{L}_{\gamma_p}^{[nt]}{\cal P}_{n^{-1/2}} f(x_n,\mathbf{z}_n)
=e^{-t\Delta/2} f({\mathbf{x}}) := \int_{\Gamma \otimes \mathbb R}
G_{t}(\mathbf{x}-\mathbf{y})
f({\bf{y}}) d {\bf{y}} \qquad (t>0),
\label{Trotter-1-2}
\end{equation}
where 
$$ G_{t}({\mathbf{x}})=\frac{1}{(2\pi t)^{d/2}}
\exp \big( -\frac{ \vert {\mathbf{x}} \vert_{g_0}^{2} }{2t} \big) \qquad 
({\bf{x}}\in \Gamma \otimes \mathbb R)
$$
is the fundamental solution of the heat equation
\begin{equation*}
\frac{\partial}{\partial t} u(t, {\bf{x}})=-\frac{1}{2}\Delta u(t, {\bf{x}}).
\end{equation*}
\end{tm}
We mention that this theorem still holds for the approximation operator ${\cal{P}}_{\varepsilon}$
given by a general periodic realization 
$\Phi: X \rightarrow \Gamma \otimes \mathbb{R}$ (cf. Ishiwata \cite{Ishiwata}).

Now we set a reference point $x_{*} \in V$ such that $\Phi_{0}(x_{*})={\bf 0}$, and 
put
$$ \xi_{n}(c):=\Phi_0 \big (w_{n}(c) \big) \qquad (n=0,1,2,\ldots,~c\in \Omega_{x_{*}}(X)).$$
We then obtain a $\Gamma \otimes \mathbb R$-valued 
random walk $({\Omega}_{x_{*}}(X),
{\mathbb P}_{x_{*}}, \{\xi_{n} \}_{n=0}^{\infty})$. 
Let ${\cal X}_t^{(n)}$ ($t\geq 0, n\in \mathbb N$) be a map from ${\Omega}_{x_{*}}(X)$ to 
$\Gamma \otimes \mathbb{R}$ 
given by
\begin{equation*}
{\cal  X}_t^{(n)}(c) := \frac{1}{\sqrt{n}}  
\big ( \xi_{[nt]}(c)- [nt] \rho_\mathbb{R} (\gamma_p ) \big ) \qquad 
(c\in {\Omega}_{x_{*}}(X)).
\end{equation*}
Then (\ref{Trotter-1-2}) can be rewritten as
\begin{equation*}
\lim_{n\rightarrow \infty} 
\sum_{c \in {\Omega}_{x_{*}} (X)} f({\cal X}_t^{(n)}(c)) \mathbb{P}_{x_{*}} (\{ c\})
=\int_{\bf W}
f({\bf w}_t) {\bf{P}}^{W}(d{\bf w}),
\end{equation*}
where ${\bf{P}}^{W}$ is the {\it{Wiener measure}} on 
${\bf W}:=C_{\bf{0}}([0,\infty), \Gamma \otimes \mathbb R)$.
As the piecewise linear interpolation of $( {\cal X}_t^{(n)})_{t\geq 0}$, 
we define a map
${\bf X}^{(n)}: {\Omega}_{x_{*}}(X) \rightarrow {\bf W}$ by
\begin{equation*}
{\bf  X}_{t}^{(n)} ( c) :=
 \frac{1}{\sqrt{n}} \Big \{ \xi_{[nt]}(c) 
+ (nt -[nt]) \big( \xi_{ [nt]+1}(c) 
-\xi_{[nt]}(c) \big)
-nt
\rho_\mathbb{R}  (\gamma_p) \Big \}
\qquad (t\geq 0).
\end{equation*}
Let ${\bf P}^{(n)} $ be the probability measure on $({\bf{W}}, {\cal B}({\bf W}))$ induced by ${\bf  X}^{(n)}$.

The following theorem is the functional CLT (i.e., Donsker's invariance principle)
for the random walk $\{ \xi_{n} \}_{n=0}^{\infty}$.
\begin{tm} \label{FCLT-1}
The sequence $ \{ {\bf P}^{(n)} \}_{n=1}^{\infty} $
converges weakly to the Wiener measure ${\bf P}^{W}$
as $n\to \infty$. Namely, ${\bf X}^{(n)}$
converges to a $\Gamma \otimes \mathbb R$-valued standard Brownian motion 
$(B_{t})_{t\geq 0}$ with $B_{0}={\bf 0}$
in law.
\end{tm}

We next present another kind of CLT for a family of random walks.
(See Durret \cite{Durrett} and Trotter \cite{Trotter} for related results.)
For the transition probability $p$, we introduce a family of ($\Gamma$-invariant)
transition probabilities $\{ p_{\varepsilon} \}_{0\leq \varepsilon \leq 1}$ on $X$ 
by
\begin{equation}
p_\varepsilon(e):=p_0(e)+\varepsilon q(e) \qquad (e\in E),
\label{def-pe}
\end{equation}
where
\begin{equation*}
p_0(e):=\frac{1}{2} \Big( p(e)+\frac{m(t(e))}{m(o(e))}p(\overline{e}) \Big), ~~
q(e):=\frac{1}{2}\Big( p(e)-\frac{m(t(e))}{m(o(e))}p(\overline{e}) \Big).
\end{equation*} 
We note that $\{ p_\varepsilon \}_{0\leq \varepsilon \leq 1}$ is 
the interpolation between the original transition probability 
$p=p_1$ and the $m$-symmetric probability $p_0$ 
such that the homological direction
$\gamma_{p_\varepsilon}$ equals $\varepsilon \gamma_p$ for every 
$0\leq \varepsilon \leq 1$ (see Lemma \ref{pdelta} below for details).

We denote by $L_{(\varepsilon)}$ 
the transition operator associated with $p_\varepsilon$.
We denote by $g_{0}^{(\varepsilon)}$ the corresponding
Albanese metric on $\Gamma \otimes \mathbb R$.
If we need to emphasize the flat
metric $g_{0}^{(\varepsilon)}$, we write $(\Gamma \otimes \mathbb R)_{(\varepsilon)}$
for $(\Gamma \otimes \mathbb R, g_{0}^{(\varepsilon)})$.
We define another approximation operator 
$P_{\varepsilon} : C_\infty ((\Gamma\otimes \mathbb{R})_{(0)}) \rightarrow C_\infty (X)$
by 
\begin{equation*}
P_{\varepsilon}f(x):=f(\varepsilon \Phi_0^{(\varepsilon)}(x)) \qquad (x\in V),
\end{equation*}
where $\Phi_0^{(\varepsilon)}:X \to \Gamma \otimes \mathbb R$ 
is the modified harmonic realization 
associated with $p_{\varepsilon}$. 

Then the CLT of the second kind is stated as follows:
\begin{tm}\label{CLT-2}
 For $0\leq s <t$ and $f \in C_\infty ((\Gamma \otimes \mathbb{R})_{(0)})$, 
\begin{equation*}
\lim_{ n \rightarrow \infty} 
\Big \|  L_{(n^{-1/2})}^{[nt]-[ns]} P_{n^{-1/2}} f  - P_{n^{-1/2}} 
e^{-(t-s) (\frac{1}{2}\Delta_{(0)}-\langle \rho_{\mathbb R}(\gamma_p ), \nabla_{(0)} 
\rangle_{g_{0}^{(0)}} )} f  \Big \|_\infty =0,
\end{equation*}
where $\Delta_{(0)}$ and $\nabla_{(0)}$ stand for the {\rm{(}}positive{\rm{)}} Laplacian 
and the gradient on $(\Gamma \otimes \mathbb R)_{(0)}$, respectively.
In particular, for any sequence $\{x_n \}_{n=1}^\infty$ in $V$ with 
$$
\lim_{n\to \infty}
n^{-1/2}\Phi^{(n^{-1/2})}_{0}(x_n)=\mathbf{x} \in (\Gamma\otimes \mathbb{R})_{(0)}, $$
and for any $f\in C_{\infty}((\Gamma \otimes \mathbb R)_{(0)})$, we have
\begin{eqnarray}
\lefteqn{
\lim_{n\rightarrow \infty}
\sum_{y \in V}p_{n^{-1/2}}([nt], x_n, y)f(n^{-1/2}\Phi_{0}^{(n^{-1/2})}(y) )
}
\nonumber \\
&=& e^{-t(\frac{1}{2}\Delta_{(0)}-\langle \rho_{\mathbb R}(\gamma_p ), \nabla_{(0)} 
\rangle_{g_{0}^{(0)}} )} f (\mathbf{x}) 
=\int_{(\Gamma\otimes \mathbb{R})_{(0)}} H_t(\mathbf{x}-\mathbf{y}) f(\mathbf{y}) d \mathbf{y},
\label{Trotter-2-2}
\end{eqnarray}
where
\begin{equation*}
H_t(\mathbf{x})=\frac{1}{(2\pi t)^{d/2}}\exp \Big(-\frac{|\mathbf{x}-\rho_\mathbb{R}
(\gamma_p)t|_{g_0^{(0)}}^2}{2t} \Big)
\end{equation*}
is the fundamental solution of the heat equation
\begin{equation*}
\frac{\partial}{\partial t} u(t, {\bf{x}})=-\frac{1}{2}\Delta_{(0)} u(t, {\bf{x}})+
\langle \rho_\mathbb{R}(\gamma_p), \nabla_{(0)} u(t, {\bf{x}})\rangle_{g_0^{(0)}}.
\end{equation*}
\end{tm}

Just like Theorems  \ref{CLT-1} and  \ref{FCLT-1}, this theorem also
implies a functional CLT.
Set a reference point $x_{*} \in V$ such that $\Phi_{0}^{(\varepsilon)}(x_{*}) = {\bf 0}$
for all $0\leq \varepsilon \leq 1$.
Putting 
$$ \xi^{(\varepsilon)}_{n}(c):=\Phi^{(\varepsilon)}_{0} \big (w_{n}(c) \big) \qquad 
(0\leq \varepsilon \leq 1,~ n=0,1,2,\ldots, ~c\in {\Omega}_{x_{*}}(X) ),$$
we also obtain a $\Gamma \otimes \mathbb R$-valued random walk
$({\Omega}_{x_{*}}(X), {\mathbb P}_{x_{*}}^{(\varepsilon)} , 
\{ \xi_{n}^{(\varepsilon)} \}_{n=0}^{\infty})$ 
associated with $p_\varepsilon$.
Let 
${\cal Y}_t^{(\varepsilon, n)}$ ($t\geq 0$, $n\in \mathbb N$, $0\leq \varepsilon \leq 1$)
be a map from
$ {\Omega}_{x_{*}} (X) $ to $\Gamma\otimes \mathbb{R}$ given
by
\begin{equation*}
{{\cal {Y}}}_{t}^{(\varepsilon, n)}(c) := \frac{1}{\sqrt{n}}  \xi_{[nt]}^{(\varepsilon)} (c)  \qquad 
(c\in {\Omega}_{x_{*}}(X) ).
\end{equation*}
Then (\ref{Trotter-2-2}) can be rewritten as
\begin{equation*}
\lim_{n\rightarrow \infty} 
\sum_{c\in {\Omega}_{x_{*}}(X)} f \big ({\cal{Y}}_t^{(n^{-1/2}, n)}(c) \big) \mathbb{P}_{x_{*}} (\{c\})
=\int_{{\bf W}_{(0)}} f({\bf w}_t) {\bf{Q}} (d{\bf w}),
\end{equation*}
where $\bf {Q}$ is the probability measure on ${\bf {W}}_{(0)}:=
C_{\bf{0}}([0,\infty), (\Gamma \otimes \mathbb{R})_{(0)})$ induced by 
$(B_t+\rho_{\mathbb{R}}(\gamma_p)t )_{t\geq 0}$.
Here $(B_t)_{t\geq 0}$ is a $(\Gamma \otimes \mathbb{R})_{(0)}$-valued standard 
Brownian motion with $B_0={\bf{0}}$.
We define a map
${\bf Y}^{(\varepsilon, n)}: {\Omega}_{x_{*}} (X) \rightarrow {\bf W}_{(0)}$ by
\begin{equation*}
{\bf  Y}_{t}^{(\varepsilon, n)} ( c):=
\frac{1}{\sqrt{n}} \Big \{ \xi_{[nt]}^{(\varepsilon)} (c) + (nt -[nt])
\big (\xi_{ [nt]+1}^{(\varepsilon)} (c) 
-\xi_{[nt]}^{(\varepsilon)} (c)  \big)
 \Big \} \qquad (t\geq 0,~c\in {\Omega}_{x_{*}}(X) ),
\end{equation*}
which is the piecewise linear interpolation of $({\cal Y}^{(\varepsilon, n)})_{t\geq 0}$.
Let ${\bf  Q}^{(\varepsilon, n)} $ be the probability measure on $({\bf W}_{(0)}, {\cal B}({\bf W}_{(0)}))$
induced by ${\bf Y}^{(\varepsilon, n)}$.

The following theorem is a functional CLT for the family of random walks
$\{ \xi_{n}^{(\varepsilon)} \}_{n=0}^{\infty}$ ($0\leq \varepsilon \leq 1$).
\begin{tm} \label{FCLT-2}
The sequence $ \{ {\bf Q}^{(n^{-1/2}, n)} \}_{n=1}^{\infty} $
converges weakly to $\bf{Q}$ as $n\to \infty$.
Namely, $\big ({\bf Y}^{(n^{-1/2}, n)}_{t} \big )_{t\geq 0}$
converges to 
a $(\Gamma \otimes \mathbb R)_{(0)}$-valued standard Brownian motion with
drift $\rho_{\mathbb R}(\gamma_{p})$ starting from $\bf 0$  
as $n\to \infty$ in law.
\end{tm}

Next, we discuss the precise asymptotic behavior of the $n$-step transition probability
$p(n,x,y)$ ($x,y\in V$) of the non-symmetric random walk $\{w_{n} \}_{n=0}^{\infty}$ on $X$ 
as $n \to \infty$. Here we impose the condition that the random walk $\{w_{n} \}_{n=0}^{\infty}$ 
on $X$ is {\it{irreducible}}.
Note that this condition implies the {\it irreducibility on $X_{0}$} 
imposed above.
Conversely, the irreducibility of the random walk on $X_0$ 
does not imply the irreducibility of the random walk on $X$. 
We define the {\it{period}} of the random walk $\{w_{n} \}_{n=0}^{\infty}$ by
$K:={\rm{gcd}} \{ n\in \mathbb N \vert~p(n,x,x)>0 \}$. It should be noted that $K$ is 
independent of $x\in V$ by the irreducibility condition. Sunada \cite[page 121]{Sunada-Lecture}
presented the LCLT
\begin{eqnarray}
\lefteqn{(2\pi n)^{d/2}p(n,x,y)m(y)^{-1} }
\nonumber \\
& &
\sim
K {\rm{vol}}({\rm{Alb}}^{\Gamma}) 
\exp \left( 
-\frac{ \vert \Phi_{0}(y)-\Phi_{0}(x) -n\rho_\mathbb{R}(\gamma_p) \vert_{g_{0}}^2}{2n} 
\right) \qquad (x,y\in V) 
\label{Sunada-LCLT}
\end{eqnarray}
as $n \to \infty$, where ${\rm{vol}}({\rm{Alb}}^{\Gamma})$ stands for the volume of
the $\Gamma$-Albanese torus $(\Gamma \otimes \mathbb R / \Gamma \otimes \mathbb Z, g_{0})$.

As a refinement of (\ref{Sunada-LCLT}), 
we obtain the following 
precise asymptotics of $p(n,x,y)$.
\begin{tm} \label{AE-LCLT}
Suppose that the random walk $\{ w_{n} \}_{n=0}^{\infty}$ on $X$ is irreducible with period $K$. 
Let $V=\coprod_{k=0}^{K-1} A_k$ be the corresponding $K$-partition of $V$. 
Then for any $x\in A_i$ and $y \in A_j$ {\rm{(}}$0\leq i,j \leq K-1${\rm{)}}, 
we have 
\begin{equation*}
p(n,x,y)= 0  \qquad (n\neq Kl+j-i)
\end{equation*}
and 
\begin{align}
(2\pi n)^{d/2} p(n ,x,y)m(y)^{-1}
&=
K {\rm vol}({\rm{Alb}}^{\Gamma})
\exp \Big( -\frac{\big \vert \Phi_{0}(y)-\Phi_{0}(x) 
-n{\rho}_{\mathbb R}(\gamma_{p}) \big \vert^{2}_{g_{0}}}
{2n} \Big) 
\nonumber \\
&\mbox{ }
\hspace{-30mm}
\times 
\Big (1+a_{1} \big (\pi(x),\pi(y),\gamma_{p}; \Phi_{0}(y)-\Phi_{0}(x) 
-n{\rho}_{\mathbb R}(\gamma_{p}) 
\big)n^{-1}
\Big) 
+
O(n^{-3/2})
\label{1-AE}
\end{align}
as $n=Kl+j-i \rightarrow \infty$ uniformly for $x,y$ with
$\vert \Phi_{0}(y)-\Phi_{0}(x)-n\rho_{\mathbb R}(\gamma_{p})
\vert_{g_{0}} \leq Cn^{1/6}$.
\end{tm}

The coefficient $a_{1}=a_{1}\big (\pi(x),\pi(y), \gamma_{p}; \Phi_{0}(y)-\Phi_{0}(x) 
-n{\rho}_{\mathbb R}(\gamma_{p}) 
\big)$ is nothing but
the one of the term $n^{-1}$
in the power series expansion of
$$ U_{n}:=\frac{  
(2\pi n)^{d/2} p(n ,x,y)m(y)^{-1}}{
K {\rm vol}({\rm{Alb}}^{\Gamma})
\exp \big( -\big \vert \Phi_{0}(y)-\Phi_{0}(x) 
-n{\rho}_{\mathbb R}(\gamma_{p}) \big \vert^{2}_{g_{0}}/2n \big) }
$$ in $n^{-1}$, and its explicit expression will be given in Theorem 
\ref{a_1} below. 
Using the same argument as in Section 6, 
we can also give the coefficient $a_{j}$ of the term $n^{-j}$ for any $j\geq 2$ in the
power series expansion of $U_{n}$. We do not write the explicit asymptotic
expansion of $p(n,x,y)$ in the present paper, because
of its complication.

We should also mention that Uchiyama and his coauthor \cite{Kazami-Uchiyama, Uchiyama} recently
obtained this kind of asymptotic expansion formula for non-symmetric random walks
on periodic graphs (i.e., crystal lattices) realized in Euclidean space under the 
{\it{zero mean condition}}. Their approach is purely probabilistic, and roughly speaking, 
the zero mean condition is translated into
$\rho_{\mathbb R}(\gamma_{p})={\bf 0}$ in our geometric framework.
Thus Theorem \ref{AE-LCLT} can be regarded as a generalization of 
their asymptotic expansion including the case $\rho_{\mathbb R}(\gamma_{p})\neq {\bf 0}$.
\section{Preliminaries}
\label{section2}
\subsection{Discrete geometric analysis on graphs}
In this subsection, we review some basic facts of {\it{discrete geometric analysis}} 
on graphs. We refer to Kotani--Sunada \cite{KS-contemp}, Sunada \cite{Sunada-DGA, S} and 
references therein for details.

First of all, we consider an irreducible random walk on 
$X_{0}=(V_{0}, E_{0})$ with the transition probability $p: E_{0} \to [0,1]$.
By the Perron--Frobenius theorem, 
there exists a unique positive
function $m: V_{0}\to \mathbb R$, called the invariant probability measure, satisfying
(\ref{m-PF-intro}).
We put $$ {\widetilde{m}}(e)=m(o(e))p(e) \qquad   (e\in E_{0}).$$
It is easy to see that ${\widetilde{m}}:E_{0}\to \mathbb R$ satisfies 
$\sum_{e\in E_{0}} {\widetilde{m}}(e)=1$
and
\begin{equation}
m(x)=\sum_{e\in (E_{0})_{x}} {\widetilde{m}}(e)
=\sum_{e\in (E_{0})_{x}} {\widetilde{m}}({\overline{e}}) \qquad (x\in V_{0}).
\label{m-inv}
\end{equation}
When $ {\widetilde{m}}(e)={\widetilde{m}}({\overline{e}})$, the random walk is said to be 
($m$-){\it{symmetric}}.

We consider the 0-chain group 
$$C_{0}(X_{0},\mathbb R)=\big \{ \sum_{x\in V_{0} }a_{x}x \vert~a_{x}\in \mathbb R  \big \}$$
and the 1-chain group 
$$C_{1}(X_{0},\mathbb R)=\big \{ \sum_{e\in E_{0} }a_{e}e\vert~a_{e}\in \mathbb R  \big \},$$
where the relation ${\overline{e}}=-e$ is imposed for $e\in E_{0}$.
The boundary operator ${\partial}:C_{1}(X_{0},\mathbb R) \to C_{0}(X_{0},\mathbb R)$ is defined 
by $\partial (e):=t(e)-o(e)$. The first homology group ${\rm H}_{1}(X_{0},{\mathbb R})$ 
is the kernel ${\rm Ker}(\partial)
\subset C_{1}(X_{0},\mathbb R)$. ${\rm H}_{1}(X_{0},{\mathbb Z})$ 
is also defined by replacing $\mathbb R$ by ${\mathbb Z}$. 

We define the 0-cochain group $$C^{0}(X_{0},\mathbb R):=\{ f:V_{0} \to \mathbb R \}$$ 
with the inner product 
$$ \langle f_{1}, f_{2} \rangle_{0} =\sum_{x\in V_{0}} f_{1}(x) f_{2}(x)
 \qquad (f_{1},f_{2}\in C^{0}(X_{0},\mathbb R) ),$$
and the 
1-cochain group $$C^{1}(X_{0},\mathbb R):=\{ \omega :E_{0} \to \mathbb R  \vert~\omega({\overline{e}})= -\omega(e) \}$$
with the inner product
$$ \langle  \omega_{1}, \omega_{2} \rangle_{1} = \frac{1}{2}
\sum_{e\in E_{0}} \omega_{1}(e) \omega_{2}(e)
 \qquad ( \omega_{1}, \omega_{2}\in C^{1}(X_{0},\mathbb R) ).$$
A $1$-cochain is occasionally called a $1$-form on $X_{0}$.

We define the difference operator $d: C^{0}(X_{0},\mathbb R) \to C^{1}(X_{0},\mathbb R)$
by $$df(e):=f(t(e))-f(o(e)) \qquad (e \in E_{0} ),$$
and the first cohomology group 
${\rm H}^{1}(X_{0},{\mathbb R}):=C^{1}(X_{0},\mathbb R) / {\rm Im}(d)$. Note that 
${\rm H}^{1}(X_{0},{\mathbb R})$
is the dual of 
the first homology group 
${\rm H}_{1}(X_{0},{\mathbb R})$.
We also define the transition operator $L:C^{0}(X_{0},\mathbb R) \to C^{0}(X_{0},\mathbb R)$ by
$$ Lf(x):=(I-\delta_{p} d)f(x)=\sum_{e\in (E_{0})_{x}} p(e) f(t(e)) \qquad (x\in V_{0} ),$$
where $\delta_{p}: C^{1}(X_{0},\mathbb R) \to C^{0}(X_{0},\mathbb R)$ is given by 
$$ (\delta_{p} \omega)(x):=-\sum_{e\in (E_{0})_{x}} p(e)\omega(e) \qquad (x\in V_{0}).$$
We easily see 
$$p(n,x,y)=L^{n}\delta_{y}(x) \qquad  (n\in \mathbb N, x,y\in V_{0}), $$
where $\delta_{y}(\cdot)$ is the Dirac delta function with pole
at $y$. 

Let $\hspace{-1mm}\mbox{ }^{t}L:C^{0}(X_{0},\mathbb R) \to C^{0}(X_{0},\mathbb R)$ be the transposed
operator of the transition
operator $L$ with respect to the inner product $\langle \cdot, \cdot \rangle_{0}$. 
The explicit form is given by
$$ \hspace{-1mm}\mbox{ }^{t}Lf(x):=\sum_{e\in (E_{0})_{x}} p({\overline{e}}) f(t(e)) \qquad (x\in V_{0}).$$
We note that (\ref{m-PF-intro}) is written as 
$$\hspace{-1mm}\mbox{ }^{t}Lm(x)=m(x), \quad \sum_{x\in V_{0}} m(x)=1.$$
It means that the invariant measure $m$ is an eigenfunction of $\hspace{-1mm}\mbox{ }^{t}L$ for the
maximal positive eigenvalue $1$.

Now we provide the notion of {\it{modified harmonic}} 1-{\it{forms}} on $X_{0}$ associated with 
the transition probability $p$. 
We introduce the random variables $\eta_{i}: {\Omega}_{x}(X_{0}) \to C_{1}(X_{0},\mathbb R)$~
($i \in \mathbb N, x\in V_{0}$) defined by
$\eta_{i}(c):=e_{i}$ for $c=(e_{1},e_{2},\ldots) \in {\Omega}_{x}(X_{0})$, and set
\begin{equation}
\gamma_{p}:=\sum_{e\in E_{0}} {\widetilde{m}}(e) e \in C_{1}(X_{0}, \mathbb R).
\label{homological direction}
\nonumber
\end{equation}
A simple application of the ergodic theorem leads to
\begin{equation}
\lim_{n\to \infty} \frac{1}{n} \sum_{i=1}^{n} \eta_{i}(c)=\gamma_{p}, \quad {\mathbb P}_{x}{\mbox{-a.e.}}~c\in {\Omega}_{x}(X_{0}).
\label{hd-ergord}
\nonumber
\end{equation}
By virtue of (\ref{m-inv}), we observe that
$\partial \gamma_{p}=0$
and hence
$\gamma_{p}\in {\rm H}_{1}(X_{0},\mathbb R)$. It provides a quantity to
measure a homological drift of the random walk and we call it the {\it{homological direction}} of 
the given random walk. Straightforwardly,  $\gamma_{p}=0$ if and only if $p$ gives a 
symmetric random walk, i.e., $ {\widetilde{m}}(e)={\widetilde{m}}({\overline{e}})$.
A $1$-form $\omega$ is said to be {\it{modified harmonic}} if
\begin{equation}
\delta_{p} \omega(x)+ \langle \gamma_{p},\omega \rangle=0 \qquad (x\in V_{0} ),
\label{M-harmonic}
\end{equation}
where it should be noted that $\langle \gamma_{p},\omega \rangle:=
\subscripts
    {C_{1}(X_{0},\mathbb R)}
    {\left\langle \gamma_{p},\omega \right\rangle}
    {C^{1}(X_{0},\mathbb R)}$ is constant as a function on $V_{0}$.
We denote by ${\cal H}^{1}(X_{0}) $ the set of modified harmonic $1$-forms, and equip ${\cal H}^{1}(X_{0}) $
with the inner product 
\begin{equation}
\langle \hspace{-0.8mm} \langle \omega_{1}, \omega_{2} \rangle \hspace{-0.8mm} \rangle 
:= 
\sum_{e\in E_{0}}
\omega_{1}(e) \omega_{2}(e) {\widetilde{m}}(e)
-\langle \gamma_{p}, \omega_{1} \rangle
\langle \gamma_{p}, \omega_{2} \rangle
\qquad ( \omega_{1}, \omega_{2} \in {\cal H}^{1}(X_{0}) ).
\label{flatmetric-norm}
\end{equation}
Then the corresponding norm $\Vert \cdot \Vert$ is given by
\begin{equation} 
\Vert \omega \Vert^{2}:=
\langle \hspace{-0.8mm} \langle \omega, \omega \rangle \hspace{-0.8mm} \rangle
=\sum_{e\in E_{0}} \vert \omega(e) \vert^{2}{\widetilde{m}}(e)
-\langle \gamma_{p}, \omega \rangle^{2} \qquad (\omega \in
{\cal H}^{1}(X_{0}) ).
\nonumber
\end{equation}
By the discrete Hodge--Kodaira theorem   
(cf. \cite[Lemma 5.2]{KS06}), we may identify ${\rm H}^{1}(X_{0},\mathbb R)$ and 
${\rm H}^{1}(X_{0},\mathbb Z)$ with
${\cal H}^{1}(X_{0})$ and 
$$ \Big \{ \omega \in {\cal H}^{1}(X_{0})
\big \vert~\int_{c} \omega :=\sum_{i=1}^{n} \omega (e_{i})
\in \mathbb Z \mbox {  for every closed path } c =(e_{1}, \ldots, e_{n}) \mbox{ in }X_{0} \Big \},$$
respectively.
Using this identification, we obtain an inner product $\langle \hspace{-0.8mm} \langle
\cdot, \cdot \rangle \hspace{-0.8mm} \rangle$ on ${\rm H}^{1}(X_{0},\mathbb R)$.
\begin{re} \label{metric}
In {\rm{\cite{KS00}}}, 
Kotani and Sunada considered 
a symmetric random  walk {\rm{(}}i.e., $\gamma_{p}=0${\rm{)}} and defined
the Albanese metric on
$\Gamma \otimes \mathbb R$ associated with the inner product 
$$ 
\frac{1}{2} \sum_{e\in E_{0}} \omega_{1}(e)
\omega_{2}(e) {\widetilde m}(e) \qquad (\omega_{1}, \omega_{2} \in  {\cal H}^{1}(X_{0}) ).$$
Due to this difference of the definition of the 
Albanese metric,
the Laplacian $\Delta$ on $\Gamma \otimes \mathbb R$ appears
in long time asymptotics 
instead of $\Delta/2$.
\end{re}
 
Let $X=(V,E)$ be a ($\Gamma$-)crystal lattice, that is,
$X$ is an abelian covering graph of a finite graph $X_{0}$ in which covering transformation group is 
$\Gamma$.
Let $p: E \to \mathbb R$ and $m:V \to \mathbb R$ be the lift of
the transition probability
$p: E_{0} \to \mathbb R$ and the invariant measure $m:V_{0} \to \mathbb R$, respectively. Namely,
$$ p(\sigma e)=p(e), \quad m(\sigma x)=m(x) \qquad (e\in E, x\in V, \sigma\in \Gamma). $$

We denote by $\pi: X\to X_{0}$ the covering map, and by
$\rho: {\rm H}_{1}(X_{0},\mathbb Z) \to \Gamma$ the surjective homomorphism associated 
with the covering map
$\pi$. We extend $\rho$ to the surjective linear map
$ {\rho}_{\mathbb R}: {\rm H}_{1}(X_{0},\mathbb R) \to \Gamma \otimes \mathbb R$.
Then we may consider 
the injective linear map
$\hspace{-1mm}\mbox{ }^{t}\rho_{\mathbb R}: {\rm Hom}(\Gamma, \mathbb R) \to 
{\rm H}^{1}(X_{0},\mathbb R)$ by
$$ \hspace{-1mm}\mbox{ }^{t}\rho_{\mathbb R}: {\omega} \in {\rm Hom}(\Gamma, \mathbb R)
\mapsto \hspace{-1mm}\mbox{ }^{t}\rho_{\mathbb R}(\omega)(\cdot):=\omega(\rho_{\mathbb R}(\cdot)) 
\in {\rm H}^{1}(X_{0},\mathbb R), $$
where ${\rm Hom}(\Gamma, \mathbb R)$ denotes the linear space
of homomorphisms of $\Gamma$ into $\mathbb R$.
Using the maps $\hspace{-1mm}\mbox{ }^{t}\rho_{\mathbb R}$ and $\rho_{\mathbb R}$, 
we identify ${\rm Hom}(\Gamma, \mathbb R)$
with the subspace ${\rm Image}(\hspace{-1mm}\mbox{ }^{t}\rho_{\mathbb R})$ in 
${\rm H}^{1}(X_{0},\mathbb R)$ and $\Gamma \otimes \mathbb R$ with the quotient linear subspace
of ${\rm H}_{1}(X_{0},\mathbb R)$. 
Throughout the present paper, we shall denote 
$\hspace{-1mm}\mbox{ }^{t}\rho_{\mathbb R}(\omega) \in {\rm H}^{1}(X_{0}, {\mathbb R})$ 
by the same symbol $\omega$ for brevity.
We restrict the inner product
$\langle \hspace{-0.8mm} \langle
\cdot, \cdot \rangle \hspace{-0.8mm} \rangle$ on ${\rm H}^{1}(X_{0},\mathbb R)$
to the subspace ${\rm Hom}(\Gamma, \mathbb R)$, and 
then take up
the dual inner product $\langle \cdot, \cdot \rangle_{alb}$
on $\Gamma \otimes \mathbb R$. The flat metric on $\Gamma \otimes \mathbb R$
induced from this inner product is called the {\it{Albanese metric}} and 
is denoted by $g_{0}$. 
This procedure is summarized in the following diagram:
\begin{align*}
&(\Gamma \otimes  \mathbb{R}, g_{0})  & &\stackrel{\rho_\mathbb{R}}{\longleftarrow \!\!\longleftarrow }
& &{\rm H}_1 (X_0 ,\mathbb{R})& \\
& \updownarrow \mbox{\rm{dual}} & &&  &\updownarrow \mbox{\rm{dual}} &\\
&\mathrm{Hom} (\Gamma, \mathbb{R}) &  & ~~
\stackrel{\,^t\rho_{\mathbb{R}}}{\hookrightarrow}& & {\rm H}^1(X_0 , \mathbb{R}) 
\cong
\big (
\mathcal{H}^1 (X_0), \langle \hspace{-0.8mm} \langle \cdot, \cdot \rangle \hspace{-0.8mm} \rangle \big )
&
\end{align*}
We write ${\rm Alb}^{\Gamma}$ for $(\Gamma \otimes \mathbb R
/ \Gamma \otimes \mathbb Z, g_{0})$, and call it the $\Gamma$-{\it{Albanese torus}} associated with
$(X, \Gamma)$.

Now we realize $X$ in $\Gamma \otimes \mathbb R$ 
equipped with the Albanese metric $g_0$ in a standard way. 
A (piecewise linear) map ${\Phi}: X \to \Gamma \otimes \mathbb R$ is said to
be a {\it{periodic realization}} of $X$ if it satisfies
\begin{equation}
{\Phi}(\sigma x)={\Phi}(x)+\sigma \otimes 1 \qquad (x\in X, \hspace{0.8mm} \sigma \in \Gamma).
 \label{periodicity}
\nonumber
\end{equation}
We may define a special periodic realization 
${{\Phi}}_{0}:X \to \Gamma \otimes \mathbb R$ 
by
${{\Phi}}_{0}(x_{*})={\bf 0}$ for a fixed base point $x_{*}\in V$ and
\begin{equation}
\subscripts
 {{\rm Hom}(\Gamma,\mathbb R)}
{\big \langle \omega , {{\Phi}}_{0}(x) \big \rangle}
{\Gamma \otimes \mathbb R}
=\int_{x_{*}}^{x} {\widetilde{\omega}} \qquad ( \omega \in {\rm Hom}(\Gamma, \mathbb R) ),
\label{def-sr}
\end{equation}
where $ \widetilde{\omega}$ is the lift of
$\omega=\hspace{-1mm}\mbox{ }^{t}\rho_{\mathbb R}(\omega) \in {\rm H}^{1}(X_{0}, {\mathbb R})$
to $X$. Here
$$ \int_{x_{*}}^{x} {\widetilde{ \omega}}=\int_{c} {\widetilde{ \omega}}
:= \sum_{i=1}^{n} {\widetilde{ \omega}}(e_{i})$$
for a path $c=(e_{1}, \ldots, e_{n})$ with $o(e_{1})=x_{*}$ and $t(e_{n})=x$.
It should be noted that this line integral does not depend on the choice of a path $c$.

One of the special properties of ${{\Phi}}_{0}$
is that
it is a vector-valued {\it{modified}}-{\it{harmonic function}}
on $X$ in the sense that
\begin{equation}
L{\Phi}_{0} (x)-{\Phi}_{0} (x)
= \rho_{\mathbb R} (\gamma_{p})
\qquad (x\in V).
\label{Phi-mh}
\end{equation} 
Indeed, for every
$\omega= \hspace{-1mm}\mbox{ }^{t}\rho_{\mathbb R}(\omega) \in {\rm H}^{1}(X_{0}, \mathbb R)$,
the modified harmonicity (\ref{M-harmonic}), $\Gamma$-invariance of the transition probability $p$
and the identity
(\ref{def-sr}) imply
\begin{eqnarray}
\subscripts
 {{\rm Hom}(\Gamma,\mathbb R)}
{\big \langle \omega , L{{\Phi}}_{0}(x)-
{{\Phi}}_{0}(x)
 \big \rangle}
{\Gamma \otimes \mathbb R}
&=&
\sum_{e\in E_{x}}
p(e) 
\subscripts
 {{\rm Hom}(\Gamma,\mathbb R)}
{\big \langle \omega , {{\Phi}}_{0}(t(e))-
{{\Phi}}_{0}(o(e))
 \big \rangle}
{\Gamma \otimes \mathbb R}
\nonumber \\
&=& \sum_{e\in E_{x}} p(e) {\widetilde \omega} (e)
\nonumber \\
&=& \sum_{e\in (E_{0})_{\pi (x)}} p(e) \omega (e)
\nonumber \\
&=&
 -(\delta_{p} \omega)(\pi(x))
\nonumber \\
&=&
 \langle \gamma_{p}, \omega \rangle
=
\subscripts
 {{\rm Hom}(\Gamma,\mathbb R)}
{\big \langle  \omega , \rho_{\mathbb R}(\gamma_{p})
 \big \rangle}
{\Gamma \otimes \mathbb R} \qquad (x\in V).
\nonumber
\end{eqnarray}
A periodic realization
${\Phi}: X \to \Gamma \otimes \mathbb R$ satisfying (\ref{Phi-mh})
is said to be {\it{modified harmonic}}.
Note that a modified harmonic realization
is uniquely determined up to
translation.

If we equip $\Gamma\otimes \mathbb{R}$ with the Albanese metric $g_0$, 
then we call the map $\Phi_0: X \to (\Gamma \otimes \mathbb{R}, g_0)$ 
the {\it modified standard realization} of $X$.
We readily check that the piecewise linear interpolation of ${{\Phi}}_{0}$ by
line segments descends to a
piecewise geodesic map $\Phi_{0}:X_{0} \to 
{\rm Alb}^{\Gamma}$. 
We call $\Phi_{0}$ the {\it{Albanese map}} associated with $(X, \Gamma)$. 
Namely, standard realization is a lift of the Albanese map.
In the symmetric case, Kotani--Sunada \cite{KS00, KS-contemp}
gave a characterization of the Albanese map $(\Phi_{0},g_{0})$ 
through a ``minimal principle" for the {\it{energy}} 
$$ {\cal E}(\Phi, g):=\frac{1}{2} \sum_{e\in E_{0}} \vert \Phi(t(\widetilde e))
-\Phi(o(\widetilde e)) \vert_{g}^{2} {\widetilde m}(e)
$$
under the fixed volume condition ${\rm vol}(\Gamma \otimes \mathbb R/ \Gamma \otimes \mathbb Z,g)
={\rm vol}({\rm Alb}^{\Gamma})$, where $g$ is a flat metric on  
$\Gamma \otimes \mathbb R/ \Gamma \otimes \mathbb Z$.
In the forthcoming paper \cite{IK}, we will discuss such 
a variational characterization of the Albanese map $(\Phi_{0}, g_{0})$ in the non-symmetric case.
\subsection{Ergodic theorems}
%
In this subsection, we make a preparation from ergodic theorems 
which will be used in the proof of main results. 
We set $\ell^{2}(X_{0}):=\{ f: V_{0} \to \mathbb C \}$,
which is equipped with the inner product
$$\langle f_{1},f_{2} \rangle_{\ell^{2}(X_{0})}
:=\sum_{x\in V_{0}} f_{1}(x) {\overline{f_{2}(x)}} \qquad (f_{1}, f_{2} \in \ell^{2}(X_{0}) ).$$
Although the following ergodic theorem is standard,
we give a proof for the sake of completeness.
\begin{tm} \label{Ergord-tm1}
Let $L$ be the transition operator on $X_{0}$ associated with
the transition probability $p$. Then it holds
\begin{equation}
\frac{1}{N}\sum_{j=0}^{N-1} L^{j}
f(x)=\sum_{x\in V_{0}} f(x)m(x)+O\big(\frac{1}{N} \big)
\qquad (x\in V_{0},~ f\in \ell^{2}(X_{0}) ).
\label{thm-ergord-1}
\nonumber
\end{equation}
\end{tm}
{\bf Proof.}~We denote by $K_{0}$ ($1\leq K_{0} \leq \vert V_{0} \vert$) the period of the random walk
on $X_{0}$, and set 
$$ \alpha_{k}:=\exp \big (\frac{2\pi k}{K_{0}} {\sqrt{-1}} \big) \qquad (k=0, \ldots, K_{0}-1 ),
$$
where $\vert V_{0} \vert$ stands for the number of elements of $V_{0}$.
By virtue of the Perron--Frobenius theorem, the transition operator
$L: \ell^{2}(X_{0}) \to \ell^{2}(X_{0})$ has the maximal simple eigenvalues
$\alpha_{0}, \ldots, \alpha_{K_{0}-1}$ with the corresponding normalized right eigenfunctions
$\phi_{0}, \ldots, \phi_{K_{0}-1}$ and left eigenfunctions $\psi_{0}, \ldots, \psi_{K_{0}-1}$. Namely,
$$ L\phi_{k}=\alpha_{k} \phi_{k},~~\hspace{-2mm}\mbox{ }^{t}L \psi_{k}=
{\overline{\alpha_{k}}} \psi_{k},~~  
\Vert \phi_{k} \Vert_{\ell^{2}(X_{0})}
=\langle \phi_{k}, \psi_{k} \rangle_{\ell^{2}(X_{0})}=1 \qquad (k=0,\ldots, K_{0}-1).
$$
In particular, we obtain
$\phi_{0}\equiv \vert V_{0} \vert^{-1/2}$ and $\psi_{0}(x)= \vert V_{0} \vert^{1/2} m(x)$  for $x\in V_{0}$.

Now we set $${\ell^{2}_{K_{0}}(X_{0})}:=\{ f\in \ell^{2}(X_{0}) \vert~\langle f, \psi_{k} \rangle_{\ell^{2}(X_{0})}=0
~~(k=0, \ldots, K_{0}-1) \}.$$
Note that $\ell^{2}_{K_{0}}(X_{0})$ is preserved by $L$. Thus $f\in \ell^{2}(X_{0})$ is decomposed as
\begin{eqnarray}
f&=&\sum_{k=0}^{K_{0}-1}  \langle f, \psi_{k} \rangle_{\ell^{2}(X_{0})}  \phi_{k} +f_{{\ell^{2}_{K_{0}}(X_{0})}}
\nonumber \\
&=& \langle f, m \rangle_{\ell^{2}(X_{0})}
+\sum_{k=1}^{K_{0}-1}  \langle f, \psi_{k} \rangle_{\ell^{2}(X_{0})}  \phi_{k} +f_{{\ell^{2}_{K_{0}}(X_{0})}}
\qquad 
(f_{{\ell^{2}_{K_{0}}(X_{0})}}\in {\ell^{2}_{K_{0}}(X_{0})}).
\label{f-koyu-bunkai-L}
\end{eqnarray}

Besides, it follows from the Perron--Frobenius theorem that
there exists some $0\leq \lambda <1$
such that $\Vert L\vert_{{\ell^{2}_{K_{0}}(X_{0})}}
 \Vert \leq \lambda$.
Furthermore by noting
$\sum_{j=0}^{K_{0}-1} \alpha_{k}^{j} =0$ 
$(k=1,\ldots, K_{0}-1)$
and (\ref{f-koyu-bunkai-L}),
we have
\begin{eqnarray}
\lefteqn{
\Big \Vert \frac{1}{N} \sum_{j=0}^{N-1}L^{j} f - \langle f, m \rangle_{\ell^{2}(X_{0})} 
\Big \Vert_{\ell^{2}(X_{0})}
}
\nonumber \\
&=&  
\Big \Vert
\frac{1}{N}
\sum_{k=1}^{K_{0}-1} 
\langle f, \psi_{k} \rangle_{\ell^{2}(X_{0})} 
\big (
\sum_{j=0}^{N-1} \alpha_{k}^{j} \big) \phi_{k}
+
\frac{1}{N}
\sum_{j=0}^{N-1} L^{j} 
f_{{\ell^{2}_{K_{0}}(X_{0})}}
\Big \Vert_{\ell^{2}(X_{0})}
\nonumber \\
&\leq & \frac{1}{N} 
\Big \{ \sum_{k=1}^{K_{0}-1}
\big \vert \langle f, \psi_{k} \rangle_{\ell^{2}(X_{0})}  \big \vert
\Big(
\big \vert \sum_{j=0}^{N-1} \alpha_{k}^{j} \big \vert 
\Big)
\Big\}
+\Big( \frac{1}{N} \sum_{j=0}^{N-1}\Vert L\vert_{{\ell^{2}_{K_{0}}(X_{0})}} \Vert^{j} \Big) 
\Vert 
f_{{\ell^{2}_{K_{0}}(X_{0})}}
\Vert_{\ell^{2}(X_{0})}
\nonumber \\
&\leq & \frac{K_{0}-1}{N}
\Big( \sum_{k=1}^{K_{0}-1}
\vert \langle f, \psi_{k} \rangle_{\ell^{2}(X_{0})} \vert \Big)+
\frac{1}{(1-\lambda)N}
\Vert 
f_{{\ell^{2}_{K_{0}}(X_{0})}}
\Vert_{\ell^{2}(X_{0})}
=O(\frac{1}{N}).
\label{f-fourier}
\nonumber
\end{eqnarray}
Thus we complete the proof.
\qed
\vspace{2mm} 

In the proof of the CLT of the second kind (Theorem \ref{CLT-2}), the following ergodic theorem plays
a crucial role.
\begin{tm}\label{Ergord-tm2}
Let $L_{(\varepsilon)}$ be the transition operator on $X_{0}$ associated with
the transition probability $p_{\varepsilon}$ given in {\rm{(\ref{def-pe})}}. 
Then there exists sufficiently small
$\varepsilon_{0} >0$ 
such that 
\begin{equation}
\frac{1}{N}\sum_{j=0}^{N-1} L^{j}_{(\varepsilon)}
f(x)=\sum_{x\in V_{0}} f(x)m(x)+O_{\varepsilon_{0}}\big(\frac{1}{N} \big)
\qquad (x\in V_{0}, ~f\in \ell^{2}(X_{0}) )
\nonumber
\end{equation}
holds for all $0\leq \varepsilon \leq \varepsilon_{0}$.
\end{tm}
To prove this theorem, we need
a fundamental perturbation theory of linear operators
taken from Parry--Pollicott \cite[Proposition 4.6]{PP}
\begin{pr} \label{PP-4.6}
Let $B({\cal V})$ denote the set of linear operators on a Banach space
${\cal V}$. Asuume ${\cal L}_{0} \in B({\cal V})$ has a simple isolated eigenvalue $\alpha_{0}$
with corresponding eigenvector $\phi_{0}$. Then for any $\varepsilon >0$, there exists $\delta>0$
such that the following identities hold for all ${\cal L}\in B({\cal V})$ with $\Vert {\cal L}-{\cal L}_{0} \Vert_{B({\cal V})}<\delta$:
\\
{\rm{(1)}}~${\cal L}$ has a simple isolated eigenvalue $\alpha ({\cal L})$ and the corresponding
eigenvector $\phi({\cal L})$ with 

$\alpha ({\cal L}_{0})=\alpha$,  $\phi ({\cal L}_{0})=\phi_{0}$,
 \\
{\rm{(2)}}~${\cal L} \mapsto \alpha({\cal L})$, ${\cal L} \mapsto \phi({\cal L})$ are analytic,
\\
{\rm{(3)}}~$\vert \alpha({\cal L})
-\alpha_{0} \vert < \varepsilon$  and  ${\rm{Spec}}({\cal L}) \setminus \{ \alpha({\cal L}) \}
\subset \{ z\in {\mathbb C}:  \vert z-\alpha_{0} \vert > \varepsilon \}$.
\end{pr}
{\bf{Proof of Theorem \ref{Ergord-tm2}.}}
Let $0\leq \varepsilon <1$. It follows from
$p_{\varepsilon}(e)>0$ ($e\in E_{0}$) that
$$
K_{0}
=
\begin{cases} 1 & \text{ if $X_{0}$ is non-bipartite},
\\
2 & \text{ if $X_{0}$ is bipartite}.
\end{cases}
$$
We only consider the bipartite case because the 
argument for the non-bipartite case is same.
Let $V_{0}=A_{0} \coprod A_{1}$ be the bipartition.  
By the Perron--Frobenius theorem, the the transition operator 
$L_{(\varepsilon)}: \ell^{2}(X_{0})
\to \ell^{2}(X_{0})$ ($0\leq \varepsilon <1$) has the maximal simple eigenvalues 
$\alpha_{0}(\varepsilon)=1, \alpha_{1}(\varepsilon)=-1$ with the 
corresponding normalized right eigenfunctions 
$$\phi_{0}(\varepsilon)\equiv \vert V_{0} \vert^{-1/2}, \qquad 
\phi_{1}(\varepsilon)(x) 
=
\begin{cases}  \vert V_{0} \vert^{-1/2} 
& \text{ ($x\in A_{0}$) },
\\
-\vert V_{0} \vert ^{-1/2} 
& \text{ ($x\in A_{1}$) },
\end{cases}
$$
and left eigenfunctions 
$$ 
\psi_{0}(\varepsilon)(x)=\vert V_{0} \vert^{1/2}m(x) \quad (x\in V_{0}), 
\qquad 
\psi_{1}(\varepsilon)(x) 
=
\begin{cases}  \vert V_{0} \vert^{1/2}m(x) 
& \text{ ($x\in A_{0}$) },
\\
-\vert V_{0} \vert^{1/2}m(x)  
& \text{ ($x\in A_{1}$) }.
\end{cases}
$$
Then we see that the subspace ${\cal U}:=\{ f\in \ell^{2}(X_{0}) \vert 
\langle f, \psi_{0} \rangle_{\ell^{2}(X_{0})}=
\langle f, \psi_{1} \rangle_{\ell^{2}(X_{0})}=0\}$
is independent of $\varepsilon$ and  it is preserved by $L_{(\varepsilon)}$.
By virtue of Proposition \ref{PP-4.6}, 
there exist sufficiently small $\varepsilon_{0}>0$ and $0\leq \lambda <1$
such that $\Vert L_{(\varepsilon)} \vert_{\cal U} \Vert \leq \lambda$ holds for 
all $0\leq \varepsilon \leq \varepsilon_{0}$. 
Hence we obtain the desired
conclusion by following the proof of Theorem \ref{Ergord-tm1}.
\qed
\subsection{Trotter's approximation theorem}
In this subsection,
we quickly review an approximation
theorem due to Trotter \cite{Trotter}.
Let ${\cal V}$ and ${\cal V}_{N}$ ($N\in \mathbb N$) be Banach spaces
and $B({\cal V}, {\cal V}_{N})$ stands for the set of bounded linear operators
from ${\cal V}$ to ${\cal V}_{N}$. Suppose that 
there exist a constant $C>0$ and two families of operators
$\{ P_{N}\in B({\cal V}, {\cal V}_{N})\}_{N=1}^{\infty}$ and 
$\{ U_{N} \in B({\cal V}_{N}) \}_{N=1}^{\infty}$ satisfying
$$\sup_{N\in \mathbb N} \Vert P_{N} \Vert_{B({\cal V}, {\cal V}_{N})} 
+\sup_{N,n\in \mathbb N} \Vert U_{N}^{n} \Vert_{B({\cal V}_{N})}
\leq C,~~
\lim_{N \to \infty} \Vert P_{N} u \Vert_{{\cal V}_{N}} =\Vert u \Vert_{\cal V}
\qquad (u\in \cal V).$$
Trotter's approximation theorem is then presented as follows:
\begin{tm}
[cf. \cite{Trotter, Kotani, KS-contemp}]
\label{Trotter-tm}
Let $\{e^{-tT}
\}_{t\geq 0}$ be a continuous semigroup on $\cal V$
with the infinitesimal generator $T$ and $\{\tau_{N}\}_{N=1}^{\infty}$
be a decreasing sequence satisfying $\lim_{N\to \infty} \tau_{N}=0$. Assume
that there exists a core $D$ of $T$ such that
$$ \lim_{N\to \infty} \Vert T_{N}P_{N}u -P_{N}Tu \Vert_{{\cal V}_{N}}=0 \qquad (u\in D),$$
where $T_{N}:=(1/\tau_{N})(I-U_{N})$. Then for any sequence $\{k_{N} \}_{N=1}^{\infty}$
of non-negative integers satisfying $k_{N}\tau_{N} \to t$, we have
$$ \lim_{N\to \infty} \Vert U^{k_{N}}_{N}P_{N} u-P_{N}e^{-tT}u \Vert_{{\cal V}_{N}}=0 
\qquad (u\in D).$$
\end{tm}
%
\section{Proof of the CLT of the first kind}
\label{CLT1}
In the following, we write
$
\omega [ {\bf x} ]_{\Gamma \otimes \mathbb R}:=
\subscripts
    {{\rm Hom}(\Gamma,\mathbb R)}
    {\langle \omega, {\bf{x}}  \rangle}
    {\Gamma \otimes \mathbb R}$
($ \omega\in {\rm Hom}(\Gamma,\mathbb R),  {\bf x}\in \Gamma \otimes \mathbb R$)    
and set $$d\Phi_{0}(e):=\Phi_{0}(t(e))-\Phi_{0}(o(e)) \quad (e\in E), \qquad
 \Vert d\Phi_{0} \Vert_{\infty}:= \max_{e\in E_{0}} \vert d\Phi_{0}({\widetilde{e}})
\vert_{g_{0}}.$$
We take an orthonormal basis $\{ \omega_{1}, \ldots, \omega_{d} \}$ of 
${\rm{Hom}}(\Gamma , \mathbb R)(\subset {\rm H}^{1}(X_{0},\mathbb R) \cong {\cal H}^{1}(X_{0}))$, 
and let $\{ {\bf{v}}_{1}, \ldots, {\bf{v}}_{d} \}$ 
denote its dual basis in $\Gamma \otimes \mathbb R$. Namely,
$ \omega_{i} [ {\bf{v}}_{j} ]_{\Gamma \otimes \mathbb R}
=\delta_{ij}$~($i,j=1,\ldots, d$).
Note that $\{{\bf{v}}_{1}, \ldots, {\bf{v}}_{d} \}$ is an orthonomal basis of 
$\Gamma \otimes \mathbb R$ with respect to the 
Albanese metric $g_{0}$. Then we have
$$ {\bf{x}}=\sum_{i=1}^{d} \omega_{i} [{\bf{x}} ]_{\Gamma \otimes \mathbb R}
 {\bf{v}}_{i}, \quad 
\vert {\bf{x}} \vert_{g_{0}}^{2}=
\sum_{i=1}^{d} 
\omega_{i} [{\bf{x}} ]_{\Gamma \otimes \mathbb R}^{2}
\qquad ({\bf{x}} \in \Gamma \otimes \mathbb R).
$$
For simplicity, we also write 
$x_{i}:=\omega_{i} [{\bf{x}} ]_{\Gamma \otimes \mathbb R}$, 
$\Phi_{0}(x)_{i}:= \omega_{i}[\Phi_{0}(x)]_{\Gamma \otimes \mathbb R}$
($i=1,\ldots, d$, $x\in V$)
and identify ${\bf x} \in \Gamma \otimes \mathbb R$ with 
$(x_{1}, \ldots, x_{d}) \in \mathbb R^{d}$. 
\subsection{Proof of Theorem \ref{CLT-1}}
To prove Theorem \ref{CLT-1}, we need the following lemma:
\begin{lm}
\label{CLT-3-1}
For any $f \in C_0^\infty( \Gamma \otimes \mathbb{R})$, as 
$N \nearrow \infty$, $\varepsilon \searrow 0$ and 
$N^{2}\varepsilon  \searrow 0$, we have
\begin{equation}
\left\| \frac{1}{N\varepsilon^2} \left(I-{\cal L}_{\gamma_{p}}^N\right) {\cal P}_\varepsilon f 
-{\cal P}_\varepsilon \big( \frac{\Delta}{2}f \big)  \right\|_\infty  \to 0,
\label{conv-CLT3-1}
\nonumber
\end{equation}
where $\Delta$ is the {\rm{(}}positive{\rm{)}} 
Laplacian $-\sum_{ i = 1}^{d}  \big( \frac{ \partial}{\partial x_i} \big)^{2}$  
on $\Gamma \otimes \mathbb{R}$ with the Albanese metric $g_{0}$.
\end{lm}
{\bf Proof.}~  
First, we define $A^{N}(\Phi_{0})_{ij}:V \to \mathbb R$ ($i,j=1, \ldots, d$, $N\in \mathbb N$) 
by
\begin{eqnarray}
A^{N}(\Phi_{0})_{ij}(x)
&:= &
\sum_{c\in {\Omega}_{x,N}(X)} p(c) 
\big( \Phi_{0} ( t(c) ) - \Phi_{0} (x)-N \rho_{\mathbb R}(\gamma_p)  \big)_i
\nonumber \\
&\mbox{ }&
\times
\big( \Phi_{0} ( t(c) ) - \Phi_{0} (x)-N \rho_{\mathbb R}(\gamma_p)  \big)_j
\qquad (x\in V),
\nonumber
\end{eqnarray}
where $p(c):=p(e_1)p(e_2)\cdots p(e_N)$ for  
$c=(e_1, e_2, \ldots, e_N) \in {\Omega}_{x,N}(X)$.
Applying Taylor's expansion formula to 
$f\big (\varepsilon (\Phi_{0}(t(c))  -\rho_\mathbb{R}( \mathbf{z} +N\gamma_p)\big)$ at 
$\varepsilon (\Phi_{0} (x) - \rho_\mathbb{R} (\mathbf{z})) \in \Gamma \otimes \mathbb R$, 
we have
\begin{eqnarray}
\lefteqn{f \big (\varepsilon (\Phi_{0}(t(c))  -\rho_\mathbb{R}( \mathbf{z} +N\gamma_p) \big )
- f \big (\varepsilon (\Phi_{0} (x) -\rho_\mathbb{R}(\mathbf{z}) ) \big)
}
\nonumber \\
&=&
\sum_{i=1}^{d}
\frac{ \partial f}{\partial x_i}  \big( \varepsilon (\Phi_{0} (x) -\rho_\mathbb{R}(\mathbf{z})) \big)
\Big( \varepsilon \big ( \Phi_{0}(t(c)) - \Phi_{0} (x) -N \rho_\mathbb{R}( \gamma_p) \big) \Big)_i
\nonumber \\
&\mbox{ }&+\frac{1}{2} \sum_{i,j=1}^{d} 
\frac{\partial^2 f}{\partial x_i \partial x_j} \big( \varepsilon (\Phi_{0} (x) -\rho_\mathbb{R}(\mathbf{x})) \big)
\Big( \varepsilon \big( \Phi_{0}(t(c)) - \Phi_{0} (x) -N \rho_\mathbb{R}( \gamma_p) \big)  \Big)_i
\nonumber \\
&\mbox{ }&
\quad
\times 
\Big( \varepsilon \big( \Phi_{0}(t(c)) - \Phi_{0} (x) -N \rho_\mathbb{R}( \gamma_p) \big) \Big)_j
+O(N^3\varepsilon ^3). \nonumber 
\end{eqnarray}
This implies
\begin{eqnarray}
\lefteqn{ 
(I-{\cal L}_{\gamma_{p}}^N) {\cal P}_\varepsilon f (x , \mathbf{z}) }
\nonumber \\
&=&
- \varepsilon \sum_{i=1}^{d} \frac{ \partial f}{\partial x_i}  \big( \varepsilon (\Phi_{0} (x) 
-\rho_\mathbb{R}(\mathbf{z})) \big)
\sum_{ c \in  {\Omega}_{x,N}(X)} p(c) \big( \Phi_{0}(t(c)) - \Phi_{0} (x) 
-N \rho_\mathbb{R}( \gamma_p) \big)_i  
\nonumber \\
&\mbox{ }&
-\frac{\varepsilon^{2}}{2} \sum_{i,j=1}^{d} 
\frac{ \partial^2 f}{\partial x_i \partial x_j} 
\big( \varepsilon (\Phi_{0} (x) -\rho_\mathbb{R}(\mathbf{z})) \big)
A^{N}(\Phi_{0})_{ij}(x)
+O(N^3 \varepsilon^3 ) . 
\label{core-expand}
\nonumber
\end{eqnarray}

Recalling the modified harmonicity (\ref{Phi-mh}),
we obtain
\begin{eqnarray}
\lefteqn{\sum_{ c \in  {\Omega}_{x,N}(X)} p(c)
\big( \Phi_{0}(t(c)) - \Phi_{0} (x) -N \rho_\mathbb{R}( \gamma_p) \big)
}
\nonumber \\
&=&
\sum_{c' \in {\Omega}_{x, N-1}(X)} 
p(c')
\Big \{
\Big (
 \sum_{e \in E_{t(c^\prime )}} 
p(e) \big( \Phi_{0} ( t(e) ) -\Phi_{0} (o(e))  - \rho_\mathbb{R} (\gamma_p) \big) \Big )
\nonumber \\
&\mbox{ }& \quad +\big ( \Phi_{0}(t(c^\prime))  -\Phi_{0}(x)  -(N-1) 
\rho_{\mathbb R}(\gamma_p) \big) \Big \}
\nonumber \\
&=& \sum_{ c^\prime \in {\Omega}_{x, N-1}(X)} p(c^\prime )
\big(  \Phi_{0} ( t(c^\prime ))    - \Phi_{0} (x)-(N-1) \rho_{\mathbb R}(\gamma_p)\big)
\nonumber \\
&=& \sum_{ e \in E_{x}} p(e)
\big(  \Phi_{0} ( t(e)) - \Phi_{0} (x)-\rho_{\mathbb R}(\gamma_p)\big) ={\bf 0}.
\label{expand2-2}
\end{eqnarray}

Next we define ${\cal A}(\Phi_{0})_{ij}: V_{0} \to \mathbb R$ ($i,j=1,\ldots, d$) by
\begin{eqnarray}
{\cal A}(\Phi_{0})_{ij}(x)
&=&\sum_{e \in (E_{0})_{x}} p(e) 
\left( \Phi_{0} ( t(\widetilde{e}) ) -\Phi_{0} (o(\widetilde{e})) -\rho_\mathbb{R}(\gamma_p) \right)_i
\nonumber \\
&\mbox{  }&
\quad
\times
     \left( \Phi_{0} ( t(\widetilde{e}) ) -\Phi_{0}(o(\widetilde{e})) -\rho_\mathbb{R}(\gamma_p) \right)_j 
\qquad (x\in V_{0}),
\nonumber
\end{eqnarray}
where $\widetilde{e}$ is a lift of $e \in E_0$ to $E$. 
Because $A^{N}(\Phi_{0})_{ij}:V \to \mathbb R$ 
is $\Gamma$-invariant, we easily see
$$
{\cal A}(\Phi_{0})_{ij}(\pi(x))=A^{1}(\Phi_{0})_{ij}(x) 
\qquad (x\in V,~ i,j=1,\ldots, d).
$$

Repeating the same argument as in (\ref{expand2-2}) and recalling (\ref{Phi-mh}) again, 
we obtain
\begin{eqnarray}
\lefteqn{
A^{N}(\Phi_{0})_{ij}(x)
}
\nonumber\\
&=& \sum_{ c^\prime \in {\Omega}_{x,N-1}(X)}
p(c^\prime) 
\sum_{e  \in E_{t(c^\prime)}}
p(e) \nonumber \\
&\mbox{ }& \times 
\Big \{ \big( \Phi_{0} ( t(e) ) -\Phi_{0}(o(e)) -
\rho_{\mathbb R}(\gamma_p) \big)_{i} 
+ \big ( \Phi_{0} (o(e)) - \Phi_{0} (x)-(N-1) \rho_{\mathbb R}(\gamma_p) \big)_i \Big \}  
\nonumber\\
&\mbox{ }& \times 
\Big \{ \big( \Phi_{0} ( t(e) ) -\Phi_{0}(o(e)) -\rho_{\mathbb R}(\gamma_p) 
\big)_{j} + \big ( \Phi_{0}(o(e)) - \Phi_{0} (x)-(N-1) \rho_{\mathbb R}(\gamma_p) \big)_j \Big \}  
\nonumber\\
&=& \sum_{ c^\prime \in {\Omega}_{x,N-1}(X)}
p(c^\prime) 
\sum_{e  \in E_{t(c^\prime)}}p(e)
\nonumber \\
&\mbox{ }& \quad \times
     \big( \Phi_{0} ( t(e) ) -\Phi_{0}(o(e)) -\rho_{\mathbb R}(\gamma_p) \big)_i
     \big( \Phi_{0} ( t(e) ) -\Phi_{0} (o(e)) -\rho_{\mathbb R} (\gamma_p) \big)_j  
\nonumber \\
&\mbox{ }&  + \sum_{ c^\prime \in {\Omega}_{x,N-1}(X)}
p(c^\prime) 
\sum_{e  \in E_{t(c^\prime)}} p(e) 
\nonumber \\
&\mbox{ }& \quad \times  
\Big\{
\big( \Phi_{0} ( t(e) ) -\Phi_{0}(o(e)) -\rho_{\mathbb R}(\gamma_p)  \big)_i 
\big(  \Phi_{0} (o(e))  - \Phi_{0} (x)-(N-1) \rho_{\mathbb R} (\gamma_p) \big)_j    \nonumber \\
&\mbox{ }& \qquad +
\big( \Phi_{0} ( t(e) ) -\Phi_{0}(o(e)) -\rho_{\mathbb R}(\gamma_p)  \big)_i 
\big( \Phi_{0} (o(e)) - \Phi _{0} (x) -(N-1) \rho_{\mathbb R} (\gamma_p) \big)_j 
\Big\}  
\nonumber \\
&\mbox{ }&  +\sum_{ c^\prime \in {\Omega}_{x,N-1}(X)}
p(c^\prime) 
\sum_{e  \in E_{t(c^\prime)}} p(e) 
\big(\Phi_{0}(o(e))  - \Phi_{0} (x)-(N-1) \rho_{\mathbb R} (\gamma_p) \big)_i
\nonumber \\
&\mbox{ }& \qquad \times
\big(\Phi_{0} (o(e)) - \Phi_{0} (x)-(N-1) \rho_{\mathbb R} (\gamma_p)  \big)_j
\nonumber \\
&=& L^{N-1} \big ({\cal A}(\Phi_{0})_{ij} \big )(\pi(x))+A^{N-1}(\Phi_{0})_{ij}(x)
\nonumber \\
&=&
\sum_{k=0}^{N-1} L^{k} \big ({\cal A}(\Phi_{0})_{ij} \big )(\pi(x)) 
\qquad \qquad (x\in V).
\label{A-cal}
\nonumber
\end{eqnarray}
Applying Theorem \ref{Ergord-tm1}, we have
\begin{equation}
\frac{1}{N} \sum_{k=0}^{N-1} 
L^{k} \big ({\cal A}(\Phi_{0})_{ij} \big )(\pi(x))
=\sum_{x \in V_0} {\cal A}(\Phi_{0})_{ij}(x) m(x) +O\big (\frac{1}{N} \big).
\label{ergord-1}
\nonumber
\end{equation}
Moreover (\ref{def-sr}) and (\ref{Phi-mh}) imply
\begin{eqnarray}
\lefteqn{
\sum_{x \in X_0} {\cal A}(\Phi_{0})_{ij}(x) m(x)}
\nonumber \\
&=& 
 \sum_{e \in E_0} \big( \Phi_{0} ( t(\widetilde{e}) ) 
 -\Phi_{0} (o(\widetilde{e}))-{\rho}_{\mathbb R}(\gamma_{p})
 \big)_{i}
 \big( \Phi_{0} ( t(\widetilde{e}) ) -\Phi_{0} (o(\widetilde{e}))
 -{\rho}_{\mathbb R}(\gamma_{p})
 \big)_{j} {\widetilde{m}}(e)
 \nonumber \\
 &=&
  \sum_{e \in E_0} \big( \Phi_{0} ( t(\widetilde{e}) ) -
  \Phi_{0}(o(\widetilde{e}))
 \big)_{i}
 \big( \Phi_{0} ( t(\widetilde{e}) ) -\Phi_{0} (o(\widetilde{e}))
 \big)_{j}  {\widetilde{m}}(e)- 
 {\rho}_{\mathbb R}(\gamma_{p})_{i}{\rho}_{\mathbb R}(\gamma_{p})_{j}
\nonumber \\
 &=& 
 \sum_{e \in E_0}  \hspace{-1mm}\mbox{ }^{t}\rho_{\mathbb R}(\omega_{i}) (e)
 \hspace{-1mm}\mbox{ }^{t}\rho_{\mathbb R}(\omega_{j}) (e) 
  {\widetilde{m}}(e)
 - \omega_{i} [ \rho_{\mathbb R}(\gamma_{p})] _{\Gamma \otimes \mathbb R} 
 \hspace{1mm}
 \omega_{j} [ \rho_{\mathbb R}(\gamma_{p})] _{\Gamma \otimes \mathbb R} 
\nonumber       \\
&=& \sum_{ e \in E_0} \omega_i (e) \omega_j (e) {\widetilde{m}}(e) 
-\langle \gamma_p, \omega_{i}  \rangle \langle  \gamma_p, \omega_{j} \rangle 
\nonumber \\
&=& \langle \hspace{-0.8mm} \langle \omega_i , \omega_j \rangle \hspace{-0.8mm} \rangle 
=\delta_{ij} \qquad \quad (i,j=1,\ldots, d).
\label{B-cal}
\nonumber
 \end{eqnarray}

Putting it all together, we now obtain
\begin{align*}
&\frac{1}{N\varepsilon^2} (I-{\cal L}_{\gamma_{p}}^N) {\cal P}_\varepsilon f (x , \mathbf{z})
={\cal P}_\varepsilon \big( \frac{\Delta }{2} f\big) ( x, \mathbf{z}) 
+O(N^2\varepsilon)  \quad \mbox{as  }N\rightarrow \infty.
\end{align*}
Finally by letting $N^{2} \varepsilon \to 0$, we complete the proof.
\qed
\vspace{2mm} \\
{\bf Proof of Theorem \ref{CLT-1}.}~We follow the proof of \cite[Theorem 4]{Kotani}.
Let $N=N(n)$ be the integer with 
$n^{1/5}\leq N<1+n^{1/5}$
and $k_{N}$ and $r_{N}$ are the quotient and remainder of $([nt]-[ns])/N$, respectively. 
We put $\varepsilon_{N}:=n^{-1/2}$ and $\tau_{N}:=N\varepsilon_{N}^{2}$.
Then $n\to \infty$ implies $N\to \infty$, $N^{2}\varepsilon_{N}
\leq (1+n^{1/5})^{2}n^{-1/2} \to 0$ and $\tau_{N} \leq (1+n^{1/5})/n \to 0$.
As $r_{N}<N$, we also observe 
$r_{N}\varepsilon_{N}^{2} \leq N\varepsilon_{N}^{2}\leq (1+n^{1/5})/n
\to 0$.  Noting $k_{N}\tau_{N}=([nt]-[ns]-r_{N})\varepsilon_{N}^{2}$,
we obtain $k_{N}\tau_{N} \to (t-s)$ as $N\to \infty$.

Now we may apply Theorem \ref{Trotter-tm} to the case where 
$${\cal V}=C_{\infty}(\Gamma \otimes \mathbb R),~~
{\cal V}_{N}=C_{\infty}(X \times {\rm H}_{1}(X_{0}, \mathbb R)),~~ U_{N}={\cal L}_{\gamma_{p}}^{N}, ~~
T=\frac{\Delta}{2},~~D=C^{\infty}_{0}({\Gamma}\otimes \mathbb R),$$ 
and Lemma \ref{CLT-3-1} implies
\begin{equation}
\lim_{n \to \infty} \left \|
{\cal L}_{\gamma_{p}}^{k_{N}N}
{\cal P}_{n^{-1/2}} f - {\cal P}_{n^{-1/2}} 
e^{ -\frac{t-s}{2}\Delta} f \right\|_\infty =0 \qquad (f\in C^{\infty}_{0}({\Gamma}\otimes \mathbb R) ).\label{Trotter-2}
\end{equation}
Further, we have
\begin{eqnarray}
\lefteqn{
\Vert 
{\cal L}_{\gamma_{p}}^{[nt]-[ns]} {\cal P}_{n^{-1/2}} f-
{\cal P}_{n^{-1/2}} e^{ -\frac{t-s}{2}\Delta} f \Vert_\infty
}
\nonumber \\
&\leq &
\Vert {\cal L}_{\gamma_{p}}^{k_{N}N} ( {\cal L}_{\gamma_{p}}^{r_{N}}-I)
{\cal P}_{n^{-1/2}}f 
\Vert_{\infty}
+
\Vert {\cal L}_{\gamma_{p}}^{k_{N}N} 
{\cal P}_{n^{-1/2}}f -{\cal P}_{n^{-1/2}} e^{ -\frac{t-s}{2}\Delta} f \Vert_\infty.
\label{Trotter-23}
\end{eqnarray}
Noting $r_{N}^{2}\varepsilon_{N} \leq (1+n^{1/5})^{2}n^{-1/2} \to 0$ and
recalling 
Lemma \ref{CLT-3-1} again, we obtain
$$ \Big \Vert \frac{1}{r_{N}\varepsilon_{N}^{2}}
(I-{\cal L}_{\gamma_{p}}^{r_{N}}){\cal P}_{\varepsilon_{N}}f
-{\cal P}_{\varepsilon_{N}}
\big( \frac{\Delta}{2}f \big)  \Big \Vert_\infty  \to 0 \qquad (f\in C^{\infty}_{0}(\Gamma \otimes \mathbb R) ).
$$
This convergence and $r_{N}\varepsilon^{2}_{N} \to 0$ imply
\begin{align}
\Vert {\cal L}_{\gamma_{p}}^{k_{N}N} ( {\cal L}_{\gamma_{p}}^{r_{N}}-I)
{\cal P}_{n^{-1/2}}f 
\Vert_{\infty}
&\leq  r_{N}\varepsilon_{N}^{2}
 \Big \Vert \frac{1}{r_{N}\varepsilon_{N}^{2}}
(I-{\cal L}_{\gamma_{p}}^{r_{N}}){\cal P}_{\varepsilon_{N}}f
-{\cal P}_{\varepsilon_{N}}
\big( \frac{\Delta}{2}f \big)  \Big \Vert_\infty
\nonumber \\
&\mbox{ }+ r_{N}\varepsilon_{N}^{2} \big \Vert \frac{\Delta}{2}f  \big \Vert_\infty
\to 0.
\label{Trotter-3}
\end{align}
Hence by combining (\ref{Trotter-2}), (\ref{Trotter-23}) with (\ref{Trotter-3}), 
we obtain (\ref{Trotter-1}) for 
$f\in C^{\infty}_{0}({\Gamma}\otimes \mathbb R)$.

For $f\in C_{\infty}({\Gamma}\otimes \mathbb R)$, we can choose a sequence 
$\{ f_{m} \}_{m=1}^{\infty} \subset C^{\infty}_{0}({\Gamma}\otimes \mathbb R)$ 
such that $\Vert f-f_{m} \Vert_{\infty} \to 0$ as $m \to \infty$. 
Because
\begin{equation}
\max \Big \{
\big \Vert {\cal L}_{\gamma_{p}}^{[nt]-[ns]} {\cal P}_{n^{-1/2}} (f-f_{m}) \big \Vert_{\infty},
\Vert 
{\cal P}_{n^{-1/2}} 
e^{ -\frac{t-s}{2}\Delta} (f-f_{m}) \Vert_{\infty} \Big \}  
\leq \Vert f-f_{m} \Vert_{\infty} 
\label{228-2015}
\nonumber
\end{equation}
holds for each $n\in \mathbb N$, it is straightforward to check (\ref{Trotter-1}) for 
$f\in C_{\infty}({\Gamma}\otimes \mathbb R)$. 

In addition, we have
\begin{eqnarray}
\lefteqn{
\vert   {\cal L}_{\gamma_{p}}^{[nt]}{\cal P}_{n^{-1/2}} f(x_n,\mathbf{z}_n)
-e^{-t\Delta/2} f({\bf{x}}) \vert
}
\nonumber \\
&\leq &
\big \|  {\cal L}_{\gamma_{p}}^{[nt]} {\cal P}_{n^{-1/2}} f - {\cal P}_{n^{-1/2}} 
e^{ -\frac{t}{2}\Delta} f \big \|_\infty
\nonumber \\
&\mbox{ }&+\big \vert  e^{ -\frac{t}{2}\Delta} 
f \big( n^{-1/2}( \Phi_{0}(x_{n}) -\rho_{\mathbb R}({\bf{z}}_{n}) \big)-e^{ -\frac{t}{2}\Delta} f({\bf{x}})
\big \vert 
\label{1may3-2015}
\end{eqnarray}
for
$f\in C_{\infty}(\Gamma \otimes \mathbb R)$.
By combining the continuity of $e^{-\frac{t}{2}\Delta}f: \Gamma \otimes \mathbb R \to \mathbb R$ and
(\ref{Trotter-1}) with (\ref{1may3-2015}), we obtain (\ref{Trotter-1-2}). This completes the proof.
\qed
\subsection{
Proof of Theorem \ref{FCLT-1}}
At the beginning, we show the convergence of finite dimensional distribution of 
$\{ {\bf X}^{(n)} \}_{n=1}^{\infty}$
by using Theorem \ref{CLT-1}.
We fix $0\leq t_1 <  \ldots < t_r <\infty $ ($r\in \mathbb N$) and consider the random variable 
${\bf X}^{(n)}_{t_1, t_2, \ldots , t_r} : {\Omega}_{x_{*}}(X) 
\rightarrow (\Gamma\otimes \mathbb{R} )^r$
given by
\begin{equation*}
{\bf X}^{(n)}_{t_{1}, \ldots, t_{r}}(c) := \big( {\bf X}_{t_1}^{(n)}(c), 
\ldots , {\bf X}_{t_r}^{(n)}(c) \big).
\end{equation*}
\begin{lm}\label{FDD1}
\begin{equation*}
{\bf X}^{(n)}_{t_{1}, \ldots, t_{r}}
\stackrel{ \mathcal{D}}{\longrightarrow} ( B_{t_1}, \cdots, B_{ t_r} ) \quad \mbox{ as } n\to \infty,
\end{equation*}
where $( B_{t} )_{t\geq 0}$ is a 
$\Gamma \otimes \mathbb R$-valued standard Brownian motion with $B_{0}={\bf 0}$.
\end{lm}
{\bf Proof.}~ For simplicity, we only show the convergence for $r=2$. 
General cases differ from this one only by being notationally more cumbersome. Take 
$f_{1}=f_{1}({\bf{y}}_{1}), 
f_{2}=f_{2}({\bf{y}}_{2}) \in C_{\infty}(\Gamma \otimes \mathbb{R})$, and set 
$s=t_{1}, t=t_{2}$.
Let $f({\bf y}_{1},{\bf y}_{2}):=f_{1}({\bf y}_{1})f_{2}({\bf y}_{2})$. We note that 
$\{ f=f_{1}({\bf y}_{1})f_{2}({\bf y}_{2}) \vert~f_{1},f_{2}\in C_{\infty}(\Gamma \otimes \mathbb R) \}
\subset C_{b}((\Gamma \otimes \mathbb R)^{2})$ is a {\it{determining class}} 
on $(\Gamma \otimes \mathbb R)^{2}$.
(See e.g. Karatzas--Shreve \cite[Definition 5.4.24]{Karatzas Shreve} for the precise meaning. 
In Klenke \cite[Definition 13.9]{Klenke}, it is
called a {\it{separating family}}.) 

We aim to show
\begin{eqnarray}
\lefteqn{
\lim_{n\to \infty} \int_{(\Gamma \otimes \mathbb R)^{2}} 
f({\bf y}_{1},{\bf y}_{2}) \big(
{\mathbb P}_{x_{*}}\circ
({\bf X}^{(n)}_{s,t})^{-1} \big)(d{\bf y}_{1}d{\bf y}_{2})}
\nonumber \\
&=& \int_{(\Gamma \otimes \mathbb R)^{2}}
G_{s}({\bf y}_{1})G_{t-s}({\bf y}_{2}-{\bf y}_{1}) f({\bf y}_{1}, {\bf y}_{2})
d{\bf y}_{1}d{\bf y}_{2}.
\label{Conv-Jan1}
\end{eqnarray}
We prepare 
\begin{align}
 \sup_{c\in {\Omega}_{x_{*}}(X)}
 \vert {\bf X}_{t}^{(n)} (c) - {\cal X}_{t}^{(n)} (c) \vert_{g_{0}}
&=
\frac{(nt -[nt])}{\sqrt{n}}  \sup_{c\in {\Omega}_{x_{*}}(X)}
\vert 
\xi_{[nt]+1}(c)-\xi_{[nt]}(c) -\rho_\mathbb{R} (\gamma_p)  \vert_{g_{0}} 
\nonumber \\
&\leq  \frac{1}{\sqrt{n}} \big(  \Vert d\Phi_{0} \Vert_{\infty}+
\vert \rho_\mathbb{R} (\gamma_p) \vert_{g_{0}} \big). 
\label{Difference}
\end{align}
Then, by virtue of uniform continuity of $f$ on 
$(\Gamma \otimes \mathbb{R})^{2}$ and (\ref{Difference}),
we obtain
\begin{equation}
\lim_{n\to \infty} \sum_{c\in {\Omega}_{x_{*}}(X)} \big \vert f({\bf X}_{s,t}^{(n)} (c))
-f({\cal X}_{s,t}^{(n)} (c))
\big \vert {\mathbb P}_{x_{*}} (\{c\}) =0.
\nonumber
\end{equation}

Moreover, we have
\begin{eqnarray}
\lefteqn{ 
\sum_{c\in {\Omega}_{x_{*}}(X)}
 f({\cal X}^{(n)}_{s,t}(c)) {\mathbb P}_{x_{*}}(\{c\})
}
\nonumber \\
&=& \sum_{c_{1} \in {\Omega}_{x_{*},[ns]}(X)}p(c_{1})f_{1}
\Big(n^{-1/2} \big( \Phi_{0}(t(c_{1}))-[ns] \rho_\mathbb{R} (\gamma_p)\big) \Big)
\nonumber \\
&\mbox{ }& \times
\Big \{ \sum_{c_{2}\in {\Omega}_{t(c_{1}), [nt]-[ns]}(X)} 
p(c_{2}) f_{2}
\Big(n^{-1/2} \big( \Phi_{0}(t(c_{2}))-[nt] \rho_\mathbb{R} (\gamma_p)\big) \Big) \Big \}
\nonumber \\
&=& \sum_{c_{1} \in {\Omega}_{x_{*},[ns]}(X)}p(c_{1}) ({\cal P}_{n^{-1/2}}f_{1})
\big( t(c_{1}), [ns]
\rho_\mathbb{R} (\gamma_p) \big)
{\cal L}_{\gamma_{p}}^{ [nt]-[ns]} ({\cal P}_{n^{-1/2}} f_{2}) 
\big( t(c_{1}), [ns]\rho_\mathbb{R} (\gamma_p) \big)
\nonumber \\
&=& {\cal L}_{\gamma_{p}}^{ [ns]} 
\big\{  ({\cal P}_{n^{-1/2}} f_{1}) \cdot {\cal L}_{\gamma_{p}}^{ [nt]-[ns]}({\cal P}_{n^{-1/2}} f_{2}) \big \} 
(x_{*}, {\bf 0}).
\nonumber
\end{eqnarray}
Then Theorem \ref{CLT-1} implies 
\begin{eqnarray}
& &
\hspace{-15mm} 
\big \Vert {\cal L}_{\gamma_{p}}^{ [ns]} \big\{  ({\cal P}_{n^{-1/2}} f_{1}) 
\cdot {\cal L}_{\gamma_{p}}^{ [nt]-[ns]}({\cal P}_{n^{-1/2}} f_{2}) \big \} 
-
{\cal P}_{n^{-1/2}} e^{-\frac{s}{2}\Delta} \big ( f_{1}e^{-\frac{t-s}{2}\Delta}f_{2} \big ) 
\big \Vert_{\infty}
\nonumber \\
& & 
%
\leq 
\big \Vert {\cal L}_{\gamma_{p}}^{ [ns]} \big\{  ({\cal P}_{n^{-1/2}} f_{1}) 
\cdot {\cal L}_{\gamma_{p}}^{ [nt]-[ns]}({\cal P}_{n^{-1/2}} f_{2}) \big \} 
\nonumber \\
& &
\hspace{20mm}
-
{\cal L}_{\gamma_{p}}^{ [ns]} 
\big \{ ( {\cal P}_{n^{-1/2}} f_{1} )\cdot {\cal P}_{n^{-1/2}} \big( e^{-\frac{t-s}{2}\Delta}f_{2} \big) \big \} 
\big \Vert_{\infty}
\nonumber \\
& &
%
\quad 
+
\Vert {\cal L}_{\gamma_{p}}^{ [ns]} {\cal P}_{n^{-1/2}} \big ( f_{1}e^{-\frac{t-s}{2}\Delta}f_{2} \big )
-{\cal P}_{n^{-1/2}} e^{-\frac{s}{2}\Delta} \big ( f_{1}e^{-\frac{t-s}{2}\Delta}f_{2} \big ) \Vert_{\infty}
\nonumber \\
& &
\leq 
\Vert f_{1} \Vert_{\infty} \Vert {\cal L}_{\gamma_{p}}^{ [nt]-[ns]}({\cal P}_{n^{-1/2}} f_{2}) 
-
{\cal P}_{n^{-1/2}} \big( e^{-\frac{t-s}{2}\Delta}f_{2} \big) \Vert_{\infty} 
\nonumber \\
& &
\quad
+
\Vert {\cal L}_{\gamma_{p}}^{ [ns]} {\cal P}_{n^{-1/2}} \big ( f_{1}e^{-\frac{t-s}{2}\Delta}f_{2} \big )
-{\cal P}_{n^{-1/2}} e^{-\frac{s}{2}\Delta} \big ( f_{1}e^{-\frac{t-s}{2}\Delta}f_{2} \big ) \Vert_{\infty}
\nonumber \\
& &
\to 0
\qquad \mbox{ as } n\to \infty.
\nonumber
\end{eqnarray}
Summarizing all the above, we have
\begin{eqnarray}
\lefteqn{
\lim_{n\to \infty} \int_{(\Gamma \otimes \mathbb R)^{2}} f({\bf y}_{1},{\bf y}_{2}) 
\big(
{\mathbb P}_{x_{*}}\circ
({\bf X}^{(n)}_{s,t})^{-1} \big)(d{\bf y}_{1}d{\bf y}_{2})} 
\nonumber \\
&=& \lim_{n \to \infty}\sum_{c\in {\Omega}_{x_{*}}(X)}
 f({\cal X}^{(n)}_{s,t}(c)) {\mathbb P}_{x_{*}}(\{c\})
 \nonumber \\
&=& \lim_{n \to \infty} 
{\cal L}_{\gamma_p}^{ [ns]} \big \{  {\cal P}_{n^{-1/2}} f_{1} \cdot {\cal P}_{n^{-1/2}} 
\big( e^{-\frac{t-s}{2}\Delta}f_{2} \big) \big \} ( x_{*}, {\bf 0})
\nonumber \\
&=& \lim_{n\to \infty} {\cal L}_{\gamma_p}^{ [ns]} {\cal P}_{n^{-1/2}} 
\big ( f_{1}e^{-\frac{t-s}{2}\Delta}f_{2} \big )(x_{*},{\bf 0})
\nonumber \\
&=& \lim_{n \to \infty} {\cal P}_{n^{-1/2}} e^{-\frac{s}{2}\Delta} 
\big ( f_{1}e^{-\frac{t-s}{2}\Delta}f_{2} \big )(x_{*}, {\bf 0})
\nonumber \\
&= & e^{-\frac{s}{2}\Delta} \big( f_{1} e^{-\frac{t-s}{2}\Delta}f_{2} \big)({\bf 0})
\nonumber \\
&=&\int_{\Gamma \otimes \mathbb R} d{\bf y}_{1}
G_{s}({\bf y}_{1}) f_{1}({\bf y}_{1}) \int_{\Gamma \otimes \mathbb R} d{\bf y}_{2}
G_{t-s}({\bf y}_{2}-{\bf y}_{1}) f_{2}({\bf y}_{2})
\nonumber \\
&=& \int_{(\Gamma \otimes \mathbb R)^{2}} 
f({\bf y}_{1},{\bf y}_{2}) G_{s}({\bf y}_{1})G_{t-s}({\bf y}_{2}-{\bf y}_{1}) d{\bf y}_{1}d{\bf y}_{2}.
\nonumber
\end{eqnarray}
Furthermore
we have
 \begin{eqnarray}
{\mathbb P}_{x_{*}} \big( \vert {\bf X}_{s,t}^{(n)} 
 \vert_{(\Gamma \otimes \mathbb R)^{2}} >R \big)
 &\leq & {\mathbb P}_{x_{*}} \Big( \vert {\bf X}_{s}^{(n)} 
 \vert_{g_{0}} >\frac{R}{{\sqrt 2}} \Big)
 +{\mathbb P}_{x_{*}} \Big( \vert {\bf X}_{t}^{(n)} \vert_{g_{0}} 
 >\frac{R}{{\sqrt 2}} \Big)
 \nonumber \\
 &\leq & \frac{4}{R^{4}} \Big( {\mathbb E}^{{\mathbb P}_{x_{*}}} \big [ 
 \vert {\bf X}_{s}^{(n)} \vert_{g_{0}}^{4} \big ]
 + {\mathbb E}^{{\mathbb P}_{x_{*}}} \big [ \vert {\bf X}_{t}^{(n)} 
 \vert_{g_{0}}^{4} \big ] \Big)
 \nonumber \\
 &\leq  &
 C R^{-4} (s^{2}+t^{2}),
 \label{chebyshev}
\nonumber
 \end{eqnarray}
where we used  (\ref{tight-4moment}) 
below and Chebyshev's inequality. This estimate implies tightness of 
$\{ {\mathbb P}_{x_{0}} \circ ({\cal X}_{s,t}^{(n)})^{-1} \}_{n=1}^{\infty}$
in probability measures on $((\Gamma \otimes \mathbb R)^{2},
{\cal B}((\Gamma \otimes \mathbb R)^{2}))$.

Finally, applying \cite[Theorem 13.16]{Klenke}, we have
the desired convergence (\ref{Conv-Jan1})  for 
every $f\in C_{b}((\Gamma \otimes \mathbb R)^{2})$. 
This completes the proof.
\qed
\vspace{2mm}

Having shown the convergence of finite dimensional distributions, 
we complete the proof of Theorem \ref{FCLT-1} by showing the following lemma:
\begin{lm} \label{tight1}
$\{ {\bf P}^{(n)} \}_{n=1}^{\infty}$ 
is tight in $({\bf W}, {\cal B}({\bf W}))$.
\end{lm}
{\bf Proof.}~From \cite[Problem 4.11]{Karatzas Shreve} (see also 
\cite[Theorem 21.42]{Klenke}), it is sufficient to show 
that there exists some $C>0$ independent of $n$ such that 
\begin{equation}
\mathbb{E}^{\mathbb P_{x_*}} \big [ 
\big | {\bf X}_t^{(n)}- {\bf X}_s^{(n)} \big |_{g_{0}}^4 \big] 
\leq C \vert t-s \vert^2 \qquad (0\leq s \leq t,~ n\in \mathbb N). 
\label{tight-4moment}
\end{equation}
We distinguish two cases:
$$ {\mbox{ {\bf{(I)}}  }}~ t-s<n^{-1}, \qquad {\mbox{ {\bf{(II)}} }}~t-s \geq n^{-1}. $$

First, we consider case {\bf{(I)}}. If $ns \geq [nt]$, 
\begin{equation}
\vert {\bf X}_t^{(n)}- {\bf X}_s^{(n)} \vert_{g_{0}}
\leq n^{1/2}(t-s) \big \{ \vert \xi_{[nt]+1}-\xi_{[nt]}\vert_{g_{0}}  
+ \vert \rho_{\mathbb R}(\gamma_{p}) 
\vert_{g_{0}} \big \}.
\nonumber
\end{equation}
On the other hand, if $ns <[nt] $, 
\begin{eqnarray}
\lefteqn{
\vert {\bf  X}_t^{(n)}- {\bf X}_s^{(n)} \vert_{g_{0}}}
\nonumber \\
& \leq & \Big \vert \frac{1-(ns-[ns])}{n^{1/2}} (\xi_{[nt]}-\xi_{[nt]-1})
+\frac{nt-[nt]}{n^{1/2}} (\xi_{[nt]+1}-\xi_{[nt]})
-n^{1/2}(t-s) {\rho_{\mathbb R}}(\gamma_{p}) \Big \vert_{g_{0}}
\nonumber \\
& \leq & n^{1/2} (t-s)
\big \{ \vert \xi_{[nt]+1}-\xi_{[nt]}\vert_{g_{0}}  + 
\vert \xi_{[nt]}-\xi_{[nt]-1}\vert_{g_{0}} +
\vert \rho_{\mathbb R}(\gamma_{p}) 
\vert_{g_{0}} \big \},
\nonumber
\end{eqnarray}
where we used $[ns]=[nt]-1$ and $ns \leq [nt] \leq nt$ for the third line.

In either case, we have
\begin{eqnarray}
  \vert {\bf X}_t^{(n)}- {\bf X}_s^{(n)} \vert_{g_{0}}
  &\leq &
 n^{1/2} (t-s)
\big \{ \vert \xi_{[nt]+1}-\xi_{[nt]}\vert_{g_{0}}  + 
\vert \xi_{[nt]}-\xi_{[nt]-1}\vert_{g_{0}} +
\vert \rho_{\mathbb R}(\gamma_{p}) 
\vert_{g_{0}} \big \}
\nonumber \\
&\leq & n^{1/2}(t-s) (2 \Vert d\Phi_{0} \Vert_{\infty} +\vert \rho_{\mathbb R}(\gamma_{p}) 
\vert_{g_{0}} ).
\label{Xn-est1}
\end{eqnarray}
Now we recall $n^{2}(t-s)^{2}<1$. Then (\ref{Xn-est1}) implies
\begin{equation}
\mathbb{E}^{\mathbb P_{x_{*}}} \big [ 
\big | {\bf X}_t^{(n)}- {\bf X}_s^{(n)} \big |_{g_{0}}^4 \big] 
\leq
(2 \Vert d\Phi_{0} \Vert_{\infty} +\vert \rho_{\mathbb R}(\gamma_{p}) 
\vert_{g_{0}} )^{4} (t-s)^{2}.
\nonumber
\end{equation}
Hence, we obtain the desired estimate (\ref{tight-4moment}) for case {\bf{(I)}}.

Next, we consider case {\bf{(II)}}. Let $\mathcal{F}$ be the fundamental domain in 
$X$ containing $x_{*} \in V$ and define 
${\mathfrak M}_{i}^{l}={\mathfrak M}_{i}^{l}(\Phi_{0}, \rho_{\mathbb R}(\gamma_{p}))$ 
($i=1, \ldots, d, l=1,2,3,4$) by
$$ 
{\mathfrak M}_{i}^{l}(x):=
\sum_{e \in E_{x}}  
p(e) \big( d\Phi_{0}(e)
-\rho_{\mathbb R}(\gamma_{p}) \big)^{l}_{i}
\qquad  (x\in V). $$
Note that ${\mathfrak M}_{i}^{l}$
is $\Gamma$-invariant and
$$ \Vert  {\mathfrak M}_{i}^{l}
\Vert_{\infty} \leq ( \Vert d\Phi_{0} \Vert_{\infty} +
\vert \rho_{\mathbb R}(\gamma_{p}) \vert_{g_{0}} )^{l} 
\qquad (i=1, \ldots, d).$$
Furthermore, the modified harmonicity (\ref{Phi-mh}) yields
${\mathfrak M}_{i}^{1}\equiv 0$~($i=1,\ldots, d$).
\vspace{2mm} 

We start by giving a bound on 
$\mathbb{E}^{\mathbb P_{x_*}}  \big [ 
\big | {\cal X}_{\frac{M}{n}}^{(n)}- 
{\cal X}_{\frac{N}{n}}^{(n)} \big |_{g_{0}}^4 \big ]$ 
($n\in \mathbb N, M\geq N\in \mathbb N$).
Applying an elementary inequality 
$(a_{1}+\cdots +a_{d})^{2} \leq d(a_{1}^{2}+\cdots +a_{d}^{2})$, we have
\begin{eqnarray}
\mathbb{E}^{\mathbb P_{x_*}}  \big [ 
\big | {\cal X}_{\frac{M}{n}}^{(n)}
- {\cal X}_{\frac{N}{n}}^{(n)} \big |_{g_{0}}^4 \big ]
&\leq & d n^{-2} 
\sum_{i=1}^{d}
\mathbb{E}^{\mathbb{P}_{x_{*}}} \big [ 
\big (\xi_{M}-\xi_{N}
-(M-N)\rho_{\mathbb R}(\gamma_{p}) \big )_{i}^{4}
\big ]
\nonumber \\
&=& d^{2}n^{-2} \max_{i=1,\ldots, d}
\sum_{c_{1} \in {\Omega}_{x_{*}, N}(X)} p(c_{1})
\sum_{c_{2} \in {\Omega}_{t(c_{1}), M-N}(X)} p(c_{2})
\nonumber \\
&\mbox{ }& \times
\big ( \Phi_{0} \big (t(c_{2}) \big)-\Phi_{0} \big (t(c_{1}) \big )-(M-N)\rho_{\mathbb R}(\gamma_{p})
 \big )_{i}^{4}
 \nonumber \\
&=& d^{2}n^{-2} \max_{i=1,\ldots, d} \max_{x\in {\cal F}}
\Big \{ \sum_{c \in {\Omega}_{x, M-N}(X)} p(c)
\nonumber \\
&\mbox{ }&
\times
\big ( \Phi_{0} \big (t(c) \big)-\Phi_{0} (x) -(M-N)\rho_{\mathbb R}(\gamma_{p})
 \big )_{i}^{4}
 \Big \}.
\label{Est-June17-0}
\end{eqnarray}
%
We now fix $i=1,\ldots, d$ and $x\in V$. 
For $k=1,\ldots, M-N$,  we have
\begin{eqnarray}
\lefteqn{ \sum_{c\in {\Omega}_{x,k}(X)}p(c) \big(
\Phi_{0}(t(c))-\Phi_{0} (x)-k \rho_{\mathbb R}(\gamma_{p}) \big )_{i}^{4} }
\nonumber \\
&=& \sum_{c'\in {\Omega}_{x,k-1}(X)} p(c') \sum_{e\in E_{t(c')}} p(e) \Big \{ \big( 
\Phi_{0}(t(e))-\Phi_{0} (o(e))-\rho_{\mathbb R}(\gamma_{p}) \big )_{i} 
\nonumber \\
&\mbox{ }& \quad + \big( 
\Phi_{0}(o(e))-\Phi_{0} (x)-(k-1)\rho_{\mathbb R}(\gamma_{p}) \big )_{i}  \Big \}^{4}
\nonumber \\
&=& \sum_{c'\in {\Omega}_{x,k-1}(X)} p(c') 
{\mathfrak M}_{i}^{4}(t(c'))
\nonumber \\
&\mbox{ }&
+4\sum_{c'\in {\Omega}_{x,k-1}(X)} p(c') 
 \big ( \Phi_{0}(t(c'))-\Phi_{0} (x)-(k-1)\rho_{\mathbb R}(\gamma_{p}) \big)_{i}  
{\mathfrak M}_{i}^{3}(t(c'))
\nonumber \\
&\mbox{ }&
+6\sum_{c'\in {\Omega}_{x,k-1}(X)} p(c') 
 \big ( \Phi_{0}(t(c'))-\Phi_{0} (x)-(k-1)\rho_{\mathbb R}(\gamma_{p}) \big)_{i}^{2}  
 {\mathfrak M}_{i}^{2} (t(c'))
 \nonumber \\
 &\mbox{ }&
+4\sum_{c'\in {\Omega}_{x,k-1}(X)} p(c') 
 \big ( \Phi_{0}(t(c'))-\Phi_{0} (x)-(k-1)\rho_{\mathbb R}(\gamma_{p}) \big)_{i}^{3}  
 {\mathfrak M}_{i}^{1} (t(c'))
 \nonumber \\
 &\mbox{ }&
 +\sum_{c'\in {\Omega}_{x,k-1}(X)}p(c') \big(
\Phi_{0}(t(c'))-\Phi_{0} (x)-(k-1) \rho_{\mathbb R}(\gamma_{p}) \big )_{i}^{4}
\nonumber \\
&\leq &
\sum_{c\in {\Omega}_{x,k-1}(X)}p(c) \big(
\Phi_{0}(t(c))-\Phi_{0} (x)-(k-1) \rho_{\mathbb R}(\gamma_{p}) \big )_{i}^{4}
\nonumber \\
&\mbox{ }&
+\Vert  {\mathfrak M}_{i}^{4}
\Vert_{\infty}
+4 \big \{ \Vert d\Phi_{0} \Vert_{\infty}+(k-1) \vert \rho_{\mathbb R}(\gamma_{p}) 
\vert_{g_{0}} \big \}
\Vert  {\mathfrak M}_{i}^{3}
\Vert_{\infty}
\nonumber \\
&\mbox{ }&
+6\Vert  {\mathfrak M}_{i}^{2}
\Vert_{\infty}
\sum_{c\in {\Omega}_{x,k-1}(X)} p(c) 
 \big ( \Phi_{0}(t(c))-\Phi_{0} (x)-(k-1)\rho_{\mathbb R}(\gamma_{p}) \big)_{i}^{2} . 
 \label{Est-June17-1}
\end{eqnarray}
Furthermore, it follows from the modified harmonicity  (\ref{Phi-mh}) that
\begin{eqnarray}
\lefteqn{ \sum_{c\in {\Omega}_{x,k-1}(X)} p(c) 
 \big ( \Phi_{0}(t(c))-\Phi_{0} (x)-(k-1)\rho_{\mathbb R}(\gamma_{p}) \big)_{i}^{2}  }
 \nonumber \\
 &=& 
 \sum_{c'\in {\Omega}_{x,k-2}(X)} p(c') \sum_{e\in E_{t(c')}} p(e) \Big \{ \big( 
\Phi_{0}(t(e))-\Phi_{0} (o(e))-\rho_{\mathbb R}(\gamma_{p}) \big )_{i} 
\nonumber \\
&\mbox{ }& \quad + \big( 
\Phi_{0}(o(e))-\Phi_{0} (x)-(k-1)\rho_{\mathbb R}(\gamma_{p}) \big )_{i}  \Big \}^{2}
 \nonumber \\
&=&
\sum_{c'\in {\Omega}_{x,k-2}(X)} p(c') 
\Big \{
 \big ( \Phi_{0}(t(c'))-\Phi_{0} (x)-(k-2)\rho_{\mathbb R}(\gamma_{p}) \big)_{i}^{2} 
+ {\mathfrak M}_{i}^{2}(t(c'))
 \Big \}
 \nonumber \\
 &\leq &
 \sum_{c'\in {\Omega}_{x,k-2}(X)} p(c') 
 \big ( \Phi_{0}(t(c'))-\Phi_{0} (x)-(k-2)\rho_{\mathbb R}(\gamma_{p}) \big)_{i}^{2} 
 + \Vert  {\mathfrak M}_{i}^{2} \Vert_{\infty}
\nonumber \\ 
 & \leq & (k-1)  \Vert  {\mathfrak M}_{i}^{2} \Vert_{\infty}.
 \label{Est-June17-2}
 \end{eqnarray}
Combining (\ref{Est-June17-1}) with  (\ref{Est-June17-2}), we obtain 
\begin{eqnarray}
\lefteqn{
 \sum_{c\in {\Omega}_{x,k}(X)}p(c) \big(
\Phi_{0}(t(c))-\Phi_{0} (x)-k \rho_{\mathbb R}(\gamma_{p}) \big )_{i}^{4}
}
\nonumber \\
& \leq & 
(k-1) \Big( 4 \vert \rho_{\mathbb R}(\gamma_{p}) 
\vert_{\Gamma \otimes \mathbb R}
 \Vert  {\mathfrak M}_{i}^{3} \Vert_{\infty}
 +6\Vert  {\mathfrak M}_{i}^{2}  \Vert_{\infty}^{2} \Big)
 +\big( \Vert  {\mathfrak M}_{i}^{4}  \Vert_{\infty}
 +6 \Vert d\Phi_{0} \Vert_{\infty}
  \Vert  {\mathfrak M}_{i}^{3} \Vert_{\infty} \big)=Ck,
\nonumber
\end{eqnarray}
and it implies 
\begin{equation}
 \sum_{c\in {\Omega}_{x,M-N}(X)}p(c) \big(
\Phi_{0}(t(c))-\Phi_{0} (x)-(M-N) \rho_{\mathbb R}(\gamma_{p}) \big )_{i}^{4}
\leq 
C (M-N)^{2}.
\label{Est-June17-3}
\end{equation}
 
Therefore, by combining (\ref{Est-June17-0}) with (\ref{Est-June17-3}), we obtain
\begin{eqnarray}
\mathbb{E}^{\mathbb P_{x_*}} \big [ 
\big | {\bf X}_t^{(n)}- {\bf X}_s^{(n)} \big |_{g_{0}}^4 \big]
& \leq &
\mathbb{E}^{\mathbb P_{x_*}} \Big [  \max_{i,j=0,1}
\big | {\cal X}_{\frac{[nt]+i}{n}}^{(n)}- {\cal X}_{\frac{[ns]+j}{n}}^{(n)}
\big |_{g_{0}}^4 \Big]
\nonumber \\
& \leq & Cd^{2}n^{-2} ( [nt]-[ns]+1)^{2}
\nonumber \\
&\leq & C d^{2} \big( t-s+\frac{2}{n})^{2}
\nonumber \\
& \leq & C d^{2} \big \{ (t-s)+2(t-s) \big \}^{2}=C (t-s) ^{2},
\nonumber
\end{eqnarray}
where we used $[nt]-[ns] \leq n(t-s)+1$ for the third line and the condition 
$n^{-1}\leq (t-s)$ for the final line.
Thus we obtained our desired estimate
(\ref{tight-4moment}) for case {\bf{(II)}}.
\qed
\section{Proof of the CLT of the second kind}
\label{weak convergence 2}
\subsection{Basic properties of the family of Albanese metrics $\{ g_{0}^{(\varepsilon)} \}$ }
Let $\{ p_{\varepsilon} \}_{0\leq \varepsilon \leq 1}$ be a family 
of transition probabilities defined by (\ref{def-pe}). Thanks to 
\cite[Proposition 2.3]{KS06}, we easily obtain the following lemma:
\begin{lm} \label{pdelta}
{\rm{(1)}}~
$p_{1}=p$ and 
$\gamma_{p_{\varepsilon}}= \varepsilon \gamma_{p}$ for every
$0\leq \varepsilon \leq 1$.
\\
{\rm{(2)}}~$p_{\varepsilon}(e)>0$ for every $e\in E_{0}$ provided $0\leq \varepsilon <1$.
\\
{\rm{(3)}}~For every $0\leq \varepsilon \leq 1$,
the invariant measure of the random walk given by
$p_\varepsilon$ is $m$.
Moreover, $p_0$ and $q$ are {\rm{(}}$m$-{\rm{)}}symmetric and {\rm{(}}$m$-{\rm{)}}anti-symmetric, 
respectively.
Specifically,
\begin{equation*}
p_0(e)m(o(e))=p_0(\overline{e})m(t(e)),~ q(e)m(o(e))=- q(\overline{e}) m(t(e))
\qquad (e\in E_{0}).
\end{equation*}
{\rm{(4)}}~
\begin{equation}
 \sum_{e\in E_{0}} q(e) \omega_{1}(e) \omega_{2}(e) m(o(e))=0
\qquad (\omega_{1}, \omega_{2}
\in C^{1}(X_{0}, \mathbb R)).
\label{q-anti-1}
\end{equation}
\end{lm}
We put  ${\widetilde m}_{\varepsilon}(e):=p_{\varepsilon}(e)m(o(e))$ 
($0\leq \varepsilon \leq 1$, $e\in E_{0}$)
and denote by ${\cal H}_{(\varepsilon)}^{1}(X_0)$ 
the set of modified harmonic $1$-forms with respect to the 
transition probability $p_{\varepsilon}$. Namely, 
${\cal H}_{(\varepsilon)}^{1}(X_0)$ is the set of $1$-forms on $X_{0}$ satisfying
\begin{equation}
(\delta_{p_{\varepsilon}} \omega)(x)+ \langle \gamma_{p_{\varepsilon}}, \omega \rangle=0
\qquad (x\in V_{0}).
\nonumber
\end{equation} 
We equip ${\cal H}_{(\varepsilon)}^{1}(X_0)$ with the inner product
\begin{align}
\langle \hspace{-0.8mm} \langle \omega_{1}, \omega_{2} 
\rangle \hspace{-0.8mm} \rangle_{(\varepsilon)} 
&:= \sum_{e\in E_{0}} \omega_{1}(e) \omega_{2}(e) 
{\widetilde m}_{\varepsilon}(e)-\langle \gamma_{p_{\varepsilon}}, \omega_{1} \rangle
\langle \gamma_{p_{\varepsilon}}, \omega_{2} \rangle
\nonumber \\
&=
 \sum_{e\in E_{0}} \omega_{1}(e) \omega_{2}(e) 
{\widetilde m}_{\varepsilon}(e)-\varepsilon^{2} \langle \gamma_{p}, \omega_{1} \rangle
\langle \gamma_{p}, \omega_{2} \rangle
\qquad 
 (\omega_{1}, \omega_{2} \in {\cal H}^{1}_{(\varepsilon)}(X_{0}) ),
\nonumber
\end{align}
The corresponding norm $\Vert \cdot \Vert_{(\varepsilon)}$ is given by
\begin{equation} 
\Vert \omega \Vert_{(\varepsilon)}^{2}:=
\langle \hspace{-0.8mm} \langle \omega, \omega \rangle \hspace{-0.8mm} \rangle_{(\varepsilon)}
=\sum_{e\in E_{0}} \vert \omega(e) \vert^{2}{\widetilde{m}_{\varepsilon}}(e)
-\varepsilon^{2} \langle \gamma_{p}, \omega \rangle^{2} \qquad (\omega \in
{\cal H}^{1}_{(\varepsilon)}(X_{0}) ).
\nonumber
\end{equation}

By the discrete Hodge--Kodaira theorem mentioned in Section 3,
we may identify ${\rm H}^{1}(X_{0}, \mathbb R)$ with 
${\cal H}_{(\varepsilon)}^{1}(X_0)$ for each $0\leq \varepsilon \leq 1$.
(It should be noted that the identification map depends on
the parameter ${\varepsilon}$ and ${\cal H}_{(1)}^{1}(X_0)={\cal H}^{1}(X_0)$.)
Moreover, we also identify 
$\mathrm{Hom}(\Gamma, \mathbb{R})$ with 
${\rm Image}(\hspace{-1mm}\mbox{ }^{t}\rho_{\mathbb R})\subset
{\rm H}^{1}(X_{0}, \mathbb R)$.
Hence we may regard
$\mathrm{Hom}(\Gamma, \mathbb{R})$ as a subspace of 
${\cal H}_{(\varepsilon)}^{1}(X_0)$.
For $\omega \in \mathrm{Hom}(\Gamma, \mathbb{R})$, 
we denote
$\hspace{-1mm}\mbox{ }^{t}\rho_{\mathbb R}(\omega) \in {\rm H}^{1}(X_{0}, \mathbb R)
\cong {\cal H}_{(\varepsilon)}^{1}(X_0)$
by $\omega^{(\varepsilon)}$.
Then we have
\begin{lm}\label{Harmonic-pdelta}
\begin{equation*}
\lim_{\varepsilon \searrow 0}
\langle \hspace{-0.5mm} \langle 
\omega^{(\varepsilon)}, 
\eta^{(\varepsilon)} 
\rangle \hspace{-0.5mm} \rangle_{(\varepsilon)}
= \langle \hspace{-0.5mm} \langle 
\omega^{(0)}, \eta^{(0)} 
\rangle \hspace{-0.5mm} \rangle_{(0)}
\qquad (\omega, \eta \in \mathrm{Hom}(\Gamma, \mathbb{R}) ). 
\end{equation*}
\end{lm}
%
{\bf Proof.}~
By following the proof of   
\cite[Lemma 5.2]{KS06}, we observe that
for any $0\leq \varepsilon \leq 1$, 
there exist functions 
$f^{(\varepsilon)}, g^{(\varepsilon)}$ defined on $X_0$ such that 
$$\omega^{(\varepsilon)}=\omega^{(0)} +df^{(\varepsilon)},~ 
\eta^{(\varepsilon)}=
\eta^{(0)} +dg^{(\varepsilon)} (\in {\cal H}_{(\varepsilon)}^1(X_0) ).$$
Since $\omega^{(0)} \in {\cal H}_{(0)}^{1}(X_0)$, we have
$\delta_{p_{0}} \omega^{(0)}=0$.  
Hence
\begin{equation}
\sum_{e\in E_{0}} p_{0}(e) \omega^{(0)}(e) df(e) m(o(e))
=2 \langle \delta_{p_{0}} \omega^{(0)}, fm \rangle_{0}=0
\qquad (f\in C^{0}(X_{0}, \mathbb R) ).
\label{p0-df}
\end{equation}
Combining (\ref{q-anti-1}), (\ref{p0-df}) with the identity
$\langle \gamma_{p}, df\rangle=0$ ($f\in C^{0}(X_{0}, \mathbb R)$),
we can expand 
\begin{eqnarray}
\langle \hspace{-0.5mm} \langle 
\omega^{(\varepsilon)}, \eta^{(\varepsilon)} \rangle  \hspace{-0.5mm} \rangle_{(\varepsilon)}
&=& \sum_{e \in E_0} p_\varepsilon (e)
\big ( \omega^{(0)}(e) + df^{(\varepsilon)}(e) \big)  \big(\eta^{(0)} (e) +dg^{(\varepsilon)}(e) \big) 
m(o(e))
\nonumber \\
&\mbox{ }&-
\varepsilon^2 \langle \gamma_p, \omega^{(0)}+df^{(\varepsilon)} \rangle
\langle \gamma_p, \eta^{(0)}+dg^{(\varepsilon)} \rangle 
\nonumber \\
&=&
\langle \hspace{-0.5mm} \langle \omega^{(0)}, \eta^{(0)} 
\rangle \hspace{-0.5mm}\rangle_{(0)} +
\langle \hspace{-0.5mm} \langle df^{(\varepsilon)}, dg^{(\varepsilon)} 
\rangle \hspace{-0.5mm} \rangle_{(0)}
-
\varepsilon^2 \langle \gamma_p, \omega^{(0)} \rangle
\langle \gamma_p, \eta^{(0)}  \rangle. 
\nonumber
\end{eqnarray}
Hence it suffices to show that 
$\langle \hspace{-0.5mm} \langle 
df^{(\varepsilon)}, df^{(\varepsilon)} \rangle \hspace{-0.5mm}\rangle_{(0)}$ 
converges to $0$ as $\varepsilon \searrow 0$. 

We now define 
the operator $Q: \ell^{2} (X_{0}) \to \ell^{2}(X_{0})$
by $$Qf(x):= \sum_{e\in (E_{0})_{x}} q(e) f(t(e))
\qquad (x\in V_{0}).$$
Then the transition operator $L_{(\varepsilon)}: \ell^{2}(X_{0})
\to 
\ell^{2}(X_{0})$ associated with the transition
probability $p_{\varepsilon}$ is decomposed by
$L_{(\varepsilon)}=L_{(0)}+\varepsilon Q$.
We also note that
anti-symmetry of $q$  implies 
\begin{equation}
\langle Qf, fm \rangle_{\ell^{2}(X_{0})}
=
\langle Qf,fm \rangle_{0}=0 \qquad (f\in  \ell^{2}(X_{0}) ).
\label{Q-int}
\end{equation}
Then (\ref{Q-int}) and the identity $I-L_{(\varepsilon)}=\delta_{p_{\varepsilon}}d$
imply
\begin{eqnarray}
\langle \hspace{-0.8mm} \langle 
df^{(\varepsilon)}, df^{(\varepsilon)} \rangle \hspace{-0.8mm}\rangle_{(0)}
&=& 
\big \langle \hspace{-1mm} \big \langle 
d \big (f^{(\varepsilon)}-m(f^{(\varepsilon)}) \big ), 
d \big (f^{(\varepsilon)}-m(f^{(\varepsilon)}) \big )
\big \rangle \hspace{-1mm} \big \rangle_{(0)}
\nonumber \\
&=&
2 \big \langle (I-L_{(0)})   \big (f^{(\varepsilon)}-m(f^{(\varepsilon)}) \big ),
\big (f^{(\varepsilon)}-m(f^{(\varepsilon)}) \big ) m \big \rangle_{0}
\nonumber \\
&=& 2 \big \langle (I-L_{(\varepsilon)})f^{(\varepsilon)}, 
\big (f^{(\varepsilon)}-m(f^{(\varepsilon)}) \big )
m  \big \rangle_{0}
\nonumber \\
&=&
2 \big \langle 
\delta_{p_{\varepsilon}}(df^{(\varepsilon)}), 
\big (f^{(\varepsilon)}-m(f^{(\varepsilon)}) \big )
m \big \rangle_{0},
\label{quadratic-dfe}
\end{eqnarray}
where $$m(f^{(\varepsilon)}):=\sum_{x \in X_0} f^{(\varepsilon)}(x) m(x)$$ 
and we used $L_{(\varepsilon)}m(f^{(\varepsilon)})=m(f^{(\varepsilon)})$ for the
third line.

Recalling $\omega^{(0)}+df^{(\varepsilon)} \in  {\cal H}_{(\varepsilon)}^1(X_0)$ and 
$\delta_{p_{0}} \omega^{(0)}=0$,
we have
\begin{eqnarray}
\delta_{p_{\varepsilon}}(df^{(\varepsilon)})(x)
&=& -(\delta_{p_{\varepsilon}} \omega^{(0)})(x)-\varepsilon \langle \gamma_{p}, \omega^{(0)}
\rangle
\nonumber \\
&=& -\varepsilon \big( 
(\delta_{q} \omega^{(0)})(x)+ \langle \gamma_{p}, \omega^{(0)}
\rangle \big) \qquad (x\in V_{0} ).
\label{dpedfe}
\end{eqnarray}
Thus
\begin{equation}
f^{(\varepsilon)}(x)=L_{(\varepsilon)}f^{(\varepsilon)}(x)
-\varepsilon 
\big( 
(\delta_{q} \omega^{(0)})(x)+ \langle \gamma_{p}, \omega^{(0)}
\rangle \big) \qquad (x\in V_{0} ),
\label{fe-representation}
\end{equation}
and (\ref{fe-representation}) still holds if we replace $f^{(\varepsilon)}$ by
$f^{(\varepsilon)}-m(f^{(\varepsilon)})$. Namely,
\begin{align}
f^{(\varepsilon)}(x)-m(f^{(\varepsilon)})
&=L_{(\varepsilon)}\big( f^{(\varepsilon)}
-m(f^{(\varepsilon)}) \big) (x)
\nonumber \\
&\quad
-\varepsilon 
\big( 
(\delta_{q} \omega^{(0)})(x)+ \langle \gamma_{p}, \omega^{(0)}
\rangle \big) \qquad (x\in V_{0} ).
\label{2fe-representation}
\end{align}

We now recall $\ell^{2}(X_{0})=\langle \phi_{0} \rangle
\oplus \ell^{2}_{1}(X_{0})$, where $\phi_{0}$ and $\ell^{2}_{1}(X_{0})$
are introduced in the proof of Proposition \ref{Ergord-tm1}.
Note $f^{(\varepsilon)}-m(f^{(\varepsilon)})\in \ell^{2}_{1}(X_{0})$ and that
the transition operator $L_{(\varepsilon)}$
maps $\ell^{2}_{1}(X_{0})$ to  $\ell^{2}_{1}(X_{0})$ for all $0\leq \varepsilon \leq 1$. 
Furthermore, as $\alpha_{0}(\varepsilon)=1$ is a simple eigenvalue of 
$L_{(\varepsilon)}$ for all $0\leq \varepsilon \leq 1$,
the inverse operator of $(1-L_{(\varepsilon)}) \vert_{\ell^{2}_{1}(X_{0})}:
\ell^{2}_{1}(X_{0})  \to  \ell^{2}_{1}(X_{0}) $ exists. 
Because of $\delta_{q}\omega^{(0)}+\langle \gamma_{p}, \omega^{(0)} \rangle
\in \ell^{2}_{1}(X_{0})$, we can solve equation
(\ref{2fe-representation}) as 
\begin{equation} f^{(\varepsilon)}-m(f^{(\varepsilon)})
= -\varepsilon (1-L_{(\varepsilon)}) \vert_{\ell^{2}_{1}(X_{0})}^{-1} 
\big( \delta_{q} \omega^{(0)}+ \langle \gamma_{p}, \omega^{(0)}
\rangle \big).
\label{2fe-sol}
\end{equation}
Combining (\ref{2fe-sol}) with the identity
$$ 
(1-L_{(\varepsilon)}) \vert_{\ell^{2}_{1}(X_{0})}^{-1} 
=(1-L_{(0)}) \vert_{\ell^{2}_{1}(X_{0})}^{-1} 
\Big [ 1- \varepsilon Q \vert_{\ell^{2}_{1}(X_{0})}
(1-L_{(0)}) \vert_{\ell^{2}_{1}(X_{0})}^{-1}  \Big ]^{-1},
$$
we obtain
\begin{eqnarray}
\Vert f^{(\varepsilon)}-m(f^{(\varepsilon)}) \Vert_{\ell^{2}(X_{0})}
&\leq & \varepsilon
\big \Vert   (1-L_{(\varepsilon)}) \vert_{\ell^{2}_{1}(X_{0})}^{-1}  \big \Vert 
\Big(
\Vert \delta_{q}\omega^{(0)} \Vert_{\ell^{2}(X_{0})} +\vert
\langle \gamma_{p}, \omega^{(0)}
\rangle \vert
\Big)
\nonumber \\
&\leq & \varepsilon 
\big \Vert   (1-L_{(0)}) \vert_{\ell^{2}_{1}(X_{0})}^{-1} \big \Vert
\Big( 1-\varepsilon 
\big \Vert   
Q \vert_{\ell^{2}_{1}(X_{0})}
(1-L_{(0)}) \vert_{\ell^{2}_{1}(X_{0})}^{-1} 
\big \Vert
\Big)^{-1}
\nonumber \\
&\mbox{ }&
\times 
\Big(
\Vert \delta_{q}\omega^{(0)} \Vert_{\ell^{2}(X_{0})} +\vert
\langle \gamma_{p}, \omega^{(0)}
\rangle \vert
\Big).
\label{kiichi-est1}
\end{eqnarray}
Now we choose a sufficiently small constant $\varepsilon_{0}>0$ such that
$$ \sup_{0\leq \varepsilon \leq \varepsilon_{0}} 
\Big( 1-\varepsilon 
\big \Vert   
Q \vert_{\ell^{2}_{1}(X_{0})}
(1-L_{(0)}) \vert_{\ell^{2}_{1}(X_{0})}^{-1} 
\big \Vert
\Big)^{-1} \leq 2.
$$
Then (\ref{kiichi-est1}) implies
\begin{equation}
\sup_{0\leq \varepsilon \leq \varepsilon_{0}} 
\Vert f^{(\varepsilon)}-m(f^{(\varepsilon)}) \Vert_{\ell^{2}(X_{0})} <\infty.
\label{kiichi-est2}
\end{equation}

Finally, by combining (\ref{quadratic-dfe}), (\ref{dpedfe}) with  (\ref{kiichi-est2}),
we obtain $\lim_{\varepsilon \searrow 0} 
\langle \hspace{-0.8mm} \langle 
df^{(\varepsilon)}, df^{(\varepsilon)} \rangle \hspace{-0.8mm}\rangle_{(0)}=0$.
This completes the proof. \qed

We define the Albanese metric $g_{0}^{(\varepsilon)}$
on $\Gamma \otimes \mathbb R$ by
the dual metric of $\langle \hspace{-0.8mm} \langle \cdot, \cdot \rangle \hspace{-0.8mm} \rangle_{(\varepsilon)}$.
The following lemma on the family of Albanese metrics
$\{ g_{0}^{(\varepsilon)} \}_{0\leq \varepsilon \leq 1}$ plays a key role in the
proof of Theorems \ref{CLT-2} and \ref{FCLT-2}.
\begin{lm}
\label{dual-Harmonic-conv}
{\rm{(1)}}
\begin{equation*}
\lim_{\varepsilon \searrow 0}
\big \langle
{\bf{x}} , {\bf{y}}
\big \rangle_{g_{0}^{(\varepsilon)}}
= 
\big \langle
{\bf{x}} , {\bf{y}}
\big \rangle_{g_{0}^{(0)}}
\qquad 
({\bf{x}}, {\bf{y}} \in \Gamma \otimes \mathbb R) .
\end{equation*}
{\rm{(2)}}~Let $\Phi_{0}^{(\varepsilon)}: X \to 
\Gamma \otimes \mathbb R$ be the modified harmonic realization with respect to $p_{\varepsilon}$.
Then we have
\begin{equation*}
\lim_{\varepsilon \searrow 0}
\sum_{e\in E_{0}} 
\big \vert d\Phi_{0}^{(\varepsilon)}({\widetilde e})
\big \vert_{g_{0}^{(\varepsilon)}}^{2} 
{\widetilde m}_{\varepsilon} (e)
= 
\sum_{e\in E_{0}}
\big \vert d\Phi_{0}^{(0)}({\widetilde e})
\big \vert_{g_{0}^{(0)}}^{2} 
{\widetilde m}_{0}(e).
\end{equation*}
\\
{\rm{(3)}}~Let $\Vert d\Phi_{0}^{(\varepsilon)} \Vert_{\infty}:=\max_{e\in E_{0}}
\vert d\Phi_{0}^{(\varepsilon)} (\widetilde e) \vert_{g_{0}^{(0)}}$, 
where $\widetilde e$ is a lift of $e\in E_{0}$ to $E$. 
Then there exists a sufficiently small $\varepsilon_{0}>0$ such that
\begin{equation}
\Vert d\Phi_{0}^{(\varepsilon)} \Vert_{\infty}
\leq 2
\big( \min_{e\in E_{0}} {\widetilde m}_{0}(e) \big)^{-1/2}
\Vert d\Phi_{0}^{(0)} \Vert_{\infty}
\nonumber
\end{equation}
for all $0\leq \varepsilon \leq \varepsilon_{0}$.
\end{lm}
{\bf Proof.}~(1)
We first take an orthonormal basis 
$\{ \omega_{1}(\varepsilon), \ldots, \omega_{d}(\varepsilon) \}$ 
of ${\rm Hom}(\Gamma, \mathbb R)
(\subset {\cal H}^{1}_{(\varepsilon)}(X_{0}))$. In particular, we write 
$\omega_{i}=\omega_{i}(0)$ ($ i=1, \ldots, d$). 
Using the basis $\{ \omega_{1}, \ldots, \omega_{d} \}$, we may expand 
$\omega_{i}(\varepsilon)$ as 
$$
\omega_{i}(\varepsilon)=\sum_{j=1}^{d} a(\varepsilon)_{i}^{j} \omega_{j}
\qquad (i=1, \ldots, d).
$$
It follows from $\langle \hspace{-0.8mm} \langle 
\omega_{i}(\varepsilon)^{(\varepsilon)}, 
\omega_{j}(\varepsilon)^{(\varepsilon)}
\rangle  \hspace{-0.8mm} \rangle_{(\varepsilon)}
=\delta_{ij}$
that
\begin{equation}
\sum_{k,l=1}^{d} a(\varepsilon)_{i}^{k} a(\varepsilon)_{j}^{l}
\langle \hspace{-0.8mm} \langle 
\omega_{k}^{(\varepsilon)}, 
\omega_{l}^{(\varepsilon)}
\rangle  \hspace{-0.8mm} \rangle_{(\varepsilon)}
=\delta_{ij}
\qquad (i,j=1,\ldots, d).
\nonumber
\end{equation}
Taking $\varepsilon \searrow 0$ and using Lemma \ref{Harmonic-pdelta},
we obtain
\begin{equation}
\lim_{\varepsilon \searrow 0} 
\sum_{k=1}^{d} a(\varepsilon)_{i}^{k} a(\varepsilon)_{j}^{k}
=\delta_{ij} \qquad (i,j=1,\ldots, d).
\label{aij-conv}
\end{equation}

Now we set
$A(\varepsilon)=(a(\varepsilon)_{i}^{j})_{i,j=1}^{d} \in \mathbb R^{d} \otimes 
\mathbb R^{d}$. Since $\omega_{1}(\varepsilon), \ldots, \omega_{d}(\varepsilon)$ are
linearly independent, $A(\varepsilon)$ is invertible for every $\varepsilon \geq 0$. 
Moreover,
we can rewrite (\ref{aij-conv}) as
$$\lim_{\varepsilon \searrow 0}
\hspace{-1mm}\mbox{ }^{t}
A(\varepsilon)A(\varepsilon)=I_{\mathbb R^{d}}.$$
This convergence yields
$ \frac{1}{2}  \leq \Vert A(\varepsilon) \Vert \leq 2$
for sufficiently small $\varepsilon>0$. Then we have
\begin{eqnarray}
\Vert 
A(\varepsilon) \hspace{-1mm}\mbox{ }^{t}A(\varepsilon)-I_{\mathbb R^{d}} \Vert
&=& \big \Vert A(\varepsilon) \big( \hspace{-1mm}\mbox{ }^{t}A(\varepsilon) A(\varepsilon)
-I_{\mathbb R^{d}} \big) A(\varepsilon)^{-1} \big \Vert
\nonumber \\
&\leq & 4 \Vert   
\hspace{-1mm}\mbox{ }^{t}
A(\varepsilon)A(\varepsilon)-I_{\mathbb R^{d}}
\Vert
\to 0 \qquad \mbox{ as } \varepsilon \searrow 0.
\nonumber
\end{eqnarray}
It means
\begin{equation}
\lim_{\varepsilon \searrow 0} 
\sum_{i=1}^{d} a(\varepsilon)_{i}^{k} a(\varepsilon)_{i}^{l}
=\delta^{kl} \qquad (k,l=1,\ldots, d).
\label{akl-conv}
\end{equation}
From (\ref{akl-conv}), we have
\begin{align}
\lim_{\varepsilon \searrow 0}
\langle {\bf x}, {\bf y} \rangle_{g_{0}^{(\varepsilon)}}
&= 
\lim_{\varepsilon \searrow 0}
\sum_{i=1}^{d} 
\omega_{i}(\varepsilon)[ {\bf{x}} ]_{\Gamma \otimes \mathbb R}
\hspace{1mm} 
\omega_{j}(\varepsilon)[ {\bf{y}} ] _{\Gamma \otimes \mathbb R}
 \nonumber \\
 &=
 \sum_{k,l=1}^{d} \Big(
 \lim_{\varepsilon \searrow 0} 
 \sum_{i=1}^{d} a(\varepsilon)_{i}^{k} a(\varepsilon)_{i}^{l} \Big)
\omega_{k}[ {\bf{x}} ]_{\Gamma \otimes \mathbb R}
\hspace{1mm} 
\omega_{l}[ {\bf{y}} ]_{\Gamma \otimes \mathbb R}
\nonumber
\\
 &=
 \sum_{k=1}^{d} 
\omega_{k}[ {\bf{x}} ]_{\Gamma \otimes \mathbb R} 
\hspace{1mm}
\omega_{k}[ {\bf{y}} ]_{\Gamma \otimes \mathbb R}  = 
\langle {\bf x}, {\bf y} \rangle_{g_{0}^{(0)}}.
\label{0212-2015}
\end{align}
Thus we have shown (1). 
We also mention that 
(\ref{akl-conv}) and (\ref{0212-2015}) imply
\begin{equation} 
\frac{1}{2} \vert {\bf x} \vert_{g_{0}^{(0)}}
\leq
\vert {\bf x} \vert_{g_{0}^{(\varepsilon)}} \leq 2
\vert {\bf x} \vert_{g_{0}^{(0)}} \qquad ({\bf{x}} \in \Gamma \otimes \mathbb R)
\label{212-2015}
\nonumber
\end{equation}
for sufficiently small $\varepsilon \geq 0$.
\\
(2)~Due to the identity (\ref{def-sr}), we can expand
$\vert
d\Phi_{0}^{(\varepsilon)}({\widetilde{e}})
\vert_{g_{0}^{(\varepsilon)}}^{2}$ as
\begin{align}
\big \vert
d\Phi_{0}^{(\varepsilon)}({\widetilde{e}})
\big \vert_{g_{0}^{(\varepsilon)}}^{2}
&=
\sum_{i=1}^{d} 
\omega_{i}(\varepsilon) \big [ \Phi_{0}^{(\varepsilon)}(t({\widetilde e}))
-\Phi_{0}^{(\varepsilon)}(o({\widetilde e})) \big ]_{\Gamma \otimes \mathbb R}^{2}
\nonumber \\
&= \sum_{i=1}^{d} \omega_{i}(\varepsilon)^{(\varepsilon)}(e)^{2}
= \sum_{j,k=1}^{d} \Big ( \sum_{i=1}^{d}  a_{i}^{j}(\varepsilon) a_{i}^{k}
(\varepsilon)  \Big) \omega_{j}^{(\varepsilon)}(e) \omega_{k}^{(\varepsilon)}(e).
\label{Feb9-2015}
\end{align}
Then by combining Lemma \ref{Harmonic-pdelta} and (\ref{akl-conv}) with (\ref{Feb9-2015}), 
we obtain
\begin{eqnarray}
\lim_{\varepsilon \searrow 0}
\sum_{e\in E_{0}} 
p_{\varepsilon}(e)
\big \vert d\Phi_{0}^{(\varepsilon)}({\widetilde e})
\big \vert_{g_{0}^{(\varepsilon)}}^{2} 
m(e)
&=& \lim_{\varepsilon \searrow 0} \sum_{j,k=1}^{d}
\Big ( \sum_{i=1}^{d}  a_{i}^{j}(\varepsilon) a_{i}^{k}
(\varepsilon)  \Big) 
\langle \hspace{-0.8mm} \langle 
\omega_{j}^{(\varepsilon)}, 
\omega_{k} ^{(\varepsilon)} 
\rangle \hspace{-0.8mm} \rangle_{(\varepsilon)}
\nonumber \\
&=&
\sum_{j=1}^{d} \Vert
\omega_{j}^{(0)}
\Vert_{(0)}^{2}
=
\sum_{e\in E_{0}} 
p_{0}(e)
\big \vert d\Phi_{0}^{(0)}({\widetilde e})
\big \vert_{g_{0}^{(0)}}^{2} 
m(e).
\nonumber
\end{eqnarray}
(3)~Recall that $\widetilde{m}_{\varepsilon}(e)$ is 
continuous with respect to $\varepsilon$ and $\widetilde{m}_{0}(e)>0$
for all $e\in E_{0}$. Then by (\ref{0212-2015}) and (\ref{212-2015}),
we can choose a sufficiently small $\varepsilon_{0}>0$ such that
\begin{eqnarray}
\vert d\Phi_{0}^{(\varepsilon)} ({{e'}}) \vert_{g_{0}^{(0)}}^{2}
&\leq & 2 \vert d\Phi_{0}^{(\varepsilon)} ({{e'}}) \vert_{g_{0}^{(\varepsilon)}}^{2}
\nonumber \\
&\leq & 
2
\big( \min_{e\in E_{0}} {\widetilde m}_{\varepsilon}(e) \big)^{-1/2} 
\sum_{e\in E_{0}} 
\big \vert d\Phi_{0}^{(\varepsilon)}({\widetilde e})
\big \vert_{g_{0}^{(\varepsilon)}}^{2} 
{\widetilde m}_{\varepsilon} (e)
\nonumber \\
&\leq & 4 \big( \min_{e\in E_{0}} {\widetilde m}_{0}(e) \big)^{-1/2} 
\sum_{e\in E_{0}} 
\big \vert d\Phi_{0}^{(0)}({\widetilde e})
\big \vert_{g_{0}^{(0)}}^{2} 
{\widetilde m}_{0} (e)
\nonumber \\
&\leq &
4 \big( \min_{e\in E_{0}} {\widetilde m}_{0}(e) \big)^{-1/2} 
\Vert d\Phi_{0}^{(0)} \Vert_{\infty} 
\nonumber
\end{eqnarray}
for all $e' \in E_{0}$ and $0\leq \varepsilon \leq \varepsilon_{0}$.
This completes the proof.
\qed
\subsection{Proof of Theorem \ref{CLT-2}}
To prove Theorem \ref{CLT-2}, the following lemma is essential.
\begin{lm} \label{Conv-CLT2}
For any $f \in C_0^\infty ( (\Gamma \otimes \mathbb{R})_{(0)})$,
as $N \nearrow \infty$, $\varepsilon \searrow 0$ and $N^2 \varepsilon \searrow 0$, we have 
\begin{equation}
\left \| \frac{1}{N \varepsilon^2} ( I - {L}_{(\varepsilon)}^N ) P_\varepsilon f -
P_\varepsilon  \Big( \frac{1}{2}\Delta_{(0)} -
\big \langle \rho_{\mathbb R} (\gamma_p), \nabla_{(0)} 
\big \rangle_{g_{0}^{(0)}}
 \Big) f 
\right\|_\infty \to  0,
\label{conv-CLT02}
\nonumber
\end{equation}
where $\Delta_{(0)}$ and $\nabla_{(0)}$ stand for 
the {\rm{(}}positive{\rm{)}} Laplacian and the gradient 
on 
$(\Gamma \otimes \mathbb{R})_{(0)}$, respectively. 
\end{lm}
{\bf Proof.}~We first take an orthonomal basis $\{ \omega_{1}, \ldots, \omega_{d} \}$ of 
${\rm{Hom}}(\Gamma , \mathbb R)(\subset H^{1}(X_{0},\mathbb R) \cong 
{\cal H}^{1}_{(0)}(X_{0}))$, 
and let $\{ {\bf{v}}_{1}, \ldots, {\bf{v}}_{d} \}$ 
denote its dual basis in $\Gamma \otimes \mathbb R$.
Note that $\{ {\bf{v}}_{1}, \ldots, {\bf{v}}_{d} \}$ is an orthonomal basis of
$(\Gamma \otimes \mathbb R)_{(0)}$. As in the previous section, we denote
$\omega_{i} [{\bf x}]_{\Gamma \otimes \mathbb R}$
by $x_{i}$, and identify ${\bf x}\in \Gamma \otimes \mathbb R$
with $(x_1 , x_2, \ldots, x_d) \in \mathbb R^{d}$. 
Applying Taylor's expansion formula, we have
\begin{eqnarray}
\lefteqn{
f \big (\varepsilon {\Phi}^{(\varepsilon)}_{0} (t(c)) \big)
-f \big (\varepsilon {\Phi}^{(\varepsilon)}_{0} (x) \big)
}
\nonumber \\
&=&
\varepsilon 
\sum_{i=1}^{d} \frac{\partial f}{\partial x_{i}} \big( \varepsilon {\Phi}^{(\varepsilon)}_{0} (x) \big)
\big ( {\Phi}^{(\varepsilon)}_{0} (t(c)) -  {\Phi}^{(\varepsilon)}_{0} (x) 
\big)_{i}
\nonumber 
\\
&\mbox{ }&
+\frac{\varepsilon^2}{2} 
\sum_{i,j=1}^{d}
\frac{\partial^2 f}{\partial x_i \partial x_j} 
\big( \varepsilon {\Phi}^{(\varepsilon)}_{0} (x) \big)
\big ( {\Phi}^{(\varepsilon)}_{0} (t(c)) -  {\Phi}^{(\varepsilon)}_{0} (x) 
\big)_{i}
\big ( {\Phi}^{(\varepsilon)}_{0} (t(c)) -  {\Phi}^{(\varepsilon)}_{0} (x) 
\big)_{j}
+O(N^3\varepsilon^3).
\nonumber
\end{eqnarray}
Recalling that the modified harmonicity of ${\Phi}^{(\varepsilon)}_{0}$ 
and $\rho_{\mathbb R}(\gamma_{p_{\varepsilon}})=\varepsilon 
\rho_{\mathbb R}(\gamma_{p})$,
we see
\begin{equation}
L^{N}_{(\varepsilon)}{\Phi}^{(\varepsilon)}_{0}=
{\Phi}^{(\varepsilon)}_{0}
+N \varepsilon \rho_{\mathbb R}(\gamma_{p}).
\label{N-mh}
\end{equation}
Then we obtain
\begin{align}
( I - L_{(\varepsilon)}^N ) P_{\varepsilon} f (x)
&= -N\varepsilon^{2} \sum_{i=1}^{d} 
\rho_{\mathbb R}(\gamma_{p})_{i}
\frac{\partial f}{\partial x_{i}} \big( \varepsilon { \Phi}^{(\varepsilon)}_{0} (x) \big)
\nonumber \\
& \quad 
-\frac{\varepsilon^2}{2} 
\sum_{i,j=1}^{d}
\frac{\partial^{2} f}{\partial x_{i}\partial x_{j}} \big( \varepsilon {\Phi}^{(\varepsilon)}_{0} (x) \big)
B^{N}_{(\varepsilon)}({{\Phi}}_{0}^{(\varepsilon)})_{ij}(x)
+O(N^3\varepsilon^3),
\nonumber
\end{align}
where
$B^{N}_{(\varepsilon)}({{\Phi}}_{0}^{(\varepsilon)})_{ij}:V \to \mathbb R$
($i,j=1, \ldots, d$, $N\in \mathbb N$) is defined
by
$$ B^{N}_{(\varepsilon)}({{\Phi}}_{0}^{(\varepsilon)})_{ij}(x)
=\sum_{c\in {\Omega}_{x,N}(X)} p_{\varepsilon}(c) \big(
{{\Phi}}_{0}^{(\varepsilon)}(t(c))-
{{\Phi}}_{0}^{(\varepsilon)}(x) \big)_{i}
\big(
{{\Phi}}_{0}^{(\varepsilon)}(t(c))-
{{\Phi}}_{0}^{(\varepsilon)}(x) \big)_{j} \qquad (x\in V).
$$

Next we define ${\cal B}_{(\varepsilon)}
({\Phi}^{(\varepsilon)}_{0})_{ij}: 
V_{0} \to \mathbb R$ ($i,j=1, \ldots, d$) by
$$
{\cal B}_{\varepsilon}
({\Phi}^{(\varepsilon)}_{0})_{ij}(x)
=
\sum_{e\in (E_{0})_{x}}p_{\varepsilon} (e)
\big({\Phi}^{(\varepsilon)}_{0} (t({\widetilde e})) -
{\Phi}^{(\varepsilon)}_{0} (o({\widetilde e})) \big)_{i}
\big( {\Phi}^{(\varepsilon)}_{0} (t({\widetilde e})) -
{\Phi}^{(\varepsilon)}_{0} (o({\widetilde e})) \big)_{j}
\qquad (x\in V_{0}),
$$
where ${\widetilde e}$ is a lift of $e\in E_{0}$ to $E$.
Because ${B}^{N}_{(\varepsilon)}({{\Phi}}_{0}^{(\varepsilon)})_{ij}:V \to \mathbb R$ 
is $\Gamma$-invariant, it holds 
$$
{\cal B}_{(\varepsilon)}({{\Phi}}_{0}^{(\varepsilon)})_{ij}(\pi(x))
={B}^{1}_{(\varepsilon)}({{\Phi}}_{0}^{(\varepsilon)})_{ij}(x)
\qquad (x\in V, ~ i,j=1,\ldots, d).
$$
Using (\ref{N-mh}), we obtain
\begin{eqnarray}
\lefteqn{
{B}^{N}_{(\varepsilon)}({{\Phi}}_{0}^{(\varepsilon)})_{ij}(x)}
\nonumber \\
&=&
\sum_{ c^\prime \in {\Omega}_{x,N-1}(X)} 
p_\varepsilon (c') \sum_{ e \in E_{t(c^\prime )}} p_\varepsilon (e)
\nonumber \\
&\mbox{ }&
\times
 \Big \{
\big( {{\Phi}}_{0}^{(\varepsilon)} (t(e)) -{{\Phi}}_{0}^{(\varepsilon)}(o(e))   \big)_i
+
\big( { {\Phi}}_{0}^{(\varepsilon)} (t(c')) -{{\Phi}}_{0}^{(\varepsilon)}(x)   \big)_i \Big \}
\nonumber \\
&\mbox{ }&
\times
 \Big \{
\big( {{\Phi}}_{0}^{(\varepsilon)} (t(e)) -{{\Phi}}_{0}^{(\varepsilon)}(o(e))   \big)_j
+
\big( { {\Phi}}_{0}^{(\varepsilon)} (t(c')) -{ {\Phi}}_{0}^{(\varepsilon)}(x)   \big)_j \Big \}
\nonumber  \\
&=&
\sum_{c^\prime \in {\Omega}_{x, N-1}(X)} p_\varepsilon (c^\prime )
\sum_{ e \in E_{t(c^\prime )}} p_\varepsilon (e)
\big( {{\Phi}}_{0}^{(\varepsilon)} (t(e)) -{{\Phi}}_{0}^{(\varepsilon)}(o(e))   \big)_i
  \big(  {{\Phi}}_{0}^{(\varepsilon)} (t(e)) -{{\Phi}}_{0}^{(\varepsilon)}(o(e))   \big)_j
 \nonumber \\ 
&\mbox{  }& +
\sum_{c^\prime \in {\Omega}_{x,N-1}(X)}
p_\varepsilon (c^\prime )\sum_{ e \in E_{t(c^\prime )}} p_\varepsilon (e)
\big(   {{\Phi}}_{0}^{(\varepsilon)} (t(c')) -{{\Phi}}_{0}^{(\varepsilon)}(x)  \big)_i
  \big(   {{\Phi}}_{0}^{(\varepsilon)} (t(c')) -{{\Phi}}_{0}^{(\varepsilon)}(x)   \big)_j
  \nonumber \\
  &\mbox{  }&
  +
\sum_{c^\prime \in {\Omega}_{x,N-1}(X)} p_\varepsilon (c^\prime )
\big(  {{\Phi}}_{0}^{(\varepsilon)} (t(c')) -{{\Phi}}_{0}^{(\varepsilon)}(x)  \big)_j
\sum_{ e \in E_{t(c^\prime )}} p_\varepsilon (e)
   \big(   {{\Phi}}_{0}^{(\varepsilon)} (t(e)) -{ {\Phi}}_{0}^{(\varepsilon)}(o(e))    \big)_i
   \nonumber \\
   &\mbox{ }&
 +\sum_{c^\prime \in {\Omega}_{x,N-1}(X)} p_\varepsilon (c^\prime )
  \big(  {{\Phi}}_{0}^{(\varepsilon)} (t(c')) -{{\Phi}}_{0}^{(\varepsilon)}(x)  \big)_i
\sum_{ e \in E_{t(c^\prime )}} p_\varepsilon (e)
   \big(   {{\Phi}}_{0}^{(\varepsilon)} (t(e)) -{{\Phi}}_{0}^{(\varepsilon)}(o(e))    \big)_j
\nonumber \\
%
%
%
&=&L_{(\varepsilon)}^{N-1} \big( {\cal B}_{\varepsilon}({ {\Phi}}_{0}^{(\varepsilon)})_{ij}
\big)(\pi(x))+{B}^{N-1}_{(\varepsilon)}({{\Phi}}_{0}^{(\varepsilon)})_{ij}(x)
\nonumber \\
&\mbox{  }&+ \big( L^{N-1}_{(\varepsilon)} {{\Phi}}_{0}^{(\varepsilon)}(x) -
  {{\Phi}}_{0}^{(\varepsilon)}(x)  \big)_{j} \rho_{\mathbb R} ( \gamma_{p_{\varepsilon}} )_i 
  +\big( L^{N-1}_{(\varepsilon)} { {\Phi}}_{0}^{(\varepsilon)}(x) -
  {{\Phi}}_{0}^{(\varepsilon)}(x)  \big)_{i} \rho_{\mathbb R} ( \gamma_{p_{\varepsilon}} )_j
  \nonumber 
\end{eqnarray}
\begin{eqnarray}
&=&   L_{(\varepsilon)}^{N-1} \big( {\cal B}_{\varepsilon}({ {\Phi}}_{0}^{(\varepsilon)})_{ij}
\big)(\pi(x))+{B}^{N-1}_{(\varepsilon)}({{\Phi}}_{0}^{(\varepsilon)})_{ij}(x)
+ 2(N-1) \rho_{\mathbb R} ( \gamma_{p_{\varepsilon}} )_i 
\rho_{\mathbb R} (\gamma_{p_{\varepsilon}} )_j 
\nonumber \\
&=& \sum_{ k=0}^{N-1} L_{(\varepsilon)}^{k} 
\big( {\cal B}_{(\varepsilon)}({ {\Phi}}_{0}^{(\varepsilon)})_{ij}
\big)(\pi(x)) + 
N(N-1) \varepsilon^{2}  \rho_{\mathbb R} ( \gamma_{p} )_i
 \rho_{\mathbb R} (\gamma_{p} )_j.
\nonumber
\end{eqnarray}
By Lemma \ref{pdelta}, we find that the invariant measure of $L_{(\varepsilon)}$ is $m$.
Then applying Proposition \ref{Ergord-tm2},
we can choose a sufficiently
small $\varepsilon_{0}>0$ such that
\begin{equation}
\frac{1}{N}
\sum_{k=0}^{N-1}
L_{(\varepsilon)}^{k} 
\big( {\cal B}_{(\varepsilon)}({ {\Phi}}_{0}^{(\varepsilon)})_{ij}
\big)(\pi(x))
=\sum_{ x \in X_0} 
{\cal B}_{(\varepsilon)}({ {\Phi}}_{0}^{(\varepsilon)})_{ij}(x)
m(x)
+O_{\varepsilon_{0}} \Big( \frac{1}{N}  \Big)
\label{2-Ergord}
\nonumber
\end{equation}
for all $0\leq \varepsilon \leq \varepsilon_0$. 
Furthermore, as (\ref{B-cal}), it follows from (\ref{def-sr}) that
\begin{eqnarray}
\lefteqn{
\sum_{x \in X_0} 
{\cal B}_{\varepsilon}({ {\Phi}}_{0}^{(\varepsilon)})_{ij}(x) m(x)
}
\nonumber \\
&=& 
 \sum_{e \in E_0} p_{\varepsilon}(e)
\big(  {{\Phi}}_{0}^{(\varepsilon)}(t(e))  - {{\Phi}}_{0}^{(\varepsilon)}(o(e)) \big)_{i}
\big(  {{\Phi}}_{0}^{(\varepsilon)}(t(e))  - {{\Phi}}_{0}^{(\varepsilon)}(o(e)) \big)_{j}  m(o(e))
\nonumber \\
&=&
\Big( \sum_{ e \in E_0} 
p_{\varepsilon}(e)
\omega_i ^{(\varepsilon)} (e) \omega_j ^{(\varepsilon)} (e) m(o(e))
-\langle \gamma_{p_{\varepsilon}}, \omega_{i}  \rangle 
\langle  \gamma_{p_{\varepsilon}}, \omega_{j} \rangle 
\Big)
+
\varepsilon^{2}
\langle \gamma_{p}, \omega_{i}  \rangle 
\langle  \gamma_{p}, \omega_{j} \rangle 
\nonumber \\
&=& 
\langle \hspace{-0.8mm} \langle \omega^{(\varepsilon)}_i , 
\omega^{(\varepsilon)}_j \rangle \hspace{-0.8mm} \rangle_{(\varepsilon)} 
+\varepsilon^{2} \rho_{\mathbb R}(\gamma_{p})_{i}\rho_{\mathbb R}(\gamma_{p})_{j}.
\label{B-cal-2}
\nonumber
\end{eqnarray}
Putting it all together, we obtain
\begin{eqnarray}
\varepsilon^{2} {B}^{N}_{(\varepsilon)}({ {\Phi}}_{0}^{(\varepsilon)})_{ij}(x)
&=&N\varepsilon^{2}  \Big( \langle \hspace{-0.8mm} \langle \omega^{(\varepsilon)}_i , 
\omega^{(\varepsilon)}_j \rangle \hspace{-0.8mm} \rangle_{(\varepsilon)} 
+\varepsilon^{2} \rho_{\mathbb R}(\gamma_{p})_{i}\rho_{\mathbb R}(\gamma_{p})_{j}
\Big)
\nonumber \\
&\mbox{  }&
+N(N-1)\varepsilon^{4}
\rho_{\mathbb R} ( \gamma_{p} )_i
 \rho_{\mathbb R} (\gamma_{p} )_j +\varepsilon^{2} O_{\varepsilon_{0}}(1),
 \nonumber
\end{eqnarray}
and it implies
\begin{eqnarray}
\lefteqn{
\frac{1}{N \varepsilon^2} (I- L_{(\varepsilon)}^N) P_\varepsilon  f(x)}
\nonumber \\
&=&
- \big \langle  \rho_{\mathbb R} (\gamma_p ), \nabla_{(0)}
 f \big (\varepsilon {\Phi}^{(\varepsilon)}_{0} (x) \big) \big \rangle_{g_{0}^{(0)}}
\nonumber \\
&\mbox{  }&
-\frac{1}{2}
\sum_{i,j=1}^{d}
\Big(
\langle \hspace{-0.8mm} \langle \omega^{(\varepsilon)}_i , 
\omega^{(\varepsilon)}_j \rangle \hspace{-0.8mm} \rangle_{(\varepsilon)} 
+\varepsilon^{2} \rho_{\mathbb R}(\gamma_{p})_{i}\rho_{\mathbb R}(\gamma_{p})_{j}
\nonumber \\
&\mbox{ }&
\qquad \qquad
+(N-1)\varepsilon^{2}
\rho_{\mathbb R} ( \gamma_{p} )_i
 \rho_{\mathbb R} (\gamma_{p} )_j +O_{\varepsilon_{0}}\big( \frac{1}{N} \big)
\Big) 
\frac{ \partial^2 f }{\partial x_i \partial x_j } 
( \varepsilon  {\Phi}^{(\varepsilon)}_{0} (x)) + O(N^2 \varepsilon)
\nonumber \\
&=& - 
 \big \langle  \rho_{\mathbb R} (\gamma_p ), (\nabla_{(0)}
 f )\big (\varepsilon { \Phi}^{(\varepsilon)}_{0} (x) \big) \big \rangle_{g_{0}^{(0)}}
\nonumber \\
&\mbox{ }&
-\frac{1}{2} \sum_{i,j=1}^{d} \langle \hspace{-0.5mm} \langle \omega^{(\varepsilon)}_i , 
\omega^{(\varepsilon)}_j \rangle \hspace{-0.5mm} \rangle_{(\varepsilon)} 
 \frac{ \partial^2 f }{\partial x_i \partial x_j } ( \varepsilon {\Phi}^{\varepsilon}_{0} (x) ) 
+ O(N^2\varepsilon)+O_{\varepsilon_{0}}\big( \frac{1}{N} \big).
\nonumber
\end{eqnarray} 

Finally, applying Lemma \ref{Harmonic-pdelta},  we obtain
\begin{equation*}
 \frac{1}{N \varepsilon^2} (I- L_\varepsilon ^N)P_\varepsilon f(x)=
P_\varepsilon 
\Big( \frac{1}{2}\Delta_{(0)} -\big \langle \rho_{\mathbb R} (\gamma_p), \nabla_{(0)} 
\big \rangle_{g_{0}^{(0)}}
 \Big) f (x) +O(N^2\varepsilon) 
 +O_{\varepsilon_{0}}(\frac{1}{N})
\end{equation*}
as $N\nearrow \infty$ and $N^{2}\varepsilon \searrow 0$.
Hence we complete the proof.
\qed
\vspace{2mm} \\
{\bf{Proof of Theorem \ref{CLT-2}.}}
Because the proof is almost same as one of Theorem \ref{CLT-1}, we only give a
sketch. Let $N=N(n)$ be the integer with $n^{1/5} \leq N<1+n^{1/5}$.
We put $k_{N}:=[ ([nt]-[ns])/N]$, $\varepsilon_{N}:=n^{-1/2}$ and $\tau_{N}:=N
\varepsilon_{N}^{2}$. Then $k_{N}\tau_{N} \to (t-s)$ as $N \to \infty$.
Then by recalling Lemma \ref{Conv-CLT2} and
applying Theorem \ref{Trotter-tm} to the case where 
\begin{eqnarray}
& &
{\cal V}=C_{\infty}((\Gamma \otimes \mathbb R)_{(0)}),~~
{\cal V}_{N}=C_{\infty}(X),~~ U_{N}={L}_{(\varepsilon_{N})}^{N}, 
\nonumber \\
& &
T=\frac{\Delta_{(0)}}{2}
-\big \langle \rho_{\mathbb R} (\gamma_p), \nabla_{(0)} 
\big \rangle_{g_{0}^{(0)}},~~D=C^{\infty}_{0}(({\Gamma}\otimes \mathbb R)_{(0)}),
\nonumber 
\end{eqnarray}
we complete the proof.
\qed
\subsection{Proof of Theorem \ref{FCLT-2}}
As is in the previous section, we complete the proof of Theorem \ref{FCLT-2} 
by showing the the convergence of the finite dimensional distributions of $\{ {\bf Y}^{(n^{-1/2},n)}
\}_{n=1}^{\infty}$ (Lemma \ref{FDD-2})
and the tightness of $\{ {\bf Q}^{(n^{-1/2}, n)} \}_{n=1}^\infty$ (Lemma \ref{tight2}).

First, by recalling 
Lemma \ref{dual-Harmonic-conv}, we prepare
\begin{eqnarray} 
\lefteqn{
 \sup_{c\in {\Omega}_{x_{*}}(X)}
 \big \vert {\bf Y}_{t}^{(\varepsilon, n)} (c) - 
{\cal Y}_{t}^{(\varepsilon, n)} (c) \big \vert_{g_{0}^{(0)}}
}
\nonumber \\
&=&
\frac{(nt-[nt])}{\sqrt{n}} \sup_{c\in {\Omega}_{x_{0}}(X)}
\big \vert 
\xi_{[nt]+1}^{(\varepsilon)} (c)-\xi_{[nt]}^{(\varepsilon)} (c)  \big \vert_{g_{0}^{(0)}} 
\nonumber \\
&\leq & n^{-1/2}
\Vert d\Phi_{0}^{(\varepsilon)} \Vert_{\infty}
\nonumber \\
&\leq & 
2\big( \min_{e\in E_{0}} {\widetilde m}_{0}(e) \big)^{-1/2} n^{-1/2}
\Vert d\Phi_{0}^{(0)} \Vert_{\infty}
\rightarrow 0 
\label{2Difference}
\end{eqnarray}
as $n \rightarrow \infty$ uniformly for $0\leq \varepsilon \leq \varepsilon_0$.
We fix $0\leq t_1 < t_2 < \cdots < t_r <\infty$ ($r\in \mathbb N$)
and set the random variable 
${\bf Y}_{t_1,t_2, \ldots , t_r}^{(\varepsilon, n)}: \Omega_{x_*} (X) \rightarrow 
(\Gamma \otimes \mathbb{R} )_{(0)}^{r} $ given by
\begin{equation*}
{\bf Y}_{t_1,t_2, \ldots , t_r}^{(\varepsilon, n)}(c):=
\left( {\bf Y}_{t_1}^{(\varepsilon, n)}(c), \ldots, {\bf Y}_{t_r}^{(\varepsilon, n)} (c)
\right).
\end{equation*}
Noting (\ref{2Difference}) and Theorem \ref{CLT-2}, and 
following the proof of Lemma \ref{FDD1},
we easily obtain
\begin{lm} \label{FDD-2}
\begin{equation*}
{\bf Y}^{(n^{-1/2}, n)}_{t_{1}, \ldots, t_{r}}
\stackrel{ \mathcal{D}}{\longrightarrow} \Big ( B_{t_1}+\rho_{\mathbb R}(\gamma_{p})t_{1}, 
\cdots, B_{ t_r}+\rho_{\mathbb R}(\gamma_{p})t_{r} \Big) 
\quad \mbox{ as } n\to \infty,
\end{equation*}
where $( B_{t} )_{t\geq 0}$ is a $(\Gamma \otimes \mathbb R)_{(0)}$-valued 
standard Brownian motion with $B_{0}={\bf 0}$.
\end{lm}

Hence in this subsection, we concentrate on proving the following lemma:
\begin{lm} \label{tight2}
$\{ {\bf Q}^{(n)} \}_{n=1}^{\infty}$ is tight in  $({\bf W}_{(0)}, {\cal B}({\bf W}_{(0)}))$.
\end{lm}
{\bf Proof.}~
It is enough to show that
there exist some $n_{0}\in \mathbb N$ and $C>0$ independent of $n$ such that
\begin{equation}
\mathbb{E}^{\mathbb P_{x_{*}}} \big [ 
\big | {\bf Y}_{t}^{(n^{-1/2},n)}- {\bf Y}_s^{(n^{-1/2},n)} \big |_{g_{0}^{(0)}}^4 \big] 
\leq C (t-s)^2 \qquad (0\leq s \leq t,~ n\geq n_{0}). 
\label{2tight-4moment}
\end{equation}
We set
$C_{0}:=2 (\min_{e\in E_{0}} {\widetilde m}(e))^{-1/2} \Vert d\Phi_{0}^{(0)} \Vert_{\infty}$
and 
$$C_{1}:=C_{0}+\vert \rho_{\mathbb R}(\gamma_{p}) \vert_{g_{0}^{(0)}},~~
C_{2}:=C_{0}^{4}+C_{0}\vert \rho_{\mathbb R}(\gamma_{p}) \vert_{g_{0}^{(0)}}^{3},~~
C_{3}:=C_{0}^{3}\vert \rho_{\mathbb R}(\gamma_{p}) \vert_{g_{0}^{(0)}}
+\vert \rho_{\mathbb R}(\gamma_{p}) \vert_{g_{0}^{(0)}}^{4}.
$$
Recalling Lemma \ref{dual-Harmonic-conv}, we have
\begin{eqnarray}
& &
\hspace{-10mm}
\big \vert \sum_{e\in E_{x}} p_{\varepsilon}(e)
\big( \Phi_{0}^{(\varepsilon)}(t(e))-\Phi_{0}^{(\varepsilon)}(o(e)) +{\bf x} \big)_{i}^{l}
\big \vert
\leq  
(C_{0}+\vert {\bf x} \vert_{g_{0}^{(0)}})^{l}
\label{219-jyunbi}
\\
& & \qquad \qquad 
(x\in V,~{\bf x}\in \Gamma \otimes \mathbb R,~ i=1,\ldots, d,
~ l\in \mathbb N,~ 0\leq \varepsilon \leq \varepsilon_{0}).
\nonumber
\end{eqnarray}
As in the proof of Lemma \ref{tight1},
we distinguish two cases:
$$ {\mbox{ {\bf{(I)}}  }}~ t-s<n^{-1}, \qquad {\mbox{ {\bf{(II)}} }}~t-s \geq n^{-1}. $$

First, we consider case {\bf{(I)}}. By following the proof of Lemma \ref{tight1}, 
we easily have
\begin{align}
  \vert {\bf Y}_t^{(\varepsilon, n)}- {\bf Y}_s^{(\varepsilon, n)} \vert_{g_{0}^{(0)}}
&\leq 
 2n^{1/2} (t-s) \Vert d\Phi_{0}^{(\varepsilon)} \Vert_{\infty}
 \nonumber \\
& \leq 2C_{0} n^{1/2}(t-s) \qquad (n\in \mathbb N, ~0\leq \varepsilon \leq \varepsilon_{0}).
\nonumber
\end{align}
This yields the desired estimate
\begin{equation}
\mathbb{E}^{\mathbb P_{x_{*}}} \big [ 
\big | {\bf  Y}_t^{(n^{-1/2}, n)}- {\bf Y}_s^{(n^{-1/2}, n)} \big |_{g_{0}^{(0)}}^4 \big] 
\leq
16C_{0}^{4} (t-s)^{2} \qquad (0\leq s \leq t)
\nonumber
\end{equation}
for all $n\geq n_{0}:= 
[ \varepsilon_{0}^{-2}]+1$,
where we used $n^{2}(t-s)^{2}<1$. 
\vspace{2mm}

Next, we consider case {\bf{(II)}}. Let $\mathcal{F}$ be the fundamental domain in 
$X$ containing $x_{*} \in V$ and $M>N$ be two positive integers.
As in the proof of Lemma \ref{tight1}, we have
\begin{eqnarray}
\lefteqn{
\mathbb{E}^{\mathbb P_{x_{*}}}  \big [ 
\big | {\cal Y}_{\frac{M}{n}}^{(\varepsilon,n)}- {\cal Y}_{\frac{N}{n}}^{(\varepsilon,n)}
\big |_{g_{0}^{(0)}}^4 \big ]
}
\nonumber \\
&\leq & d^{2}n^{-2} \max_{i=1,\ldots, d} \max_{x\in {\cal F}}
\Big \{ \sum_{c \in {\Omega}_{x, M-N}(X)} p_{\varepsilon}(c)
\big ( \Phi^{(\varepsilon)}_{0} \big (t(c) \big)-\Phi^{(\varepsilon)}_{0} (x) 
 \big )_{i}^{4}
 \Big \}.
\label{Est-Feb219}
\end{eqnarray}
%
We now fix $i=1,\ldots, d, x\in V$ and set
\begin{eqnarray}
{\mathfrak M}_{\varepsilon}^{l}(k;{\bf{x}})
&:=&\sum_{c\in \Omega_{x,k}(X)}
p_{\varepsilon}(c) \big( \Phi_{0}^{(\varepsilon)}(t(c))-\Phi_{0}^{(\varepsilon)}(x) +{\bf{x}} \big)_{i}^{l}
\nonumber \\
&\mbox{ }&
\qquad  \qquad
(k=1,\ldots, M-N,~ l\in \mathbb N,~ {\bf{x}}\in \Gamma \otimes \mathbb R).
\nonumber
\end{eqnarray}
Then we have
\begin{eqnarray}
\lefteqn{
{\mathfrak M}_{\varepsilon}^{4} \big( k; j \varepsilon \rho_{\mathbb R}(\gamma_{p}) \big)
}
\nonumber \\
&=& \sum_{c'\in {\Omega}_{x,k-1}(X)} p_{\varepsilon} (c') 
\sum_{e\in E_{t(c')}} p_{\varepsilon}(e) \Big \{ \big( 
\Phi^{(\varepsilon)}_{0}(t(e))-\Phi_{0}^{(\varepsilon)} (o(e))
-\varepsilon \rho_{\mathbb R}(\gamma_{p}) \big )_{i} 
\nonumber \\
&\mbox{ }& \quad + \big( 
\Phi_{0}^{(\varepsilon)} (t(c'))-
\Phi_{0}^{(\varepsilon)}(x)+(j+1) \varepsilon \rho_{\mathbb R}(\gamma_{p})
\big )_{i}  \Big \}^{4}
\nonumber \\
&=& \sum_{c'\in {\Omega}_{x,k-1}(X)} p_{\varepsilon}(c') 
\sum_{e\in E_{t(c')}} p_{\varepsilon}(e) \big( 
\Phi^{(\varepsilon)}_{0}(t(e))-\Phi_{0}^{(\varepsilon)} (o(e))
-\varepsilon \rho_{\mathbb R}(\gamma_{p}) \big )_{i}^{4} 
\nonumber \\
&\mbox{ }&
+4\sum_{c'\in {\Omega}_{x,k-1}(X)} p_{\varepsilon}(c') 
\big( 
\Phi_{0}^{(\varepsilon)} (t(c'))-
\Phi_{0}^{(\varepsilon)}(x)+
(j+1) \varepsilon \rho_{\mathbb R}(\gamma_{p})
\big )_{i} 
\nonumber \\
&\mbox{ }& \qquad
\times
\sum_{e\in E_{t(c')}} p_{\varepsilon}(e) \big( 
\Phi^{(\varepsilon)}_{0}(t(e))-\Phi_{0}^{(\varepsilon)} (o(e))
-\varepsilon \rho_{\mathbb R}(\gamma_{p}) \big )_{i}^{3} 
\nonumber \\
&\mbox{ }&
+6
\sum_{c'\in {\Omega}_{x,k-1}(X)} p_{\varepsilon}(c') 
\big( 
\Phi_{0}^{(\varepsilon)} (t(c'))-
\Phi_{0}^{(\varepsilon)}(x)+
(j+1) \varepsilon \rho_{\mathbb R}(\gamma_{p})
\big )_{i}^{2} 
\nonumber \\
&\mbox{ }& \qquad
\times
\sum_{e\in E_{t(c')}} p_{\varepsilon}(e) \big( 
\Phi^{(\varepsilon)}_{0}(t(e))-\Phi_{0}^{(\varepsilon)} (o(e))
-\varepsilon \rho_{\mathbb R}(\gamma_{p}) \big )_{i}^{2} 
\nonumber \\
&\mbox{ }&
+4\sum_{c'\in {\Omega}_{x,k-1}(X)} p_{\varepsilon}(c') 
\big( 
\Phi_{0}^{(\varepsilon)} (t(c'))-
\Phi_{0}^{(\varepsilon)}(x)+
(j+1) \varepsilon \rho_{\mathbb R}(\gamma_{p})
\big )_{i}^{3} 
\nonumber \\
&\mbox{ }& \qquad
\times
\sum_{e\in E_{t(c')}} p_{\varepsilon}(e) \big( 
\Phi^{(\varepsilon)}_{0}(t(e))-\Phi_{0}^{(\varepsilon)} (o(e))
-\varepsilon \rho_{\mathbb R}(\gamma_{p}) \big )_{i} 
\nonumber \\
&\mbox{ }&+
{\mathfrak M}_{\varepsilon}^{4} \big( k-1;  (j+1) \varepsilon \rho_{\mathbb R}(\gamma_{p}) \big)
\qquad (k=2,\ldots, M-N,~ j=0,1,\ldots).
%
\label{Est-Feb18-1}
\end{eqnarray}
Since $\Phi_{0}^{(\varepsilon)}$ enjoys the modified
harmonicity
$
L_{(\varepsilon)} \Phi_{0}^{(\varepsilon)} (x)=
\Phi_{0}^{(\varepsilon)} (x)+\varepsilon \rho_{\mathbb R}(\gamma_{p})$
$(x\in V)$,
the fourth term on the right-hand side of
({\ref{Est-Feb18-1}}) is equal to $0$. 
Furthermore it follows from (\ref{219-jyunbi}) that
%
\begin{equation}
\sum_{e\in E_{y}} p_{\varepsilon}(e) \big( 
\Phi^{(\varepsilon)}_{0}(t(e))-\Phi_{0}^{(\varepsilon)} (o(e))
-\varepsilon \rho_{\mathbb R}(\gamma_{p}) \big )_{i}^{l} 
\leq  
C_{1}^{l} \qquad (y\in V,~ l=2,4)
 \label{Feb19-C1}
\end{equation}
and
\begin{eqnarray}
& &
\hspace{-10mm}
\big \vert
\sum_{c'\in {\Omega}_{x,k-1}(X)} p_{\varepsilon}(c') 
\big( 
\Phi_{0}^{(\varepsilon)} (t(c'))-
\Phi_{0}^{(\varepsilon)}(x)+(j+1) \varepsilon \rho_{\mathbb R}(\gamma_{p})
\big )_{i}
\nonumber \\
& & \qquad \qquad \times 
\sum_{e\in E_{t(c')}} p_{\varepsilon}(e) \big( 
\Phi^{(\varepsilon)}_{0}(t(e))-\Phi_{0}^{(\varepsilon)} (o(e))
-\varepsilon \rho_{\mathbb R}(\gamma_{p}) \big )_{i}^{3} 
\big \vert
\nonumber \\
& &\leq  
\big \{ (k-1) \Vert d\Phi_{0}^{(\varepsilon)} \Vert_{\infty} +(j+1)\varepsilon
\vert \rho_{\mathbb R}(\gamma_{p}) \vert_{g_{0}^{(0)}} \big \}
\big( \Vert d\Phi_{0}^{(\varepsilon)} \Vert_{\infty}+
\varepsilon
\vert \rho_{\mathbb R}(\gamma_{p}) \vert_{g_{0}^{(0)}} \big)^{3}
\nonumber \\
& &
\leq 4(C_{0}^{3}+\varepsilon \vert \rho_{\mathbb R}(\gamma_{p}) \vert_{g_{0}^{(0)}}^{3})
\big \{ (k-1) C_{0} +(j+1)\varepsilon
\vert \rho_{\mathbb R}(\gamma_{p}) \vert_{g_{0}^{(0)}} \big \}
\nonumber \\
& &\leq 4C_{2}(k-1)+4C_{3} (j+1).
\label{Feb19-C2C3}
\end{eqnarray}
Then by combining (\ref{Feb19-C1}), (\ref{Feb19-C2C3}) with ({\ref{Est-Feb18-1}}), we have
\begin{align}
{\mathfrak M}_{\varepsilon}^{4} \big( k; j \varepsilon \rho_{\mathbb R}(\gamma_{p}) \big)
&\leq C_{1}^{4}+4C_{2}(k-1)+4C_{3}(j+1)+
6C_{1}^{2}{\mathfrak M}_{\varepsilon}^{2} \big ( 
k-1; (j+1) \varepsilon \rho_{\mathbb R}(\gamma_{p}) \big)
\nonumber \\
&\quad+
{\mathfrak M}_{\varepsilon}^{4} \big( k-1; (j+1) \varepsilon \rho_{\mathbb R}(\gamma_{p}) \big).
\label{Est-Feb220-1}
\end{align}
Moreover we also have
\begin{eqnarray}
\lefteqn{
{\mathfrak M}_{\varepsilon}^{2} \big ( 
k-1; (j+1) \varepsilon \rho_{\mathbb R}(\gamma_{p}) \big)}
\nonumber \\
&=&
\sum_{c'\in {\Omega}_{x,k-2}(X)} p_{\varepsilon} (c') 
\sum_{e\in E_{t(c')}} p_{\varepsilon}(e) \Big \{ \big( 
\Phi^{(\varepsilon)}_{0}(t(e))-\Phi_{0}^{(\varepsilon)} (o(e))
-\varepsilon \rho_{\mathbb R}(\gamma_{p}) \big )_{i} 
\nonumber \\
&\mbox{ }& \quad + \big( 
\Phi_{0}^{(\varepsilon)} (t(c'))-
\Phi_{0}^{(\varepsilon)}(x)+ (j+2)\varepsilon \rho_{\mathbb R}(\gamma_{p})
\big )_{i}  \Big \}^{2}
\nonumber \\
&=&
\sum_{c'\in {\Omega}_{x,k-2}(X)} p_{\varepsilon} (c') 
\sum_{e\in E_{t(c')}} p_{\varepsilon}(e) \big( 
\Phi^{(\varepsilon)}_{0}(t(e))-\Phi_{0}^{(\varepsilon)} (o(e))
-\varepsilon \rho_{\mathbb R}(\gamma_{p}) \big )_{i}^{2}
\nonumber \\
&\mbox{ }&
+{\mathfrak M}_{\varepsilon}^{2} \big ( 
k-2; (j+2) \varepsilon \rho_{\mathbb R}(\gamma_{p}) \big)
\nonumber \\
&\leq & C_{1}^{2}
+{\mathfrak M}_{\varepsilon}^{2} \big ( 
k-2; (j+2) \varepsilon \rho_{\mathbb R}(\gamma_{p}) \big),
\label{Est-Feb220-2}
\end{eqnarray}
where we used the modified harmonicity 
again for the third line.
It follows from (\ref{Est-Feb220-2}) that
\begin{align}
{\mathfrak M}_{\varepsilon}^{2} \big ( 
k-1; (j+1) \varepsilon \rho_{\mathbb R}(\gamma_{p}) \big)
&\leq  C_{1}^{2}(k-2)
+{\mathfrak M}_{\varepsilon}^{2} 
\big ( 1; (j+k-1) \varepsilon \rho_{\mathbb R}(\gamma_{p}) \big)
\nonumber \\
&\leq  C_{1}^{2}(k-2)+\big( C_{0}+(j+k-1)\varepsilon 
\vert \rho_{\mathbb R}(\gamma_{p}) \vert_{g^{(0)}_{0}} \big)^{2}
\nonumber \\
&\leq   C_{1}^{2}(k-1)+2C_{0}^{2}+2(j+k-1)^{2}\varepsilon^{2}
\vert \rho_{\mathbb R}(\gamma_{p}) \vert_{g^{(0)}_{0}}^{2}.
\label{Est-Feb220-3}
\end{align}
Putting (\ref{Est-Feb220-3}) into (\ref{Est-Feb220-1}),
we obtain
\begin{eqnarray}
{\mathfrak M}_{\varepsilon}^{4} \big( k; j \varepsilon \rho_{\mathbb R}(\gamma_{p}) \big)
&\leq &(2C_{0}^{2}+C_{1}^{4})+(6C_{1}^{4}+C_{1}^{2}+4C_{2})(k-1)+4C_{3}(j+1)
\nonumber \\
&\mbox{ }&
+12(j+k-1)^{2}\varepsilon^{2} 
\vert \rho_{\mathbb R}(\gamma_{p}) \vert_{g^{(0)}_{0}}^{2}
+
{\mathfrak M}_{\varepsilon}^{4} \big( k-1; (j+1) \varepsilon \rho_{\mathbb R}(\gamma_{p}) \big)
\nonumber \\
&=:&C_{4}+C_{5}(k-1)+C_{6}(j+1)+C_{7}(j+k-1)^{2}\varepsilon^{2}
\nonumber \\
&\mbox{ }&
+
{\mathfrak M}_{\varepsilon}^{4} \big( k-1; (j+1) \varepsilon \rho_{\mathbb R}(\gamma_{p}) \big).
\label{Est-Feb220-4}
\nonumber
\end{eqnarray}
Then we have
\begin{eqnarray}
\lefteqn{
\sum_{c \in {\Omega}_{x, M-N}(X)} p_{\varepsilon}(c)
\big ( \Phi^{(\varepsilon)}_{0} \big (t(c) \big)-\Phi^{(\varepsilon)}_{0} (x) 
 \big )_{i}^{4}
 }
\nonumber \\
&=&
{\mathfrak M}_{\varepsilon}^{4} \big(M-N; {\bf{0}} \big)
\nonumber \\
&\leq &
C_{4}+C_{5}(M-N-1)+C_{6}+C_{7}(M-N-1)^{2}\varepsilon^{2}
+
{\mathfrak M}_{\varepsilon}^{4} \big( M-N-1; \varepsilon \rho_{\mathbb R}(\gamma_{p}) \big)
\nonumber \\
&\leq &C_{4}(M-N)+(C_{5}+C_{6}) \Big(\sum_{k=1}^{M-N-1}k \Big)+
C_{7}(M-N)^{3}\varepsilon^{2}+
{\mathfrak M}_{\varepsilon}^{4} \big( 1; (M-N-1)\varepsilon \rho_{\mathbb R}(\gamma_{p}) \big)
\nonumber \\
&\leq &(C_{4}+C_{5}+C_{6})(M-N)^{2}+C_{7}(M-N)^{3}\varepsilon^{2}
+(C_{0}+(M-N) \varepsilon \vert \rho_{\mathbb R}(\gamma_{p}) \vert_{g^{(0)}_{0}})^{4}
\nonumber \\
&\leq & 8C_{0}^{4}+
(C_{4}+C_{5}+C_{6})(M-N)^{2}+C_{7}(M-N)^{3}\varepsilon^{2}
+8\vert \rho_{\mathbb R}(\gamma_{p}) \vert_{g^{(0)}_{0}}^{4}
(M-N)^{4} \varepsilon^{4}.
\label{Est-Feb220-5}
\end{eqnarray}

Finally, we put $\varepsilon=n^{-1/2}$ and $C_{8}=8C_{0}^{4}+C_{4}+C_{5}+C_{6}+C_{7}
+8\vert \rho_{\mathbb R}(\gamma_{p}) \vert_{g^{(0)}_{0}}^{4}$.
By (\ref{Est-Feb219}) and (\ref{Est-Feb220-5}), we have
\begin{eqnarray}
\lefteqn{
\mathbb{E}^{\mathbb P_{x_{*}}} \big [ 
\big | {\bf Y}_t^{(n^{-1/2}, n)}- {\bf Y}_s^{(n^{-1/2}, n)} \big |_{g_{0}^{(0)}}^4 \big]
}
\nonumber \\
& \leq &
\mathbb{E}^{\mathbb P_{x_*}} \Big [  \max_{i,j=0,1}
\Big | {\cal  Y}_{\frac{[nt]+i}{n}}^{(n^{-1/2}, n)}- {\cal  Y}_{\frac{[ns]+j}{n}}^{(n^{-1/2},n)}
\Big |_{g_{0}^{(0)}}^4 \Big]
\nonumber \\
& \leq & C_{8} d^{2}n^{-2}\big \{  1+( [nt]-[ns]+1)^{2}+( [nt]-[ns]+1)^{3}n^{-1}
+( [nt]-[ns]+1)^{4}n^{-2} \big \}
\nonumber \\
&\leq & C_{8} d^{2} \big \{ 1
+3(t-s+\frac{2}{n})^{2} \big \}
\nonumber \\
& \leq & C_{8} d^{2} \big \{ (t-s)^{2} +27(t-s)^{2} \big \}=C (t-s)^{2}
\qquad (n\geq n_{0}=[ \varepsilon_{0}^{-2} ] +1),
\nonumber
\end{eqnarray}
where we used $[nt]-[ns] \leq n(t-s)+1$ for the fourth line and the condition 
$n^{-1}\leq (t-s)$ for the final line.
Thus we obtained our desired estimate
(\ref{2tight-4moment}) for case {\bf{(II)}}. This completes the proof.
\qed
\section{Asymptotic expansion of the transition probability}
\label{Asymptotic expansion}
The main purpose of this section is to prove Theorem \ref{AE-LCLT}.
Recall that $X=(V,E)$ is a crystal lattice in which covering transformation group $\Gamma$ is 
a torsion free abelian group of rank $d$ and torsion free.
Throughout this section, we always assume that the random walk $\{ w_{n} \}_{n=1}^{\infty}$
on $X$ is {\it irreducible}, and let $K$ be the period of the random walk on $X$.
As mentioned in Section 2, this assumption implies the irreducibility of the
corresponding random walk $\{ \pi (w_{n}) \}_{n=0}^{\infty}$ on $X_{0}=(V_{0},E_{0})$, and 
let $K_{0}$ be the period of
the random walk on $X_{0}$. However it should be remarked that $K$ does not coincide 
with $K_{0}$ in general. 
To overcome this difficulty, the following lemma plays a key role.
\begin{lm}[cf. \cite{KSS}]
Suppose that the random walk on the crystal lattice $X$ is of period $K$. Then 
there exists a subgroup $\Gamma_{1}$ of $ \Gamma$ with 
$\mathrm{rank} (\Gamma_1)=d$ such that the corresponding random walk on the
quotient graph $X_{1}=\Gamma_{1} \backslash X$ has the same period $K$.
\end{lm}
{\bf{Proof.}} ~Let $V=\coprod_{k=0}^{K-1} A_k$ be the corresponding 
$K$-partition of $V$. For any $x_{0} \in A_{0}$, we set
\begin{equation*}
\Gamma_{1}:=\{ \sigma \in \Gamma \, | \, \sigma x_0 \in A_0 \}.
\end{equation*}
It is easy to see that $\Gamma_{1}$ is a subgroup of $\Gamma$ with rank $d$. 
Thus it suffices to show that $\Gamma_{1}$ preserves the $K$-partition.
For arbitrary $x_{k} \in A_k$ and $\sigma \in \Gamma_{1}$, 
the step of the walk from $\sigma x_0\in A_0$ to $\sigma x_k$ is same as 
one from $x_0$ to $x_k$, which concludes that $\sigma x_k \in A_k$.
\qed
\vspace{2mm}

Henceforth, without loss of generality, 
we may assume that the covering transformation group $\Gamma$ 
preserves the $K$-partition, 
which implies that the corresponding random walk 
on the quotient graph 
$X_0=\Gamma \backslash X$ 
has the same period $K$. 
\begin{ex} Let $X=(V, E)$ be the $2$-dimensional square lattice with the covering transformation
group $\Gamma=\langle \sigma_{1}, \sigma_{2} \rangle \cong \mathbb Z^{2}$.
Its quotient graph is a $2$-bouquet graph $X_{0}$ {\rm{(}}see Figure $1${\rm{)}}. Then
the period of the simple random walk on $X$ and the one of the corresponding random walk on $X_{0}$
are $2$ and $1$, respectively. On the other hand, if we replace $\Gamma$ by 
a subgroup $\Gamma_{1}=\langle \sigma_{1}', \sigma_{2}' \rangle$ defined as in Figure $2$,
the period of the corresponding random walk on $X_{1}$ is equal to 
$2$. 
\begin{figure}[tbph]
\begin{center}
{\unitlength 0.1in%
\begin{picture}( 43.5500, 16.0000)(  5.8500,-26.0000)%
%
\special{pn 8}%
\special{pa 1000 2400}%
\special{pa 3000 2400}%
\special{fp}%
\special{pa 3000 2000}%
\special{pa 1000 2000}%
\special{fp}%
\special{pa 1000 1600}%
\special{pa 3000 1600}%
\special{fp}%
\special{pa 3000 1200}%
\special{pa 1000 1200}%
\special{fp}%
\special{pa 1200 1000}%
\special{pa 1200 2600}%
\special{fp}%
\special{pa 1600 2600}%
\special{pa 1600 1000}%
\special{fp}%
\special{pa 2000 1000}%
\special{pa 2000 2600}%
\special{fp}%
\special{pa 2400 2600}%
\special{pa 2400 1000}%
\special{fp}%
\special{pa 2800 1000}%
\special{pa 2800 2600}%
\special{fp}%
%
\special{pn 4}%
\special{ar 2000 2000 400 100  0.5191461  3.1415927}%
%
\special{pn 4}%
\special{pa 2335 2055}%
\special{pa 2347 2050}%
\special{fp}%
\special{sh 1}%
\special{pa 2347 2050}%
\special{pa 2278 2057}%
\special{pa 2298 2071}%
\special{pa 2293 2094}%
\special{pa 2347 2050}%
\special{fp}%
%
\special{pn 4}%
\special{ar 2000 2000 400 400  3.1415927  4.4674103}%
%
\special{pn 4}%
\special{pa 1889 1616}%
\special{pa 1903 1612}%
\special{fp}%
\special{sh 1}%
\special{pa 1903 1612}%
\special{pa 1833 1611}%
\special{pa 1852 1627}%
\special{pa 1844 1650}%
\special{pa 1903 1612}%
\special{fp}%
\put(20.4000,-21.3500){\makebox(0,0)[lt]{{\tiny $\sigma_1^\prime$}}}%
\put(17.3000,-17.2000){\makebox(0,0)[lt]{{\tiny $\sigma_2^\prime$}}}%
%
\special{pn 8}%
\special{pa 3000 1800}%
\special{pa 3600 1800}%
\special{fp}%
\special{sh 1}%
\special{pa 3600 1800}%
\special{pa 3533 1780}%
\special{pa 3547 1800}%
\special{pa 3533 1820}%
\special{pa 3600 1800}%
\special{fp}%
\put(40.0000,-18.0000){\makebox(0,0){$X_1=$}}%
%
\special{pn 8}%
\special{ar 4600 1800 300 100  0.0000000  6.2831853}%
\put(36.8500,-12.9500){\makebox(0,0){$\Gamma_1=\langle \sigma_1^\prime, \sigma_2^\prime \rangle$}}%
\put(32.9000,-17.2000){\makebox(0,0){$\pi$}}%
\put(5.8500,-18.0000){\makebox(0,0)[lb]{$X=$}}%
%
\special{pn 8}%
\special{ar 4600 1800 300 260  0.0000000  6.2831853}%
%
\special{sh 1.000}%
\special{ia 1600 2000 35 35  0.0000000  6.2831853}%
\special{pn 8}%
\special{ar 1600 2000 35 35  0.0000000  6.2831853}%
%
\special{sh 1.000}%
\special{ia 2400 2000 35 35  0.0000000  6.2831853}%
\special{pn 8}%
\special{ar 2400 2000 35 35  0.0000000  6.2831853}%
%
\special{sh 1.000}%
\special{ia 2000 2400 35 35  0.0000000  6.2831853}%
\special{pn 8}%
\special{ar 2000 2400 35 35  0.0000000  6.2831853}%
%
\special{sh 1.000}%
\special{ia 1200 2400 35 35  0.0000000  6.2831853}%
\special{pn 8}%
\special{ar 1200 2400 35 35  0.0000000  6.2831853}%
%
\special{sh 1.000}%
\special{ia 2800 2400 35 35  0.0000000  6.2831853}%
\special{pn 8}%
\special{ar 2800 2400 35 35  0.0000000  6.2831853}%
%
\special{sh 1.000}%
\special{ia 2800 1600 35 35  0.0000000  6.2831853}%
\special{pn 8}%
\special{ar 2800 1600 35 35  0.0000000  6.2831853}%
%
\special{sh 1.000}%
\special{ia 2000 1600 35 35  0.0000000  6.2831853}%
\special{pn 8}%
\special{ar 2000 1600 35 35  0.0000000  6.2831853}%
%
\special{sh 1.000}%
\special{ia 1200 1600 35 35  0.0000000  6.2831853}%
\special{pn 8}%
\special{ar 1200 1600 35 35  0.0000000  6.2831853}%
%
\special{sh 1.000}%
\special{ia 1600 1200 35 35  0.0000000  6.2831853}%
\special{pn 8}%
\special{ar 1600 1200 35 35  0.0000000  6.2831853}%
%
\special{sh 1.000}%
\special{ia 2400 1200 35 35  0.0000000  6.2831853}%
\special{pn 8}%
\special{ar 2400 1200 35 35  0.0000000  6.2831853}%
%
\special{sh 1.000}%
\special{ia 4300 1800 35 35  0.0000000  6.2831853}%
\special{pn 8}%
\special{ar 4300 1800 35 35  0.0000000  6.2831853}%
%
\special{sh 0}%
\special{ia 2000 2000 40 40  0.0000000  6.2831853}%
\special{pn 8}%
\special{ar 2000 2000 40 40  0.0000000  6.2831853}%
%
\special{sh 0}%
\special{ia 1600 2400 40 40  0.0000000  6.2831853}%
\special{pn 8}%
\special{ar 1600 2400 40 40  0.0000000  6.2831853}%
%
\special{sh 0}%
\special{ia 2400 2400 40 40  0.0000000  6.2831853}%
\special{pn 8}%
\special{ar 2400 2400 40 40  0.0000000  6.2831853}%
%
\special{sh 0}%
\special{ia 2800 2000 40 40  0.0000000  6.2831853}%
\special{pn 8}%
\special{ar 2800 2000 40 40  0.0000000  6.2831853}%
%
\special{sh 0}%
\special{ia 1200 2000 40 40  0.0000000  6.2831853}%
\special{pn 8}%
\special{ar 1200 2000 40 40  0.0000000  6.2831853}%
%
\special{sh 0}%
\special{ia 1600 1600 40 40  0.0000000  6.2831853}%
\special{pn 8}%
\special{ar 1600 1600 40 40  0.0000000  6.2831853}%
%
\special{sh 0}%
\special{ia 2400 1600 40 40  0.0000000  6.2831853}%
\special{pn 8}%
\special{ar 2400 1600 40 40  0.0000000  6.2831853}%
%
\special{sh 0}%
\special{ia 2800 1200 40 40  0.0000000  6.2831853}%
\special{pn 8}%
\special{ar 2800 1200 40 40  0.0000000  6.2831853}%
%
\special{sh 0}%
\special{ia 2000 1200 40 40  0.0000000  6.2831853}%
\special{pn 8}%
\special{ar 2000 1200 40 40  0.0000000  6.2831853}%
%
\special{sh 0}%
\special{ia 1200 1200 40 40  0.0000000  6.2831853}%
\special{pn 8}%
\special{ar 1200 1200 40 40  0.0000000  6.2831853}%
%
\special{sh 0}%
\special{ia 4900 1800 40 40  0.0000000  6.2831853}%
\special{pn 8}%
\special{ar 4900 1800 40 40  0.0000000  6.2831853}%
%
\special{pn 13}%
\special{pa 1600 2000}%
\special{pa 1955 2000}%
\special{fp}%
\special{sh 1}%
\special{pa 1955 2000}%
\special{pa 1888 1980}%
\special{pa 1902 2000}%
\special{pa 1888 2020}%
\special{pa 1955 2000}%
\special{fp}%
%
\special{pn 13}%
\special{pa 1605 2000}%
\special{pa 1605 1640}%
\special{fp}%
\special{sh 1}%
\special{pa 1605 1640}%
\special{pa 1585 1707}%
\special{pa 1605 1693}%
\special{pa 1625 1707}%
\special{pa 1605 1640}%
\special{fp}%
\put(15.7000,-17.6000){\makebox(0,0)[rt]{{\tiny $\sigma_2$}}}%
\put(18.4500,-19.7000){\makebox(0,0)[rb]{{\tiny $\sigma_1$}}}%
\end{picture}}%
\end{center}
\caption{Square lattice graph and the quotient preserving the bipartition}
\label{figure:square2}
\end{figure}
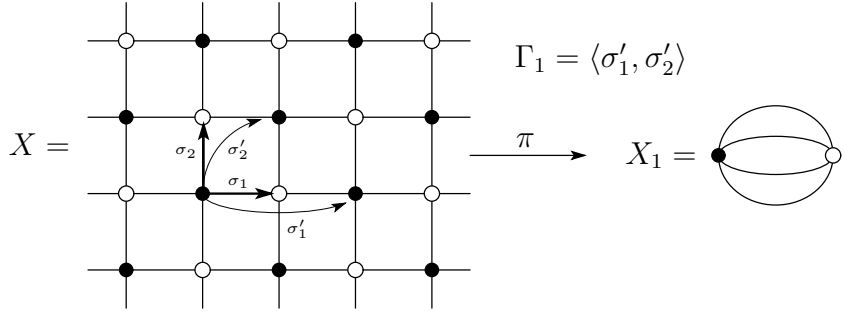
\end{ex}
\subsection{Preliminaries from the twisted transition operators}
For the reader's convenience, 
we first review some
basic results on the twisted transition operators studied in \cite{Kotani contemp, KS00, KSS}.
Let ${\widehat {\rm H}}_{1}(X_{0}, \mathbb Z)$ be the group of unitary characters of 
${\rm H}_{1}(X_{0}, \mathbb Z)$. We identify ${\widehat {\rm H}}_{1}(X_{0}, \mathbb Z)$ with
the Jacobian torus $$J(X_{0}):=\mathrm{H}^{1}(X_{0},\mathbb R)/\mathrm{H}^{1}(X_{0},\mathbb Z)$$ 
by the mapping
$$  \mathrm{H}^{1}(X_{0},\mathbb R) \ni \omega \mapsto \chi_{\omega}\in  
{\widehat {\rm H}}_{1}(X_{0}, \mathbb Z), $$
where
$$
\chi_{\omega}(\sigma):=\exp \big( 2\pi {\sqrt{-1}} \int_{c_{\sigma}} \omega \big)
\qquad (\sigma \in \Gamma)
$$
and $c_{\sigma}$ is a closed path in $X_{0}$ satisfying
$\rho([ c_{\sigma} ])=\sigma$.  We equip a flat metric on $J(X_{0})$ induced by
the metric (\ref{flatmetric-norm}) on ${\rm H}^{1}(X_{0}, \mathbb R)(\cong {\cal H}^{1}(X_{0}))$. 

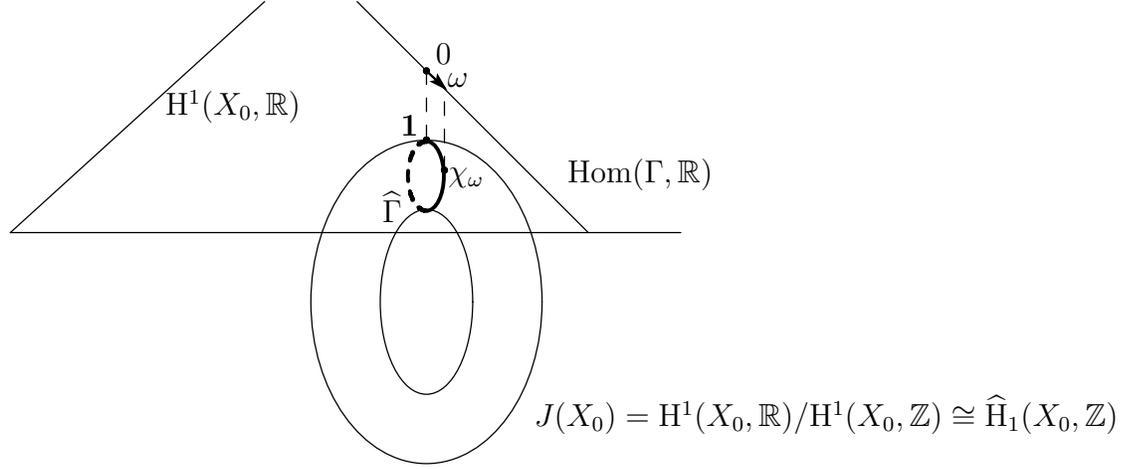
\begin{figure}[tbph]
\begin{center}
\scalebox{1.0}{
{\unitlength 0.1in%
\begin{picture}( 46.4400, 24.2000)( -0.2100,-32.2000)%
%
\special{pn 4}%
\special{sh 1}%
\special{ar 2378 1163 16 16 0  6.28318530717959E+0000}%
\special{sh 1}%
\special{ar 2378 1163 16 16 0  6.28318530717959E+0000}%
%
\special{pn 8}%
\special{ar 2378 2373 605 847  0.0000000  6.2831853}%
%
\special{pn 8}%
\special{ar 2378 2373 242 484  0.0000000  6.2831853}%
%
\special{pn 20}%
\special{ar 2378 1714 91 182  4.7123890  1.5707963}%
%
\special{pn 20}%
\special{pn 20}%
\special{pa 2378 1896}%
\special{pa 2371 1896}%
\special{pa 2371 1895}%
\special{pa 2366 1895}%
\special{pa 2365 1894}%
\special{pa 2362 1894}%
\special{pa 2362 1893}%
\special{pa 2359 1893}%
\special{pa 2359 1892}%
\special{pa 2357 1892}%
\special{pa 2357 1891}%
\special{pa 2355 1891}%
\special{pa 2354 1890}%
\special{pa 2351 1889}%
\special{pa 2350 1888}%
\special{pa 2349 1888}%
\special{pa 2349 1887}%
\special{pa 2347 1887}%
\special{pa 2347 1886}%
\special{pa 2344 1885}%
\special{pa 2344 1884}%
\special{pa 2343 1884}%
\special{pa 2343 1883}%
\special{pa 2341 1883}%
\special{pa 2341 1882}%
\special{pa 2340 1882}%
\special{pa 2340 1881}%
\special{pa 2339 1881}%
\special{pa 2339 1880}%
\special{pa 2338 1880}%
\special{pa 2336 1878}%
\special{pa 2335 1878}%
\special{pa 2335 1878}%
\special{fp}%
\special{pa 2308 1839}%
\special{pa 2308 1839}%
\special{pa 2307 1839}%
\special{pa 2307 1837}%
\special{pa 2306 1836}%
\special{pa 2306 1835}%
\special{pa 2305 1835}%
\special{pa 2305 1833}%
\special{pa 2304 1832}%
\special{pa 2304 1831}%
\special{pa 2303 1830}%
\special{pa 2303 1828}%
\special{pa 2302 1828}%
\special{pa 2302 1826}%
\special{pa 2301 1825}%
\special{pa 2301 1824}%
\special{pa 2300 1823}%
\special{pa 2300 1821}%
\special{pa 2299 1820}%
\special{pa 2299 1819}%
\special{pa 2298 1818}%
\special{pa 2297 1813}%
\special{pa 2296 1812}%
\special{pa 2296 1810}%
\special{pa 2295 1810}%
\special{pa 2295 1807}%
\special{pa 2294 1806}%
\special{pa 2294 1804}%
\special{pa 2293 1803}%
\special{pa 2293 1800}%
\special{pa 2292 1800}%
\special{pa 2291 1793}%
\special{pa 2291 1793}%
\special{fp}%
\special{pa 2282 1741}%
\special{pa 2282 1733}%
\special{pa 2281 1732}%
\special{pa 2281 1696}%
\special{pa 2282 1695}%
\special{pa 2282 1686}%
\special{fp}%
\special{pa 2291 1634}%
\special{pa 2291 1632}%
\special{pa 2292 1631}%
\special{pa 2292 1628}%
\special{pa 2293 1628}%
\special{pa 2293 1625}%
\special{pa 2294 1624}%
\special{pa 2294 1622}%
\special{pa 2295 1621}%
\special{pa 2295 1618}%
\special{pa 2296 1618}%
\special{pa 2297 1613}%
\special{pa 2298 1612}%
\special{pa 2298 1610}%
\special{pa 2299 1609}%
\special{pa 2299 1608}%
\special{pa 2300 1607}%
\special{pa 2300 1605}%
\special{pa 2301 1605}%
\special{pa 2302 1600}%
\special{pa 2303 1600}%
\special{pa 2303 1598}%
\special{pa 2304 1597}%
\special{pa 2304 1596}%
\special{pa 2305 1595}%
\special{pa 2305 1593}%
\special{pa 2306 1593}%
\special{pa 2306 1592}%
\special{pa 2307 1591}%
\special{pa 2307 1589}%
\special{pa 2308 1589}%
\special{pa 2308 1588}%
\special{fp}%
\special{pa 2335 1550}%
\special{pa 2338 1549}%
\special{pa 2338 1548}%
\special{pa 2339 1548}%
\special{pa 2339 1547}%
\special{pa 2340 1547}%
\special{pa 2340 1546}%
\special{pa 2341 1546}%
\special{pa 2341 1545}%
\special{pa 2343 1545}%
\special{pa 2343 1544}%
\special{pa 2344 1544}%
\special{pa 2344 1543}%
\special{pa 2347 1542}%
\special{pa 2347 1541}%
\special{pa 2349 1541}%
\special{pa 2349 1540}%
\special{pa 2352 1539}%
\special{pa 2353 1538}%
\special{pa 2354 1538}%
\special{pa 2355 1537}%
\special{pa 2357 1537}%
\special{pa 2357 1536}%
\special{pa 2359 1536}%
\special{pa 2359 1535}%
\special{pa 2362 1535}%
\special{pa 2362 1534}%
\special{pa 2365 1534}%
\special{pa 2366 1533}%
\special{pa 2371 1533}%
\special{pa 2371 1532}%
\special{pa 2378 1532}%
\special{fp}%
%
\special{pn 8}%
\special{pa 2378 1163}%
\special{pa 2378 1544}%
\special{da 0.070}%
%
\special{pn 4}%
\special{sh 1}%
\special{ar 2378 1526 16 16 0  6.28318530717959E+0000}%
\special{sh 1}%
\special{ar 2378 1526 16 16 0  6.28318530717959E+0000}%
\special{sh 1}%
\special{ar 2378 1526 16 16 0  6.28318530717959E+0000}%
\special{sh 1}%
\special{ar 2378 1526 16 16 0  6.28318530717959E+0000}%
%
\special{pn 8}%
\special{pa 2015 800}%
\special{pa 3225 2010}%
\special{fp}%
%
\special{pn 8}%
\special{pa 3709 2010}%
\special{pa 200 2010}%
\special{fp}%
%
\special{pn 8}%
\special{pa 200 2010}%
\special{pa 1531 800}%
\special{fp}%
\put(24.2600,-11.2700){\makebox(0,0)[lb]{$0$}}%
\put(17.1900,-14.3500){\makebox(0,0)[rb]{$\mathrm{H}^1(X_0,\mathbb{R})$}}%
\put(31.1600,-17.8600){\makebox(0,0)[lb]{$\mathrm{Hom}(\Gamma,\mathbb{R})$}}%
\put(29.4100,-28.6300){\makebox(0,0)[lt]{$J(X_0)= \mathrm{H}^1(X_0,\mathbb{R})/\mathrm{H}^1(X_0,\mathbb{Z})\cong\widehat{\mathrm{H}}_1(X_0,\mathbb{Z})$}}%
\put(22.4800,-17.9200){\makebox(0,0)[rt]{$\widehat{\Gamma}$}}%
\put(23.3600,-15.0200){\makebox(0,0)[rb]{${\bf 1}$}}%
%
\special{pn 8}%
\special{pa 4623 3105}%
\special{pa 3352 3105}%
\special{pa 3352 2567}%
\special{pa 4623 2567}%
\special{pa 4623 3105}%
\special{ip}%
%
\special{pn 8}%
\special{pa 2472 1683}%
\special{pa 2472 1260}%
\special{da 0.070}%
%
\special{pn 13}%
\special{pa 2381 1166}%
\special{pa 2472 1251}%
\special{fp}%
\special{sh 1}%
\special{pa 2472 1251}%
\special{pa 2437 1191}%
\special{pa 2433 1215}%
\special{pa 2410 1220}%
\special{pa 2472 1251}%
\special{fp}%
\put(24.9000,-12.4800){\makebox(0,0)[lb]{$\omega$}}%
%
\special{pn 4}%
\special{sh 1}%
\special{ar 2472 1683 16 16 0  6.28318530717959E+0000}%
\special{sh 1}%
\special{ar 2472 1683 16 16 0  6.28318530717959E+0000}%
\put(24.9300,-16.6800){\makebox(0,0)[lt]{$\chi_\omega$}}%
\end{picture}}%
}
\end{center}

\caption{$\widehat{\Gamma}\subset J(X_{0})$
and $\mathrm{Hom}(\Gamma, \mathbb{R})
\subset \mathrm{H}^1(X_0, \mathbb{R})$.}
\label{figure:torus}
\end{figure}

Let $\widehat \Gamma$ be the group of unitary characters of the covering transformation group
$\Gamma$. By the above mapping, we can also identify $\widehat \Gamma$ with the 
$\Gamma$-Jacobian torus
$$
\mathrm{Jac}^\Gamma := \mathrm{Hom}(\Gamma, \mathbb{R}) 
/ \mathrm{Hom}(\Gamma, \mathbb{Z}).
$$
The canonical surjective homomorphism $\rho: {\rm H}_{1}(X_{0}, \mathbb Z) \to \Gamma$ gives 
rise to an injective homomorphism $\mathrm{Jac}^\Gamma$ into $J(X_{0})$.  We regard 
$\mathrm{Jac}^\Gamma$ as the flat torus with the metric induced by that on $J(X_{0})$.
The tangent space
$T_{\bf 1}{\widehat \Gamma}$ at the trivial character
$\bf 1$ coincides with $\{ \omega \in {\rm H}^{1}(X_{0}, \mathbb R) \vert~\chi_{\omega}\in \widehat \Gamma \}$,
and it is identified with $\mathrm{Hom}(\Gamma, \mathbb{R})$ (see Figure \ref{figure:torus}).
Since the lattice group $\Gamma \otimes \mathbb Z$ in $\Gamma \otimes \mathbb R$
and the lattice group $\mathrm{Hom}(\Gamma, \mathbb{Z})$ in $\mathrm{Hom}(\Gamma, \mathbb{R})$
are dual each other, we observe that the $\Gamma$-Albanese torus 
${\rm Alb}^{\Gamma}=(\Gamma \otimes \mathbb R / \Gamma \otimes \mathbb Z,g_{0})$ is the dual
flat torus of ${\rm Jac}^{\Gamma}$, and hence 
${\rm vol}({\widehat \Gamma})=
{\rm vol}({\rm Jac}^{\Gamma})
={\rm vol}({\rm Alb}^{\Gamma})^{-1}$.

To analyze the $n$-step transition probability $p(n,x,y)$ for the random walk on the crystal lattice
$X=(V,E)$, we need to introduce the {\it{twisted transition operator}} $L_{\chi}$ for
a unitary character $\chi \in {\widehat \Gamma}$.
For each $\chi \in \widehat \Gamma$, we consider the 
$\vert V_{0} \vert$-dimensional inner product space
$$ \ell^{2}_{\chi}= \big \{ f:X \to \mathbb C \big \vert~f(\sigma x)=\chi(\sigma)f(x)~\mbox{for}~\sigma \in \Gamma \}$$
with the inner product
\begin{equation*}
\langle f, g \rangle_\chi= \sum_{x \in \mathcal{F}} f(x) \overline{g(x)} ,
\end{equation*}
where $\mathcal F \subset V$ is a fundamental domain of $X$ for $\Gamma$.
We note that the inner product is independent of the choice of a fundamental domain $\cal F$.

As the transition operator $L$ and its transpose $\,^t L$ preserve $\ell^2_\chi$ 
(see \cite{KSS}), 
we define the {\it{twisted transition operator}} 
$L_{\chi}: \ell^{2}_{\chi} \to \ell^{2}_{\chi}$ and its 
transposed operator $^t L_{\chi}: \ell^{2}_{\chi} \to \ell^{2}_{\chi}$ 
by the restriction of $L$ and $\,^t L$, respectively.
For the trivial character $\chi={\bf 1}$, $(L_{\bf 1},  \ell^{2}_{\bf 1})$
and $(^t L_{\bf 1}, \ell^{2}_{\bf 1})$ are identified with $(L, \ell^{2}(X_{0}))$ and 
$(^t L, \ell^{2}(X_{0}))$, respectively. The family $\{ L_{\chi} \}_{\chi \in {\widehat \Gamma}}$ gives rise
to the direct integral decomposition
$$ (L, \ell^{2}(X) )=\int_{{\widehat \Gamma}}^{\oplus} 
(L_{\chi}, \ell^{2}_{\chi} ) d\chi,
$$
where $d\chi$ denotes the normalized Haar measure on $\widehat \Gamma$. 
As in \cite[Section 7]{KSS},
this decomposition
implies an integral expression of the $n$-step transition probability
\begin{equation}
p(n,x,y) =\int_{\widehat{\Gamma}} \big \langle L^n_\chi f_y, f_x \big \rangle_{\chi} d\chi
\qquad (n \in \mathbb N, ~x,y\in V),
\label{HK-integral}
\end{equation}
where $f_x \in \ell^2_\chi$ is the modified delta function defined by
\begin{equation}
f_{x}(z)
:=
\begin{cases}
\chi (\sigma)
& \text{ if~~  $z=\sigma x$},
\\ 
0
& \text{ otherwise}.
\end{cases}
\nonumber
\end{equation}
It follows from $\langle L^n_\chi f_y, f_x \big \rangle_{\chi}=0$ ($x\in A_i, y \in A_j, n\neq Kl+j-i$)
that $p(n,x,y)=0$.

By virtue of the Perron--Frobenius theorem for the random walk with period $K$, 
the twisted transition operator $L_{\chi}$ has $K$ maximal
eigenvalues $\mu_{0}(\chi), \ldots, \mu_{K-1}(\chi)$. For the eigenvalue $\mu_{k}(\chi)$ $(k=0,\ldots, K-1)$,
we denote the corresponding right eigenfunction 
and left eigenfunction 
by $\phi_{k, \chi}$ and $\psi_{k, \chi}$, respectively.
Specifically,
\begin{equation*}
L_\chi \phi_{k, \chi}= \mu_k (\chi) \phi_{k, \chi}, \quad 
\,^t L_{\chi} \psi_{k, \chi}=\overline{ \mu_k(\chi)} \psi_{k, \chi}.
\end{equation*}
Normalizing $\phi_{k, \chi}$ and $\psi_{k, \chi}$, we may assume 
$$\langle \phi_{k, \chi},\phi_{k, \chi} \rangle_\chi= 
\langle \phi_{k, \chi}, \psi_{k, \chi} \rangle_\chi=1. $$
It should be noted that the function $\mu_{k}=\mu_{k}(\chi)$ is a
continuous function in $\chi \in \widehat \Gamma$ and 
$$ \mu_k({\bf 1})=\exp \big( \frac{2\pi k}{K} \sqrt{-1} \big) \qquad (k=0, \ldots, K-1). $$
Applying Proposition \ref{PP-4.6},
we observe that the eigenvalue $\mu_{k}(\chi)$ is simple and 
$\mu_{k}(\chi)$, $\phi_{k, \chi}$, $\psi_{k, \chi}$
are smooth in a small neighborhood $U(\bf 1)(\subset {\widehat \Gamma})$ of
the trivial character $\bf 1$. Moreover, the operator norm
$\| L_\chi \| $ is equal to $1$ if and only if $\chi={\bf 1}$. (see \cite{Kotani contemp, KSS}
for details.)

We next decompose $\ell ^2_\chi$ as $\ell ^2_\chi= \oplus_{k=0}^{K-1}
\langle \phi_{k, \chi} \rangle \oplus {\cal V}_\chi$,
where $${\cal V}_\chi :=\left\{ f \in \ell ^2_\chi \, \big | \, 
\langle f, \psi_{k, \chi} \rangle_\chi =0~~(0\leq k \leq  K-1) \right\}.$$
More precisely, for $f \in \ell^2_\chi$, 
there exists a unique $f_{\cal V_{\chi}}\in {\cal V}_{\chi}$ such that
\begin{equation}
f= \sum_{k=0}^{K-1} \langle f, \psi_{k, \chi} \rangle_\chi \phi_{k, \chi} 
+ f_{\cal V_{\chi}}.  
\label{spectral decomposition}
\end{equation}
Combining (\ref{HK-integral}) with (\ref{spectral decomposition}), we obtain 
\begin{equation}
p(n,x,y)=
\int_{\widehat{\Gamma}}  \Big ( \sum_{k=0}^{K-1} \mu_k (\chi)^n \phi_{k, \chi} (x) 
\overline{ \psi_{k, \chi} (y)} +
\langle L^n_\chi \{ (f_y)_{\cal V_{\chi}} \},f_x \rangle_\chi   \Big ) d\chi
\label{decomp}
\end{equation}
for $x\in A_i$, $y \in A_j$ and $n=Kl+j-i$.
Since $L_{\chi}$ preserves ${\cal V}_\chi$, and 
$\Vert L_\chi \vert_{\cal V_\chi} \Vert <1-\varepsilon$ for some $\varepsilon >0$ 
uniformly in $\chi \in \widehat{\Gamma}$ (see \cite{Kotani contemp, PP}), we have
\begin{equation}
\Big \vert \int_{\widehat{\Gamma}} 
\langle L^n_\chi \{ (f_y)_{\cal V_{\chi}} \} ,f_x \rangle_\chi  d\chi \Big \vert \leq C (1-\varepsilon)^{n}
\label{exp-falloff1}
\end{equation}
for some positive constant $C$ independent of $x$ and $y$. Therefore 
it suffices to discuss an precise asymptotic behavior of the first term of 
(\ref{decomp}) for proving Theorem \ref{AE-LCLT}.

To look more closely at the first term of the integrand in (\ref{decomp}),
we make use of the correspondence $\omega \mapsto \chi_{\omega}=
\exp( 2\pi {\sqrt{-1}} \langle \omega, \cdot \rangle)$ between small 
$\omega \in {\rm Hom}(\Gamma, \mathbb R)$ and $\chi_{\omega}\in U({\bf 1})$.
For each $\omega \in {\rm Hom}(\Gamma, \mathbb R) (\subset \mathrm{H}^{1}(X_{0}, \mathbb R)$), 
we define a function $s_{\omega}\in {\ell}^{2}_{\chi_{\omega}}$ by
$$ s_{\omega}(x):=\exp \Big( 2\pi {\sqrt{-1}} \int_{x_*}^x \widetilde{ \omega} \Big)=
\exp \Big ( 2\pi {\sqrt{-1}} \langle \omega, \Phi_{0}(x) \rangle \Big )
\qquad (x\in V), $$
where $x_{*}\in V$ is a fixed base point with $\Phi_{0}(x_{*})={\bf 0}$
and $\widetilde \omega$ stands for the lift of $\omega$ to $X$.
%
We define a unitary map 
$S: \ell^{2}(X_{0}) \to \ell_{\chi_{\omega}}^{2}$ by 
$$(Sf)(x):={\widetilde f}(x)s_{\omega}(x) \qquad (x\in V),$$
where $\widetilde f$ is the lift of $f \in \ell^{2}(X_{0})$ to $X$.

We now introduce operators $H_{\omega}:=S^{-1}L_{\chi_{\omega}}S$
and 
$H_{\omega}^*:=S^{-1} (^{t}L_{\chi_{\omega}})S$ acting on $\ell^{2}(X_{0})$.
Then 
we have
\begin{equation}
\begin{split}
H_{\omega}f(x_0) &:=\sum_{e\in (E_{0})_{x_0}} p(e) \exp ({2\pi \sqrt{-1}}\omega(e)) f(t(e))
\qquad (x_0 \in V_0),
\\
H_{\omega}^* f(x_0) &:=\sum_{e\in (E_0)_{x_0}} p({\overline{e}}) 
\exp ({2\pi \sqrt{-1}}\omega (e)) f(t(e)) \qquad (x_0 \in V_0).
\\
\end{split}
\nonumber
\end{equation}
We put $\phi_{k, \omega} :=S^{-1}\phi_{k, \chi_{\omega}}$ and 
$\psi_{k, \omega}:=S^{-1}\psi_{k, \chi_{\omega}}$ for $k=0,\ldots, K-1$. 
Then we easily see that $\phi_{k, \omega}$ and $\psi_{k, \omega}$ 
are the eigenfunctions of the operators $H_\omega$ and $H^*_\omega$
corresponding to the eigenvalues $\mu_k (\chi_{\omega})$ and 
$\overline{ \mu_k (\chi_{\omega})}$, respectively. It is obvious
that $\mu_{k} (\chi_{\omega})$,
$\phi_{k, \omega}$ and  $\psi_{k, \omega}$ depend smoothly on $\omega$ around $0$.
Applying the Perron--Frobenius theorem and Proposition \ref{PP-4.6}, 
we easily see
$$ \mu_{k}(\chi_{\omega})=\exp 
\Big(\frac{2k \pi \sqrt{-1}}{K} \Big) \mu_{0}(\chi_{\omega}), \quad  
\phi_{k, \omega}=T^k \phi_{0, \omega}, \quad \psi_{k, \omega}=T^k \psi_{0, \omega}$$ 
for sufficiently small $\omega\in {\rm Hom}(\Gamma, \mathbb R)$,
where $T:\ell^2(X_0) \rightarrow \ell^2(X_0)$ is defined by
\begin{equation*}
Tf(x_{0})=\exp \Big(\frac{2k\pi \sqrt{-1}}{K} \Big) f(x_{0}) \qquad (x_{0}\in \pi (A_{k}),~  k=0,\ldots, K-1).
\end{equation*} 
Then we obtain
\begin{eqnarray}
\lefteqn{
\mu_{k}(\chi_{\omega})^{n}  \phi_{k, \chi_{\omega}}(x) {\overline{\psi_{k, \chi_{\omega}}(y) }}
}
\nonumber \\
&=&
\exp
\Big(\frac{2k(n+i-j) \pi \sqrt{-1}}{K} \Big)
\mu_{0}(\chi_{\omega})^{n} 
\phi_{0,\omega}(\pi(x)) {\overline{\psi_{0,\omega}(\pi(y))}}
\exp \Big( 2\pi {\sqrt{-1}} \int_{y}^x \tilde{ \omega} \Big)
\nonumber \\
&=&
\mu_{0}(\chi_{\omega})^{n} 
\phi_{0,\omega}(\pi(x)) {\overline{\psi_{0,\omega}(\pi(y))}}
\exp \Big( -2\pi {\sqrt{-1}} \langle \omega, \Phi_{0}(y)-\Phi_{0}(x) \rangle \Big)
\label{chi-omega-henkan}
\end{eqnarray}
for $x\in A_{i}, y\in A_{j}, n=Kl+j-i$, $k=0,\ldots, K-1$
and sufficiently small $\omega\in {\rm Hom}(\Gamma, \mathbb R)$.

Now, we are going to analyze the behavior of $\mu_0( \chi_{\omega})$,
$\phi_{0, \chi_{\omega}}$ and $\psi_{0, \chi_{\omega}}$
around $\omega=0$.
For this sake, we take a path $\chi_{t}=\chi_{t\omega}$ in 
$\widehat \Gamma$ through the trivial character $\bf 1$ at $t=0$.
We put 
$$H_t=H_{\chi_{t\omega}},~H^*_t=H^*_{\chi_{t\omega}},~
\phi_t =\phi_{t\omega}=\phi_{0, \chi_{t\omega}},~\psi_t=\psi_{t\omega}=\psi_{0,\chi_{t\omega}},~ 
\lambda_{\omega} (t)= - \log \mu_{0}(\chi_{t\omega }) $$
near $t=0$. In the sequel, we denote the $k$-th derivative of functions 
$\lambda_{\omega}(t)$, $\phi_t(x_0)$ and $\psi_t(x_0)$ in $t$ by 
$\lambda^{(k)}_{\omega}(t)$, 
$\phi_t^{(k)}(x_0)$ and $\psi_t^{(k)}(x_0)$, respectively. 

To prove Theorem \ref{AE-LCLT}, we need  
$\lambda^{(i)}(0)$~($0\leq i \leq 4$) and $\phi_0^{(j)}(x_0), \psi_0^{(j)}(y_0)$~($0\leq j \leq 2$). 
Differentiating both sides of $H_{t}\phi_{t}=\exp(-\lambda_{\omega}(t)) \phi_{t}$ in $t$ at $t=0$,
and following the argument in \cite[Lemma 3.2]{KS00}, we obtain 
\begin{lm} \label{eigen-bibun-1}
\begin{eqnarray*}
\lambda (0)&=&0, 
\\
\lambda^\prime (0)&=&-2\pi {\sqrt{-1}} \langle \gamma_p, \omega \rangle, \\
\lambda^{\prime \prime}(0)&=&
4\pi^2 \Big( \sum_{e \in E_0} p(e) \omega(e)^2 m(o(e)) - \langle \gamma_p ,\omega \rangle^2 \Big) 
=4\pi^2 \| \omega \|^2, \\
\lambda^{(3)}(0)&=& 
8\pi^3 \sqrt{-1} \sum_{e \in E_0} p(e)\omega (e)^3 m(o(e)) -24\pi^2 \sqrt{-1}
\| \omega \|^2 \langle \gamma_p, \omega \rangle -8\pi^3 \sqrt{-1} 
\langle \gamma_p, \omega \rangle^3 
\nonumber \\
& & \qquad -6\pi \sqrt{-1} |V_0|^{1/2}  \sum_{e \in E_0} p(e)\omega (e) d \phi_0^{\prime \prime} (e) m(o(e)),  
\\
\lambda^{(4)} (0)  &=&
-16\pi^4 \sum_{ e \in E_0} p(e) \omega (e)^4 m(o(e)) +48 \pi^4 \| \omega \|^4 
\nonumber \\
&\mbox{ }&
+64 \pi^4 \langle \gamma_p , \omega \rangle \sum_{e \in E_0} p(e) \omega (e)^3 m(o(e)) 
%
-96\pi^4 \langle \gamma_p, \omega \rangle^2 \| \omega \|^2 -48 \pi^4 
\langle \gamma_p, \omega \rangle^4 
\nonumber \\
&\mbox{ }&
-48\pi^2|V_0|^{1/2} \langle \gamma_p, \omega \rangle \sum_{e \in E_0 } p(e) \omega (e) d\phi_0^{\prime \prime} (e) m(o(e))  \nonumber \\
&\mbox{  }&
+24\pi^2 |V_0|^{1/2} \sum_{e \in E_0} p(e) \omega(e)^2 
\Big( \phi_0^{\prime \prime } (t(e)) -\sum_{z \in V_0} \phi_0^{\prime \prime } (z) m(z) \Big)
m(o(e)) \nonumber \\
&\mbox{ }&
-8\pi \sqrt{-1} |V_0|^{1/2} \sum_{e \in E_0} p(e) \omega(e) d\phi_0^{(3)} (e) m(o(e)).  
\end{eqnarray*}
\end{lm}
\begin{re}
Differentiating both sides of $H^*_t \psi_{t}=\exp(-{\overline{\lambda_{\omega}(t)}}) \psi_{t}$ 
four times in $t$ at $t=0$, we also obtain 
\begin{eqnarray*}
\lambda^{(3)}(0)&=& 8\pi^3 \sqrt{-1} \sum_{e \in E_0} p(e) \omega(e)^3 m(o(e)) - 24\pi^2\sqrt{-1}
\| \omega \|^2 \langle \gamma_p, \omega \rangle -8\pi^3 \sqrt{-1} \langle \gamma_p, \omega \rangle^3
\nonumber \\
& &\qquad +12\pi^2 |V_0|^{-1/2} \sum_{e \in E_0} p(e) \omega(e)^2  \overline{\psi_0^\prime (o(e))},
\end{eqnarray*}
and
\begin{eqnarray*}
\lambda^{(4)}(0)&=&-16\pi^4\sum_{e \in E_0} p(e)\omega(e)^4 m(o(e)) +48\pi^4 \| \omega \|^4 
\nonumber \\
&\mbox{ }&
+64\pi^4 \langle \gamma_p, \omega \rangle \sum_{e \in E_0}p(e) \omega(e)^3 m(o(e)) 
-96\pi^4 \langle \gamma_p , \omega \rangle^2 \| \omega \|^2 
\nonumber \\
&\mbox{ }&
-48\pi^4 \langle \gamma_p, \omega \rangle^4 
-96\pi^3\sqrt{-1} |V_0|^{-1/2} \langle \gamma_p, \omega \rangle 
\sum_{ e \in E_0} p(e) \omega (e)^2 \overline{ \psi_0^\prime (o(e)) } \nonumber \\
&\mbox{  }&+32\pi^3 \sqrt{-1} |V_0|^{-1/2} \sum_{e \in E_0} p(e) \omega(e)^3 \overline{ \psi_0^\prime (o(e)) }   \nonumber \\
&\mbox{  }&+24\pi^2 |V_0|^{-1/2} \sum_{ e \in  E_0} p(e) \omega(e)^2 
\Big( \psi_0^{\prime \prime} (o(e)) -\sum_{z \in V_0} \psi_0^{\prime \prime } (z) \cdot m(o(e)) \Big).
\end{eqnarray*}
\end{re}

Changing eigensections $s_{\omega}$ as in the proof of
Lemma \ref {lemma of expansion} below, if necessary, we obtain the following:
\begin{lm}   \label{eigen-bibun-2}
\begin{align*}
& \phi_0(x_0)= |V_0|^{-1/2}, \quad \psi_0(x_0) =|V_0|^{1/2} m(x_0), \quad
 \phi^{\prime }_0(x_0 )=0 \qquad (x_{0}\in V_{0}).  
 \end{align*}
Furthermore $\psi_0^{ \prime}$ is a purely imaginary-valued first order polynomial of $\omega$
satisfying
\begin{equation*}
\begin{cases} 
(I-\,^t L) \psi_0^{ \prime} (x_0) =2 \pi \sqrt{-1} |V_0|^{1/2}
\Big( \sum_{e \in (E_0)_{x_0}} p(\overline{e}) \omega(e) m(t(e)) 
\\
\qquad \qquad \qquad \qquad \qquad \quad
-m(x_{0}) \sum_{e \in E_0} p(\overline{e}) \omega(e) m(t(e)) \Big)  \qquad \quad (x_{0}\in V_{0})
\\
\sum_{z \in V_0}  \psi_0^\prime (z)=0, 
\end{cases}
\end{equation*}
and
$\phi_0^{\prime \prime}$ and $\psi_0^{\prime \prime}$ are real-valued second 
oder polynomials of $\omega$, satisfying
\begin{equation*}
\begin{cases}
(I-L)\phi_0^{\prime \prime} (x_0)=-4\pi^2 |V_0|^{-1/2}
\Big( \sum_{e \in (E_0)_{x_0}} p(e) \omega(e)^2 
\\
\qquad \qquad \qquad \qquad \qquad \qquad \quad
-\sum_{e \in E_0} p(e) \omega(e)^2 m(o(e)) \Big)  \qquad \quad (x_{0}\in V_{0})
\\ 
\sum_{z \in V_0} \phi_0^{\prime \prime} (z)=0
\end{cases}
\label{phi^2}
\end{equation*}
and
\begin{equation*}
\begin{cases}
(I-\,^t L) \psi_0^{\prime \prime} (x_0)=
4\pi \sqrt{-1} \Big( \sum_{e \in (E_0)_{x_0}} 
p(\overline{e}) \omega(e) \psi_0^\prime (t(e)) 
+\langle \gamma_p , \omega \rangle \psi_0^\prime (x_{0}) \Big)
\\
\qquad \qquad \qquad \qquad
-4\pi^2 \Big( \sum_{e \in (E_0)_{x_0}} p( \overline{e}) \omega(e)^2 \psi_0(t(e))
\\
\qquad \qquad \qquad \qquad \qquad \quad
-m(x_{0}) \sum_{e  \in E_0} p(\overline{e}) \omega(e)^2 \psi_0 (t(e))  \Big) \qquad (x_{0}\in V_{0})
\label{psi^2}
\\ 
\sum_{z \in V_0} \psi_0^{\prime \prime} (z)=-\vert V_{0} \vert \sum_{z\in V_{0}} \phi''_{0}(z)m(z),
\end{cases}
\end{equation*}
respectively.
\end{lm}

The following lemma plays a key role in the proof of Theorem \ref{AE-LCLT}.
\begin{lm}[cf. \mbox{\cite[Lemma 3.2]{KS00}}]\label{lemma of expansion}
For any $k \in \mathbb{N}$, 
by changing eigensections $s_\omega$ if necessary, 
$ \lambda^{(i)}_{\omega}(0)$, $\phi_0^{(i)}(x_0)$, $\psi_0^{(i)}(x_0)$ 
{\rm{(}}$1\leq i \leq k$, $x_0 \in V_0${\rm{)}}
are the $i$-th order real coefficient homogeneous polynomials 
of $\sqrt{-1} \omega$. Here 
the polynomial of $\sqrt{-1}\omega$ 
means that the polynomial of $\sqrt{-1}u_1,\ldots, \sqrt{-1}u_d$ when
$\omega=u_1\omega_1+u_1 \omega_2+ \cdots +u_d \omega_d$, where 
$\{ \omega_{1}, \ldots, \omega_{d} \}$ is an orthonomal basis of 
${\rm Hom}(\Gamma, \mathbb R)$.
\end{lm}
{\bf Proof.} 
We proceed by induction on $k$. 
By Lemmas \ref{eigen-bibun-1} and
\ref{eigen-bibun-2}, we easily see that the desired assertion holds for $k=1$.
 Assume that the assertion is true for all $i \leq k-1$. 
  We recall that the eigenvalue
$\lambda_{\omega}(t)$ and the corresponding eigenfunctions $\phi_t$, $\psi_t$ satisfy
$$
H_t \phi_t (x_0) = e^{-\lambda (t)} \phi_t (x_0), \quad 
H_t^* \psi_t (x_0) = e^{- \overline{ \lambda (t)}} \psi_t (x_0) \qquad (x_0 \in V_0).
$$
Taking the $k$-th derivative of both sides and letting $t \to 0$, we have
\begin{eqnarray}
(I-L) \phi_0^{(k)} (x_0) &=& \sum_{i=0}^{k-1} \binom{k}{i} \Big \{ 
\sum_{e \in (E_0)_{x_0}} p(e) \left( 2\pi \sqrt{-1} \omega(e) \right)^{k-i} \phi_ 0^{(i)} (t(e))
\nonumber \\
&\mbox{  }&
- \big( \frac{d}{dt} \big)^{k-i} \Big \vert_{t=0}
e^{ -\lambda (t) }  \cdot \phi_0^{(i)} (x_0) \Big \}
\qquad (x_{0}\in V_{0}),
 \label{derivative of phi}\\
(I-\,^tL )\psi_0^{(k)} (x_0) &=& \sum_{ i=0} ^{k-1} \binom{k}{i} \Big \{
\sum_{e \in (E_0)_{x_0}} p(\overline{e})  \left( 2\pi \sqrt{-1} \omega(e)  \right)^{k-i}  \psi_ 0^{(i)} (t(e))
\nonumber \\
&\mbox{  }&
-\big( \frac{d}{dt} \big)^{k-i} \Big \vert_{t=0} 
e^{ -\overline{\lambda (t) }}  \cdot \psi_0^{(i)} (x_0)
\Big \}
\qquad (x\in V_{0}). \label{derivative of psi}
\end{eqnarray}
Taking sum of both sides of (\ref{derivative of phi}) for $x_0 \in V_0$ with the invariant measure $m$, the 
left-hand side vanishes. By the assumption, 
we find that $\lambda^{(k)} (0)$ is a desired polynomial.
We can obtain the same result by 
taking sum of the both sides of (\ref{derivative of psi}) for $ x_0 \in V_0$.

In the case $k$ is even, substituting $\lambda^{(k)}(0)$ into  (\ref{derivative of phi}) and 
(\ref{derivative of psi}), we observe that 
\begin{equation*}
(I-L) \Im \phi_0^{(k)} =0, \quad 
(I-\,^t L)\Im \psi_0^{(k)} =0.
 \end{equation*}
These imply that $\Im \phi_0^{(k)}=a_k(\omega)$ and 
$\Im \psi_0^{(k)} =b_k(\omega) m$ for some constant functions 
$a_k(\omega) $ and $b_k(\omega)$ on $V_0$ depending on 
$\omega \in \mathrm{Hom}(\Gamma, \mathbb{R})$. 
Clearly,
$a_k(s\omega)=s^k a_k(\omega)$ and 
$b_k(s \omega)=s^k b_k (\omega)$ hold.
Here we note that 
 \begin{align}
0=& 
\big( \frac{d}{dt} \big)^{k} \Big \vert_{t=0} 
\langle \phi_t, \psi_t \rangle_{\ell^{2}(X_{0})} \nonumber \\
=& |V_0|^{1/2}\sum_{z \in V_0} \phi_0^{(k)} (z) m(z) + 
\sum_{ i=1}^{k-1} \binom{k}{i} \sum_{z \in V_0}
\phi_0^{(k-i)} (z) \overline{ \psi_0^{(i)} (z) }
+ |V_0|^{-1/2} \sum_{ z \in V_0} \overline{ \psi_0^{(k)} (z) }. \label{<phi,psi>}
 \end{align}
 Taking the imaginary part, by the assumption of $\phi_0^{(i)}$ and $\psi_0^{(i)}$ 
for $i \leq k$, we obtain $b_k(\omega) = |V_0| a_k(\omega)$.
By replacing eigensection $s_\omega$ with
$
s_\omega \exp \left( \sqrt{-1} \frac{a_k(\omega) }{ k!} |V_0| ^{1/2} \right)
$, and eigenfunctions $\phi_\omega$, $\psi_\omega $ with
$\phi_\omega \exp \left(  -\sqrt{-1} \frac{a_k(\omega) }{ k!} |V_0| ^{1/2} \right)$, 
$\psi_\omega \exp \left(  -\sqrt{-1} \frac{a_k(\omega) }{ k!} |V_0| ^{1/2} \right)$, 
respectively, we can assume that $a_k=b_k=0$. 
Then $\phi_0^{(k)}$ and $\psi_0^{(k)}$ are real functions satisfying 
(\ref{derivative of phi}), (\ref{derivative of psi}), respectively. 
Since the right-hand side of (\ref{derivative of phi}) and (\ref{derivative of psi}) are 
$k$-th order homogeneous polynomials of $\omega$:
\begin{equation*}
 \sum_{|\alpha|=k} \mathcal{C}^1_\alpha(x_0) u^\alpha, ~
 \sum_{|\alpha|=k} \mathcal{C}^2_\alpha(x_0) u^\alpha,
\end{equation*}
where $u^{\alpha}=u_{1}^{\alpha_{1}}\cdots u_{d}^{\alpha_{d}}$ and
$\alpha=(\alpha_1,  \ldots ,\alpha_d)\in {\mathbb Z}_{+}^{d}$ is a multi-index.
$\phi_0^{(k)}$ and $\psi_0^{(k)}$ are also $k$-th order homogeneous polynomials 
of $ \omega$ 
up to $\mathrm{Ker} (I-L)$ and $\mathrm{Ker} (I-\,^t L)$, respectively. 
Namely, for $\omega=u_{1}\omega_{1}+\cdots +u_{d} \omega_{d}$, 
$$
\phi_0^{(k)} (x_0) = \sum_{ |\alpha|=k}  \phi_\alpha (x_0) u^\alpha+ c_k(\omega) ,  \quad 
\psi_0^{(k)} (x_0) = \sum_{ |\alpha|=k}  \psi_\alpha (x_0) u^\alpha+ d_k(\omega) m (x_0),
$$
where $\phi_\alpha$ and $\psi_\alpha$ are respectively solutions of
\begin{equation*}
(I-L)\phi_\alpha (x_0)= \mathcal{C}^1_\alpha (x_0), \quad 
(I-\,^t L)\psi_\alpha (x_0)=\mathcal{C}^2_\alpha (x_0) \qquad (x\in V_{0}).
\end{equation*}
By the assumption of $\langle \phi_t, \phi_t \rangle_{\ell^{2}(X_{0})}=1$, we have
\begin{align}
0=&
\big( \frac{d}{dt} \big)^{k} \Big \vert_{t=0} 
\langle \phi_t, \phi_t \rangle_{\ell^{2}(X_{0})}
\nonumber \\
=& |V_0|^{-1/2} \sum_{z \in V_0} \phi_0^{(k)} (z)+
\sum_{i =1}^{k-1} \binom{k}{i} \sum_{z \in V_0} \phi_0^{(k-i)}(z) \overline {\phi_0^{(i)}(z) }  
+|V_0|^{-1/2} \sum_{z \in V_0} \overline{ \phi_0^{(k)} (z) }  . \label{<phi,phi>}
\end{align}
Together with (\ref{<phi,psi>}), we find that $c_k(\omega)$ and $d_k(\omega)$ are
$k$-th order homogeneous polynomials of $\sqrt{-1} \omega$.
Hence,  we conclude that $\phi_0^{(k)}(x_0)$ and 
$\psi_0^{(k)}(x_0)$ are $k$-th order homogeneous polynomial of $\sqrt{-1} \omega$.

In the case $k$ is odd, we see that 
\begin{equation*}
(I-L) \Re \phi_0^{(k)} =0, \quad 
(I-\,^t L)\Re \psi_0^{(k)} =0,
 \end{equation*}
which implies that $\Re \phi_0^{(k)}=a_k(\omega)$, 
$\Re \psi_0^{(k)}=b_k(\omega) m$ for some constant function $a_k$, $b_k$
on $V_0$ depending on $\omega \in \mathrm{Hom}(\Gamma, \mathbb{R})$. 
As before, we see that $a_k( s \omega ) = s^k a_k (\omega)$ and 
$b_k (s \omega) = s^k b_k (\omega )$. 
Taking the real part of 
$(d/dt)^{k} \vert_{t=0} \langle \phi_t, \psi_t \rangle_{\ell^{2}(X_{0})}$ 
as above, we obtain $a_k=|V_0| b_k$.
Additionally, using (\ref{<phi,phi>}), we have
$$2 |V_0|^{-1/2} \sum_{ z \in V_0} \Re \phi_0^{(k)}(z)=0,$$
which implies that $a_k=b_k=0$. Then we see that $\phi_0^{(k)}$ and 
$\psi_0^{(k)}$ are imaginary-valued functions on $V_0$ written as
\begin{align*}
\phi_0^{(k)} (x_0) =& \sum_{ |\alpha|=k} \sqrt{-1} a_\alpha \phi_\alpha (x_0) + 
\sqrt{-1} c_k( \omega) , \\
\psi_0^{(k)} (x_0) =& \sum_{ |\alpha|=k} \sqrt{-1} b_\alpha \psi_\alpha (x_0) + 
\sqrt{-1} d_k( \omega) m (x_0)
\end{align*}
for some constant functions $c_k(\omega), d_k(\omega)$ on $V_0$ depending on 
$\omega \in \mathrm{Hom}(\Gamma, \mathbb{R})$. 
As before, it is easy to see that 
$c_k(s \omega)= s^k c_k(\omega)$, $d_k(s \omega)= s^k d_k(\omega)$. 
Then, by replacing eigensection $s_\omega$ with
$
s_\omega \exp \left( \sqrt{-1} \frac{c_k(\omega) }{ k!} |V_0| ^{1/2} \right)
$, 
and eigenfunctions $\phi_\omega$, $\psi_\omega $ with
$\phi_\omega \exp \left(  -\sqrt{-1} \frac{c_k(\omega) }{ k!} |V_0| ^{1/2} \right)$, 
$\psi_\omega  \exp \left(  -\sqrt{-1} \frac{c_k(\omega) }{ k!} |V_0| ^{1/2} \right)$, 
respectively, 
we can assume that $c_k=0$, which implies that $\phi_0^{(k)}(x_0)$ is a $k$-th order 
homogeneous polynomial of $\sqrt{-1} \omega$. 

Finally, from (\ref{<phi,psi>}) and the assumption of the induction, 
we obtain that $d_k$ is a $k$-th order homogeneous 
polynomial, which completes the proof. 
\qed
\subsection{Basic facts on the Fourier transform}
In this subsection, we 
quickly review several basic facts on the Fourier transform based on \cite[Section 4]{KS00}. We define
the function $\mathfrak{F} f (\xi )=\mathfrak{F}_{2\pi^{2}}f(\xi)$, for a rapidly decreasing function $f=f(u)$ 
($u\in \mathbb R^{d}$), by
\begin{equation*}
\mathfrak{F} f (\xi ):= \int_{\mathbb{R}^d} f(u) \exp \Big( -2\pi^2 |u|_{\mathbb R^{d}}^2 -
2\pi \sqrt{-1}  \big(u, \xi \big)_{\mathbb R^{d}}  \Big) du \qquad (\xi \in \mathbb R^{d}). 
\end{equation*}
It is easy to see
\begin{equation}
\mathfrak{F} (1) (\xi ) = (2\pi )^{-d/2}  \exp \big( -\frac{|\xi |_{\mathbb R^{d}}^2 }{2} \big)
\qquad (\xi \in \mathbb R^{d}),
\label{Fourier-basic1}
\end{equation}
and 
\begin{equation*}
u^\alpha \exp \big( - 2\pi \sqrt{-1} (u, \xi)_{\mathbb R^{d}} \big) =
\frac{1}{ (-2\pi \sqrt{-1} )^{|\alpha  | } } 
\partial_\xi ^\alpha \exp \big( - 2\pi \sqrt{-1} (u,\xi)_{\mathbb R^{d}} \big),
\end{equation*}
where $\partial^{\alpha}_{\xi}={\partial}_{1}^{\alpha_{1}}\cdots {\partial}_{d}^{\alpha_{d}}=
(\frac{\partial}{{\partial \xi}_{1}})^{\alpha_{1}}
\cdots 
(\frac{\partial}{{\partial \xi}_{d}})^{\alpha_{d}}
$ for $\xi=(\xi_{1},\ldots, \xi_{d})$, $\alpha=(\alpha_{1}, \ldots, \alpha_{d})\in \mathbb Z_{+}^{d}$.
Differentiating of the left-hand side of (\ref{Fourier-basic1}) with respect to $\xi$ gives
\begin{equation}
\mathfrak{F} (u^\alpha ) (\xi) = \frac{1}{ (2 \pi )^{d/2} ( - 2 \pi \sqrt{-1} )^{|\alpha | } } \partial_\xi^\alpha 
\exp \big(- \frac{  |\xi |_{\mathbb R^{d}}^2} {2} \big) \qquad (\xi \in \mathbb R^{d}).
\label{Fourier-basic2}
\end{equation}
On the other hand, we also have
\begin{align*}
\partial_i  \exp \big( -\frac{|\xi |_{\mathbb R^{d}}^2 }{2} \big) 
&= -\xi_i  \exp \big( -\frac{|\xi |_{\mathbb R^{d}}^2 }{2} \big), \\
\partial_{i} \partial_{j}\exp \big( -\frac{|\xi |_{\mathbb R^{d}}^2 }{2} \big) 
&= (\xi_{i}\xi_{j}-\delta_{ij}) \exp \big( -\frac{|\xi |_{\mathbb R^{d}}^2 }{2} \big), \\
\partial_{i} \partial_{j} \partial_{k} \exp \big( -\frac{|\xi |_{\mathbb R^{d}}^2 }{2} \big) 
&= (-\xi_{i}\xi_{j}\xi_{k}+\delta_{ij}\xi_{k}+\delta_{jk}\xi_{i}+\delta_{ki}\xi_{j}) 
\exp \big( -\frac{|\xi |_{\mathbb R^{d}}^2 }{2} \big), \\
\partial_{i} \partial_{j} \partial_{k} \partial_{l} \exp \big( -\frac{|\xi |_{\mathbb R^{d}}^2 }{2} \big) 
&= \big( \xi_{i}\xi_{j}\xi_{k}\xi_{l}-\delta_{ij}\xi_{k}\xi_{l}-\delta_{jk}\xi_{l}\xi_{i}-\delta_{kl}\xi_{i}\xi_{j}
-\delta_{li}\xi_{j}\xi_{k}
\\
&\qquad 
-\delta_{ik}\xi_{j}\xi_{l}-\delta_{jl}\xi_{i}\xi_{k}+\delta_{ij}\delta_{kl}+\delta_{ik}\delta_{jl}
+\delta_{li}\delta_{jk} \big) 
\exp \big( -\frac{|\xi |_{\mathbb R^{d}}^2 }{2} \big). \\
\end{align*}
Combining these identities with (\ref{Fourier-basic2}), we obtain
\begin{align}
\mathfrak{F} (u_{i}) (\xi) &= -(2 \pi )^{-d/2} \frac{{\sqrt{-1}}
\xi_{i}}{2\pi}
\exp \big(- \frac{  |\xi |_{\mathbb R^{d}}^2} {2} \big),
\nonumber \\
\mathfrak{F} (u_{i}u_{j}) (\xi)  &= -(2 \pi )^{-d/2}
\frac{\xi_{i}\xi_{j}-\delta_{ij}}{4\pi^{2}}
\exp \big(- \frac{  |\xi |_{\mathbb R^{d}}^2} {2} \big),
\nonumber \\
\mathfrak{F} (u_{i}u_{j}u_{k}) (\xi)  &= -(2 \pi )^{-d/2} 
\frac{{\sqrt{-1}}}{8\pi^{3}}
\big( -\xi_{i}\xi_{j}\xi_{k}+\delta_{ij}\xi_{k}+\delta_{jk}\xi_{i}+\delta_{ki}\xi_{j} \big) 
\exp \big( -\frac{|\xi |_{\mathbb R^{d}}^2 }{2} \big), 
\nonumber \\
\mathfrak{F} (u_{i}u_{j}u_{k}u_{l}) (\xi) 
&= 
(2 \pi )^{-d/2} \frac{1}{16\pi^{4}}
\Big( \xi_{i}\xi_{j}\xi_{k}\xi_{l}-\delta_{ij}\xi_{k}\xi_{l}-\delta_{jk}\xi_{l}\xi_{i}-\delta_{kl}\xi_{i}\xi_{j}
-\delta_{li}\xi_{j}\xi_{k}
\nonumber \\
&\qquad 
-\delta_{ik}\xi_{j}\xi_{l}-\delta_{jl}\xi_{i}\xi_{k}+\delta_{ij}\delta_{kl}+\delta_{ik}\delta_{jl}
+\delta_{li}\delta_{jk} \Big) 
\exp \big( -\frac{|\xi |_{\mathbb R^{d}}^2 }{2} \big). 
\nonumber
\end{align}
Repeating this argument twice, we further obtain
\begin{align}
\mathfrak{F} (u_{i}u_{j}u_{k}u_{l}u_{m}u_{n}) (\xi) 
&= -(2 \pi )^{-d/2} \frac{1}{64\pi^{6}} \Big( \xi_{i}\xi_{j} \xi_{k} \xi_{l} \xi_{m} \xi_{n}
-{\mathfrak f}_{4}(i,j,k,l,m,n)(\xi) 
\nonumber \\
&
\hspace{-4mm} +{\mathfrak f}_{2}(i,j,k,l,m,n)(\xi)-
{\mathfrak g}(i,j,k,l,m,n) \Big)
\exp \big( -\frac{|\xi |_{\mathbb R^{d}}^2 }{2} \big),
\nonumber
\end{align}
where ${\mathfrak f}_{r}(i,j,k,l,m,n)(\xi)$ ($r=2,4$) is a homogeneous polynomial of degree $r$
in the variables $\xi_{i},\xi_{j},\ldots, \xi_{n}$ and the constant
${\mathfrak g}(i,j,k,l,m,n)$ is given by  
\begin{align}
{\mathfrak g}(i,j,k,l,m,n)&=\delta_{ij}(\delta_{kl}\delta_{mn}+\delta_{km}\delta_{ln}+\delta_{kn}\delta_{lm})
+\delta_{ik}(\delta_{jl}\delta_{mn}+\delta_{jm}\delta_{ln}+\delta_{jn}\delta_{lm})
\nonumber \\
&\quad +
\delta_{il}(\delta_{jk}\delta_{mn}+\delta_{jm}\delta_{kn}+\delta_{jn}\delta_{km})
+\delta_{im}(\delta_{jk}\delta_{ln}+\delta_{jl}\delta_{kn}+\delta_{jn}\delta_{kl})
\nonumber \\
&\quad +\delta_{in}(\delta_{jk}\delta_{lm}+\delta_{jl}\delta_{km}+\delta_{jm}\delta_{kl}).
\nonumber
\end{align}
\subsection{Proof of Theorem \ref{AE-LCLT}}
Since $p(n,x,y)=0$ for 
$x\in A_{i}, y\in A_{j}$ and $n\neq Kl+j-i$ is obvious, 
we only consider the case $x\in A_{i}, y\in A_{j}$ and $n=Kl+j-i$.
For later use, we take an extension of 
$\lambda_{\omega} (1)=-\log \mu_{0}(\chi_{\omega})$ to the whole space
${\rm{Hom}}(\Gamma, \mathbb R)$ such that 
${\rm Re} \lambda_{\omega}(1) \geq b\Vert \omega \Vert^{2}$
for some constant $b>0$. We also extend 
$\phi_\omega=\phi_{0, \chi_{\omega}}$ and 
$\psi_\omega=\psi_{0, \chi_{\omega}}$ to smooth
compact supported functions on $\mathrm{Hom}(\Gamma, \mathbb{R})$,
and set ${\mathfrak p}(\omega; x,y)=\phi_{\omega}(\pi(x)) {\overline{\psi_{\omega}(\pi(y))}}$.
Since $\mathfrak p$ has compact support in the variable $\omega$, it holds
$$ \Vert {\mathfrak p} \Vert_{\infty}:=\sup \big\{ \vert 
{\mathfrak p}(\omega; x,y) \vert :~\omega\in {\rm Hom}(\Gamma, \mathbb R),~ x,y\in V\big \}<\infty.
$$
We take an orthonormal basis $\{ \omega_{1}, \ldots, \omega_{d} \}$ of 
$\mathrm{Hom}(\Gamma, \mathbb{R}) (\subset {\rm H}^{1}(X_{0}, \mathbb R))$ 
and write $\omega=u_{1}\omega_{1}+
\cdots+u_{d}\omega_{d}$. Then the normalized Haar measure $d\chi$ on $\widehat \Gamma$
is written as 
$$d\chi=d\chi_{\omega}={\rm vol}({\widehat \Gamma})^{-1}d\omega
={\rm vol}({\rm Alb}^{\Gamma}) d\omega,$$
where $d\omega=du=du_{1}\cdots du_{d}$ is the Lebesgue measure on 
$\mathrm{Hom}(\Gamma, \mathbb{R})$.

We divide the proof of Theorem \ref{AE-LCLT} into several steps.
\vspace{2mm} \\
{\bf Step 1.}~Recall the integral representation (\ref{decomp}).
As mentioned in Subsection 6.1, we may choose 
a constant $0<\eta <1$ and
a sufficiently small neighborhood $U(\bf 1)$ of the trivial character $\chi={\bf 1}$
such that $\vert \mu_{0}(\chi) \vert <\eta$ for all
$\chi\in \widehat{\Gamma} \setminus U(\bf 1)$. By
(\ref{decomp}), (\ref{exp-falloff1}) and (\ref{chi-omega-henkan}), we have
\begin{eqnarray*}
p(n,x,y) &=& K {\rm vol}({\rm Alb}^{\Gamma})
\int_{{\rm Hom}(\Gamma, \mathbb R)} \exp(-n\lambda_{\omega}(1))
\phi_{\omega}(\pi(x)) {\overline{\psi_{\omega}(\pi(y))}}
\nonumber \\
&\mbox{ }&
\times
\exp \Big( -2\pi {\sqrt{-1}} \langle \omega, \Phi_{0}(y)-\Phi_{0}(x) \rangle \Big) d\omega
+C{\eta}^{n}
\nonumber \\
&=&
K {\rm vol}({\rm Alb}^{\Gamma})n^{-d/2}
\int_{{\rm Hom}(\Gamma, \mathbb R)} \exp \big ( -n\lambda_{\frac{\omega}{\sqrt n}}(1) \big )
\phi_{\omega/{\sqrt n}}(\pi(x)) {\overline{\psi_{\omega/{\sqrt n}}(\pi(y))}}
\nonumber \\
&\mbox{ }&
\times
\exp \Big( -2\pi {\sqrt{-1}} \langle \frac{\omega}{\sqrt n}, \Phi_{0}(y)-\Phi_{0}(x) \rangle \Big) d\omega
+C{\eta}^{n}
\nonumber \\
&=:& 
K {\rm vol}({\rm Alb}^{\Gamma})n^{-d/2} I(n)+C{\eta}^{n}
\label{AE-proof-1}
\end{eqnarray*}
for some positive constant $C$ independent of $x$ and $y$. 
For simplicity, 
we now set 
$${\mathfrak d}_{n}(x,y; \rho_{\mathbb R}(\gamma_{p}))
:=n^{-1/2}(\Phi_{0}(y)-\Phi_{0}(x)-n{\rho_{\mathbb R}}(\gamma_{p})).$$
By expanding  
$n \lambda_{\omega /{\sqrt n}}(1)$
as
\begin{align}
n \lambda_{\omega /{\sqrt n}}(1)&=n \sum_{k=0}^{4} \frac{1}{k!} 
\lambda^{(k)}_{\omega /{\sqrt n}}(0) +n O  \Big( \frac{ \Vert \omega \Vert^{5} }{n^{5/2}} \Big) 
\nonumber \\
&= n \Big( -2\pi {\sqrt{-1}} \big \langle \gamma_{p}, \frac{\omega}{\sqrt n} \big \rangle 
+2\pi^{2} \Big \Vert \frac{\omega}{\sqrt n} \Big \Vert^{2} \Big)
+n \Big( \frac{1}{6} \lambda^{(3)}_{\omega/{\sqrt n}} (0)+\frac{1}{24} \lambda^{(4)}_{\omega/{\sqrt n}} (0) \Big ) 
\nonumber \\
&\quad
+n O  \Big( \frac{ \Vert \omega \Vert^{5} }{n^{5/2}} \Big),
\nonumber
\end{align}
we have
\begin{eqnarray*}
I(n)&=&
\int_{{\rm Hom}(\Gamma, \mathbb R)}
e^{-2\pi^{2} \Vert \omega \Vert^{2}}
\exp \Big( -\frac{n}{6} \lambda^{(3)}_{\omega/{\sqrt n}} (0)-
\frac{n}{24} \lambda^{(4)}_{\omega/{\sqrt n}} (0) \big ) \Big)
\nonumber \\
&\mbox{  }& \times 
\exp \Big \{ O  \Big( \frac{ \Vert \omega \Vert^{5} }{n^{3/2}} \Big) \Big \}
{\mathfrak p}(\omega/ {\sqrt n};x,y)
\exp \big( -2\pi {\sqrt{-1}}  \langle \omega, {\mathfrak d}_{n}(\gamma_{p};x,y)
\rangle \big) d\omega
\nonumber \\
&=& \int_{\Vert \omega \Vert \leq n^{1/6}} 
%
e^{-2\pi^{2} \Vert \omega \Vert^{2}}
\exp \Big( -\frac{n}{6} \lambda^{(3)}_{\omega/{\sqrt n}} (0)-
\frac{n}{24} \lambda^{(4)}_{\omega/{\sqrt n}} (0) \big ) \Big)
\nonumber \\
&\mbox{  }& \qquad \times 
\exp \Big \{ O  \Big( \frac{ \Vert \omega \Vert^{5} }{n^{3/2}} \Big) \Big \}
{\mathfrak p}(\omega/ {\sqrt n};x,y)
\exp \big( -2\pi {\sqrt{-1}}  \langle \omega, {\mathfrak d}_{n}(\gamma_{p};x,y)
\rangle \big) d\omega
\nonumber \\
&\mbox{ }&+
\int_{\Vert \omega \Vert > n^{1/6}}
\exp \big ( -n\lambda_{\frac{\omega}{\sqrt n}}(1) \big )
{\mathfrak p}(\omega/ {\sqrt n};x,y)
\nonumber \\
&\mbox{ }&
\qquad
\times
\exp \Big( -2\pi {\sqrt{-1}} \langle \frac{\omega}{\sqrt n}, \Phi_{0}(y)-\Phi_{0}(x) \rangle \Big) d\omega
\nonumber \\
&=:& I_{1}(n)+I_{2}(n).
\end{eqnarray*}
Recalling $n {\rm Re} \lambda_{\omega/{\sqrt n}}(1) \geq b \Vert \omega \Vert^{2}$ 
($\omega \in {\rm Hom}(\Gamma, \mathbb R)$), we obtain
\begin{equation*}
\vert I_{2}(n) \vert \leq \Vert {\mathfrak p} \Vert_{\infty} \int_{\Vert \omega \Vert > n^{1/6}}
e^{-b\Vert \omega \Vert^{2}} d\omega
\label{I2-expfall}
\end{equation*}
Thus $I_{2}(n)$ converges to $0$ as $n\to \infty$ exponentially fast uniformly for $x,y\in V$.
\vspace{2mm} \\
{\bf Step 2.} We divide the integral $I_{1}(n)$ into
\begin{eqnarray*}
I_{1}(n)&=&
\int_{\Vert \omega \Vert \leq n^{1/6}} 
%
e^{-2\pi^{2} \Vert \omega \Vert^{2}}
\exp \Big( -\frac{n}{6} \lambda^{(3)}_{\omega/{\sqrt n}} (0)-
\frac{n}{24} \lambda^{(4)}_{\omega/{\sqrt n}} (0) \big ) \Big)
\nonumber \\
&\mbox{  }& \qquad \times 
{\mathfrak p}(\omega/ {\sqrt n};x,y)
\exp \big( -2\pi {\sqrt{-1}}  \langle \omega, {\mathfrak d}_{n}(\gamma_{p};x,y)
\rangle \big) d\omega
\nonumber \\
&\mbox{ }&
+
\int_{\Vert \omega \Vert \leq n^{1/6}} 
%
e^{-2\pi^{2} \Vert \omega \Vert^{2}}
\exp \Big( -\frac{n}{6} \lambda^{(3)}_{\omega/{\sqrt n}} (0)-
\frac{n}{24} \lambda^{(4)}_{\omega/{\sqrt n}} (0) \big ) \Big)
{\mathfrak p}(\omega/ {\sqrt n};x,y)
\nonumber \\
&\mbox{  }& \qquad \times 
\Big( \exp \Big \{ O  \Big( \frac{ \Vert \omega \Vert^{5} }{n^{3/2}} \Big) \Big \}-1 \Big)
\exp \big( -2\pi {\sqrt{-1}}  \langle \omega, {\mathfrak d}_{n}(\gamma_{p};x,y)
\rangle \big) d\omega
\nonumber \\
&=:& I_{3}(n)+I_{4}(n).
\end{eqnarray*}
Thanks to Lemma \ref{lemma of expansion}, we find 
$$
n\lambda_{\omega/{\sqrt n}}^{(3)}(0)=
n O \Big( \frac{\Vert \omega \Vert}{n^{1/2}} \Big)^{3}=O(1), \quad 
n\lambda_{\omega/{\sqrt n}}^{(4)}(0)=n O \Big( \frac{\Vert \omega \Vert}{n^{1/2}} \Big)^{4}
=O(n^{-1/2})
$$
and
\begin{align*}
\big \vert 
\exp \Big \{ O  \Big( \frac{ \Vert \omega \Vert^{5} }{n^{3/2}} \Big) \Big \}-1
\big \vert &=
\exp \Big \{ O  \Big( \frac{ \Vert \omega \Vert^{5} }{n^{3/2}} \Big) \Big \}
O  \Big( \frac{ \Vert \omega \Vert^{5} }{n^{3/2}} \Big)
\nonumber \\
&= \exp \big( O(n^{-2/3}) \big ) 
O  \Big( \frac{ \Vert \omega \Vert^{5} }{n^{3/2}} \Big).
\end{align*}
provided $\Vert \omega \Vert \leq n^{1/6}$. 
Thus we can estimate $I_{4}(n)$ as
\begin{equation*}
\vert I_{4}(n) \vert \leq C \Vert {\mathfrak p} \Vert_{\infty} 
\int_{ {\rm Hom}(\Gamma, \mathbb R)}
O  \Big( \frac{ \Vert \omega \Vert^{5} }{n^{3/2}} \Big)
e^{-2\pi^{2}\Vert \omega \Vert^{2}} d\omega =O(n^{-3/2}).
\end{equation*}
%
{\bf Step 3.}~As a direct consequence of Lemma \ref{lemma of expansion}, we easily obtain
\begin{lm} \label{q-coefficient}
There exist real constants ${\mathfrak q}_{\alpha}={\mathfrak q}_{\alpha}(\pi(x),\pi(y);\gamma_{p})$
{\rm{(}}$\alpha=(\alpha_{1}, \ldots, \alpha_{r}) \in \{1,\ldots, d\}^{r}, \vert \alpha \vert=r=1,\ldots,4${\rm{)}} such that
${\mathfrak q}_{\sigma (\alpha)}={\mathfrak q}_{\alpha}$ 
{\rm{(}}$\sigma \in {\cal S}_{r}${\rm{)}}
and
\begin{align} 
&
\vert V_{0} \vert^{-1/2} {\overline{\psi_{\omega(u)/\sqrt n}'(\pi(y))}}
=
\frac{{\sqrt{-1}}}{n^{1/2}}\sum_{i=1}^{d} {\mathfrak q}_{i} u_{i}
=: \frac{{\sqrt{-1}}}{n^{1/2}} {\cal Q}_{1}(u),
\nonumber \\
&
|V_0|^{1/2} \phi''_{\omega(u)/\sqrt n}
(\pi(x)) +|V_0|^{-1/2} {\overline{\psi^{''}_{\omega(u)/\sqrt n} (\pi(y)) }} 
\nonumber \\
&\qquad \qquad \qquad \qquad \qquad 
=
\Big( \frac{{\sqrt{-1}}}{n^{1/2}} \Big)^{2} \sum_{i,j=1}^{d} {\mathfrak q}_{ij} u_{i}u_{j}
=:
\Big( \frac{{\sqrt{-1}}}{n^{1/2}} \Big)^{2} {\mathcal Q}_{2}(u),
\nonumber 
\end{align}
\begin{align}
&
\lambda_{\omega(u)/n^{1/2}}^{(3)}(0)
=\Big( \frac{{\sqrt{-1}}}{n^{1/2}} \Big)^{3} \sum_{i,j,k=1}^{d} {\mathfrak q}_{ijk} u_{i}u_{j}u_{k}
=:
\Big( \frac{{\sqrt{-1}}}{n^{1/2}} \Big)^{3} {\mathcal Q}_{3}(u),
\nonumber \\
&
\lambda_{\omega(u)/n^{1/2}}^{(4)}(0)
=
\Big( \frac{{\sqrt{-1}}}{n^{1/2}} \Big)^{4} \sum_{i,j,k,l=1}^{d} {\mathfrak q}_{ijkl} u_{i}u_{j}u_{k}u_{l}
=:
\Big( \frac{{\sqrt{-1}}}{n^{1/2}} \Big)^{4} {\cal Q}_{4}(u)
\nonumber 
\end{align}
hold for all $\omega=\omega(u)=u_{1}\omega_{1}+\cdots +u_{d}\omega_{d} \in 
{\rm Hom}(\Gamma, \mathbb R)$.
\end{lm}

Then, under the condition $\Vert \omega \Vert \leq n^{1/6}$, we can expand as
\begin{eqnarray}
\lefteqn{
\exp \Big( -\frac{n}{6} \lambda^{(3)}_{\omega/{\sqrt n}} (0)-
\frac{n}{24} \lambda^{(4)}_{\omega/{\sqrt n}} (0) \big ) \Big)
}
\nonumber \\
&=&
1 -\Big( \frac{ \lambda^{(3)}_{\omega/{\sqrt n}}(0)} {6}+
\frac{ \lambda^{(4)}_{\omega/{\sqrt n}}(0)} {24} \Big)n 
+\frac{1}{2} \Big (
\frac{ \lambda^{(3)}_{\omega/{\sqrt n}}(0)} {6}+
\frac{ \lambda^{(4)}_{\omega/{\sqrt n}}(0)} {24} \Big)^{2}
n^2 
\nonumber \\
&\mbox{ }&
+
\frac{1}{6}\exp \big( O(1)+O(n^{-2/3}) \big ) 
\Big (
\frac{ \lambda^{(3)}_{\omega/{\sqrt n}}(0)} {6}+
\frac{ \lambda^{(4)}_{\omega/{\sqrt n}}(0)} {24} \Big)^{3}
n^3 
\nonumber \\
&=&1+\big( \frac{\sqrt{-1}}{n^{1/2}}\big) \frac{{\mathcal Q}_{3}(u)}{6}
+\big( \frac{\sqrt{-1}}{n^{1/2}}\big)^{2} \Big( \frac{ {\mathcal Q}_{4}(u)}{24}
+\frac{{\mathcal Q}_{3}(u)^{2}}{72} \Big)
+
O\big( \frac{\vert u \vert_{\mathbb R^{d}}^{9}}{n^{3/2}} \big),
\label{exp2-expansion}
\end{eqnarray}
and
\begin{eqnarray}
{\mathfrak p}(\omega/{\sqrt n};x,y)&=&
\big (\phi_{0}(\pi(x))+\phi'_{0}(\pi(x))+\cdots \big) \big(
{\overline{\psi_{0}(\pi(y))}}+{\overline{\psi'_{0}(\pi(y))}}+\cdots \big)
\nonumber \\
&=&
m(\pi(y))+|V_0|^{-1/2} \overline{ \psi_{\omega/{\sqrt n}}^\prime (\pi(y))} 
\nonumber \\
&\mbox{  }&
+\frac{1}{2} \left( |V_0|^{1/2} \phi_{\omega/{\sqrt n}}^{(2)} (\pi(x)) 
m(\pi(y)) +|V_0|^{-1/2}  {\overline{\psi_{\omega/{\sqrt n}}^{(2)} (\pi(y))}}
\right) 
\nonumber \\
&\mbox{  }&
+O\Big( \frac{ \Vert \omega \Vert^{3}}{n^{3/2}} \Big)
\nonumber \\
&=& 
m(\pi(y))+\big( \frac{\sqrt{-1}}{n^{1/2}}\big) {\mathcal Q}_{1}(u)+ 
\big( \frac{\sqrt{-1}}{n^{1/2}}\big)^{2} \frac{ {\mathcal Q}_{2}(u)}{2}
+O\Big( \frac{ \vert u \vert^{3}_{\mathbb R^{d}}}{n^{3/2}} \Big).
\label{p-expansion}
\end{eqnarray}
%
Multiplying (\ref{exp2-expansion}) and (\ref{p-expansion}) together, 
we obtain
\begin{align*}
I_{3}(n)&= \int_{\vert u \vert_{\mathbb R^{d}}\leq n^{1/6}} 
\sum_{i=0}^{2}
\big (\frac{{\sqrt -1}}{n^{1/2}} \big)^{i}
{\mathcal A}_{i}(u)
\nonumber \\
&\quad \times
\exp \Big( -2\pi^{2} \vert u \vert_{\mathbb R^{d}}^{2} 
-2\pi {\sqrt{-1}}  \big( u, {\mathfrak d}_{n}(x,y; \rho_{\mathbb R}(\gamma_{p}))
\big)_{\mathbb R^{d}} \Big) du +O(n^{-3/2}),
\end{align*}
where 
\begin{align*}
{\cal A}_{0}(u)&=m(\pi(y)),
\nonumber \\
{\cal A}_{1}(u)&={\cal Q}_{1}(u)+\frac{m(\pi(y))}{6}{\cal Q}_{3}(u),
\nonumber \\
{\cal A}_{2}(u)&=\frac{1}{2}{\cal Q}_{2}(u)+\frac{1}{6}{\cal Q}_{1}(u) {\cal Q}_{3}(u)
+m(\pi(y)) \big( \frac{{\cal Q}_{4}(u)}{24}+ \frac{{\cal Q}_{3}(u)^{2}}{72} \big).
\end{align*}
We also find that the integral
$$
\int_{\vert  u \vert_{\mathbb R^{d}}>n^{1/6}} 
\sum_{i=0}^{2}
\big (\frac{{\sqrt -1}}{n^{1/2}} \big)^{i } {\cal A}_{i}(u)
\exp \Big( -2\pi^{2} \vert u \vert_{\mathbb R^{d}}^{2} 
-2\pi {\sqrt{-1}}  \big( u, {\mathfrak d}_{n}(x,y; \rho_{\mathbb R}(\gamma_{p}))
\big)_{\mathbb R^{d}} \Big) du
$$
converges to $0$ as $n\to \infty$ exponentially fast in the same way as Step 1.

Putting it all together, we now obtain
\begin{eqnarray}
I(n)&=&
 \int_{\mathbb R^{d}} 
\sum_{i=0}^{2}
\big (\frac{{\sqrt -1}}{n^{1/2}} \big)^{i} {\cal A}_{i}(u)
\nonumber \\
&\mbox{ }& \times
\exp \Big( -2\pi^{2} \vert u \vert_{\mathbb R^{d}}^{2} 
-2\pi {\sqrt{-1}}  \big( u, {\mathfrak d}_{n}(x,y; \rho_{\mathbb R}(\gamma_{p}))
\big)_{\mathbb R^{d}} \Big) du +O(n^{-3/2})
\nonumber \\
&=& \sum_{i=0}^{2} \big (\frac{{\sqrt -1}}{n^{1/2}} \big)^{i}
{\mathfrak F}({\cal A}_{i})({\mathfrak d}_{n}(x,y; \rho_{\mathbb R}(\gamma_{p})))+O(n^{-3/2}).
\label{final-remark}
\end{eqnarray}
{\bf Step 4.}  We are now in a position to calculate 
$\big (\frac{{\sqrt -1}}{n^{1/2}} \big)^{i}{\mathfrak F}({\cal A}_{i})({\mathfrak d}_{n}(x,y; \rho_{\mathbb R}(\gamma_{p})))$
($n=0,1,2$). It follows directly from (\ref{Fourier-basic1}) that
\begin{equation*}
 \mathfrak{F} ({\cal A}_{0}) \big( {\mathfrak d}_{n}(x,y; \rho_{\mathbb R}(\gamma_{p})) \big)
 =  (2\pi)^{-d/2} m(\pi(y)) 
\exp \left( -\frac{| \Phi_{0}(y)-\Phi_{0}(x)
-n{\rho}_{\mathbb R}(\gamma_{p}) |_{g_{0}}^2 }{2n} \right).
\end{equation*}
Using (\ref{Fourier-basic2}) and noting the condition 
$| \Phi_{0}(y)-\Phi_{0}(x)
-n{\rho}_{\mathbb R}(\gamma_{p}) |_{g_{0}} \leq Cn^{1/6}$, we have
\begin{align}
\big (\frac{{\sqrt -1}}{n^{1/2}} \big) &
\mathfrak{F} ({\cal A}_{1}) \big( {\mathfrak d}_{n}(x,y; \rho_{\mathbb R}(\gamma_{p})) \big)
\nonumber \\
&= (2\pi)^{-d/2} m(\pi(y))
\exp \left( -\frac{| \Phi_{0}(y)-\Phi_{0}(x)
-n{\rho}_{\mathbb R}(\gamma_{p}) |_{g_{0}}^2 }{2n} \right)
\nonumber \\
&\quad \times \Big \{ \Big ( \frac{1}{2\pi} 
m(\pi(y))^{-1}
\sum_{i=1}^{d} {\mathfrak q}_{i} 
\big( \Phi_{0}(y)-\Phi_{0}(x)-\rho_{\mathbb R}(\gamma_{p}) \big)_{i} 
\nonumber \\
&\qquad \quad +  \frac{1}{16\pi^{3}} \sum_{i,j=1}^{d} {\mathfrak q}_{iij}
\big( \Phi_{0}(y)-\Phi_{0}(x)-\rho_{\mathbb R}(\gamma_{p}) \big)_{j} \Big ) n^{-1}+O(n^{-3/2}) \Big \}.
\nonumber
\end{align}
In the same way, 
we also obtain 
\begin{align}
\big (\frac{{\sqrt -1}}{n^{1/2}} \big)^{2} &
\frac{1}{2} \mathfrak{F} ({\cal Q}_{2}) \big( {\mathfrak d}_{n}(x,y; \rho_{\mathbb R}(\gamma_{p})) \big)
\nonumber \\
&=(2\pi)^{-d/2} m(\pi(y))
\exp \left( -\frac{| \Phi_{0}(y)-\Phi_{0}(x)
-n{\rho}_{\mathbb R}(\gamma_{p}) |_{g_{0}}^2 }{2n} \right)
\nonumber \\
&\qquad \times \Big \{ \big(-\frac{m(\pi(y))^{-1}}{8\pi^{2}}
\sum_{i=1}^{d} {\mathfrak q}_{ii}  \big)n^{-1}+O(n^{-5/3}) \Big \},
\nonumber
\end{align}
\begin{align}
\big (\frac{{\sqrt -1}}{n^{1/2}} \big)^{2} &
\frac{1}{6} \mathfrak{F} ({\cal Q}_{1}{\cal Q}_{3}) \big( {\mathfrak d}_{n}(x,y; \rho_{\mathbb R}(\gamma_{p})) \big)
\nonumber \\
&=
(2\pi)^{-d/2} m(\pi(y))
\exp \left( -\frac{| \Phi_{0}(y)-\Phi_{0}(x)
-n{\rho}_{\mathbb R}(\gamma_{p}) |_{g_{0}}^2 }{2n} \right)
\nonumber \\
&
\qquad \times \Big \{ \big(-\frac{m(\pi(y))^{-1}}{32\pi^{4}}\sum_{i,j=1}^{d} {\mathfrak q}_{i}{\mathfrak q}_{ijj}  
\big)n^{-1}+O(n^{-5/3}) \Big \},
\nonumber \\
\big (\frac{{\sqrt -1}}{n^{1/2}} \big)^{2} &
\frac{m(\pi(y))}{24} \mathfrak{F} ({\cal Q}_{4}) \big( {\mathfrak d}_{n}(x,y; \rho_{\mathbb R}(\gamma_{p})) \big)
\nonumber \\
&=
(2\pi)^{-d/2} m(\pi(y))
\exp \left( -\frac{| \Phi_{0}(y)-\Phi_{0}(x)
-n{\rho}_{\mathbb R}(\gamma_{p}) |_{g_{0}}^2 }{2n} \right)
\nonumber \\
&
\qquad \times
 \Big \{  \big( -\frac{1}{128\pi^{4}} \sum_{i,j=1}^{d} {\mathfrak q}_{iijj} \big)
n^{-1}+O(n^{-5/3}) \Big \},
\nonumber 
\\
\big (\frac{{\sqrt -1}}{n^{1/2}} \big)^{2} &
\frac{m(\pi(y))}{72} \mathfrak{F} ({\cal Q}_{3}^{2}) \big( {\mathfrak d}_{n}(x,y; \rho_{\mathbb R}(\gamma_{p})) \big)
\nonumber \\
&=
(2\pi)^{-d/2} m(\pi(y))
\exp \left( -\frac{| \Phi_{0}(y)-\Phi_{0}(x)
-n{\rho}_{\mathbb R}(\gamma_{p}) |_{g_{0}}^2 }{2n} \right)
\nonumber \\
&
\qquad \times \Big \{  \big( -\frac{5}{1536 \pi^{6}} \sum_{i,j,k=1}^{d} {\mathfrak q}_{iij} {\mathfrak q}_{jkk} \big)
n^{-1}+O(n^{-5/3}) \Big \}.
\nonumber
\end{align}

Recalling $m(\pi(y))=m(y)$ and summarizing all above arguments, we complete the 
proof of Theorem \ref{AE-LCLT}. Moreover, we also obtain the explicit expression of
the coefficient $a_{1}=a_{1}\big (\pi(x),\pi(y), \gamma_{p}; \Phi_{0}(y)-\Phi_{0}(x) 
-n{\rho}_{\mathbb R}(\gamma_{p}) 
\big)$ in (\ref{1-AE}) as follows:
\begin{tm}\label{a_1}
\begin{align*}
& a_{1}\big (\pi(x),\pi(y), \gamma_{p}; \Phi_{0}(y)-\Phi_{0}(x) 
-n{\rho}_{\mathbb R}(\gamma_{p}) 
\big)
\nonumber \\
&=\frac{m(\pi(y))^{-1}}{2\pi} \sum_{i=1}^{d} {\mathfrak q}_{i} \big( 
\Phi_{0}(y)-\Phi_{0}(x)-n{\rho}_{\mathbb R}(\gamma_{p}) \big)_{i} 
\nonumber \\
&\quad 
+\frac{1}{16\pi^{3}} \sum_{i,j=1}^{d} {\mathfrak q}_{ij} \big( 
\Phi_{0}(y)-\Phi_{0}(x)-n{\rho}_{\mathbb R}(\gamma_{p}) \big)_{j} 
-\frac{m(\pi(y))^{-1}}{8\pi^{2}} \sum_{i=1}^{d} {\mathfrak q}_{ii}
\nonumber \\
&\quad
-\frac{m(\pi(y))^{-1}}{32\pi^{4}} \sum_{i,j=1}^{d} {\mathfrak q}_{i}{\mathfrak q}_{ijj}
-\frac{1}{128\pi^{4}} \sum_{i,j=1}^{d} {\mathfrak q}_{iijj}-\frac{5}{1536\pi^{6}} \sum_{i,j,k=1}^{d} {\mathfrak q}_{iij}
{\mathfrak q}_{jkk},
\end{align*}
where the coefficients ${\mathfrak q}_{\alpha}={\mathfrak q}_{\alpha}(\pi(x),\pi(y);\gamma_{p})$
{\rm{(}}$\alpha=(\alpha_{1}, \ldots, \alpha_{r}) \in \{1,\ldots, d\}^{r}, r=1,\ldots,4${\rm{)}}
are given in Lemma {\rm{\ref{q-coefficient}}}.
\end{tm}
\begin{re} As we observed in Subsection {\rm{6.1}}, the coefficients ${\mathfrak q}_{\alpha}$ 
{\rm{(}}$\vert \alpha \vert=1,3${\rm{)}} are equal to $0$ provided the random walk is 
$m$-symmetric {\rm{(}}i.e., $\gamma_{p}=0${\rm{)}}. In that case, we have
\begin{align*}
& a_{1}\big (\pi(x),\pi(y); \Phi_{0}(y)-\Phi_{0}(x) \big)
\nonumber \\
&=\frac{1}{16\pi^{3}} \sum_{i,j=1}^{d} {\mathfrak q}_{ij} \big( 
\Phi_{0}(y)-\Phi_{0}(x) \big)_{j} 
-\frac{m(\pi(y))^{-1}}{8\pi^{2}} \sum_{i=1}^{d} {\mathfrak q}_{ii}
-\frac{1}{128\pi^{4}} \sum_{i,j=1}^{d} {\mathfrak q}_{iijj}. \end{align*}
\end{re}

Before closing this subsection, we should mention that 
\begin{equation}
I(n)={\mathfrak F}({\cal A}_{0})({\mathfrak d}_{n}(x,y; \rho_{\mathbb R}(\gamma_{p})))+O(n^{-1/2})
\label{final-remark2}
\end{equation}
also holds uniformly for $x,y\in V$ in Step 3 of the above proof of Theorem \ref{AE-LCLT}. 
Then by replacing
(\ref{final-remark}) by (\ref{final-remark2}), 
we easily obtain the following LCLT:
\begin{co}[Sunada \cite{Sunada-Lecture}] \label{co-sunada}
Suppose that the random walk $\{ w_{n} \}_{n=0}^{\infty}$ on $X$ is irreducible with period $K$. 
Let $V=\coprod_{k=0}^{K-1} A_k$ be the corresponding $K$-partition of $V$. 
Then for any $x\in A_i$ and $y \in A_j$ {\rm{(}}$0\leq i,j \leq K-1${\rm{)}}, 
we have
\begin{align*}
p(n,x,y)= 0 \qquad (n\neq Kl+j-i), 
\end{align*}
and
\begin{align}
(2\pi n)^{d/2} p(n ,x,y)m(y)^{-1}
-
K {\rm vol}({\rm{Alb}}^{\Gamma})
\exp \Big( -\frac{\big \vert \Phi_{0}(y)-\Phi_{0}(x) 
-n{\rho}_{\mathbb R}(\gamma_{p}) \big \vert^{2}_{g_{0}}}
{2n} \Big)  \to 0
\nonumber
\label{Sunada-LCLT}
\end{align}
as $n=Kl+j-i \rightarrow \infty$ uniformly for $x$ and $y$. 
\end{co}
\subsection{Application: 
Another approach to the CLT of the first kind}
In this subsection, applying the LCLT (Corollary \ref{co-sunada}),
we give another proof of (\ref{Trotter-1-2}) in Theorem \ref{CLT-1}
under the irreducibility condition of the random walk on $X$.
For simplicity of the argument, we only consider the case $K=1$. 
(In other cases, the proof goes through in a very similar way with a
slight modification.)

First of all, we take a sequence $\{(x_n, \mathbf{z}_n)\}_{n=1}^\infty $ in 
$V\times \mathrm{H}_1(X_0, \mathbb{R})$ satisfying
\begin{equation*}
\lim_{n\rightarrow \infty} n^{-1/2}
\left( \Phi_0 (x_n)-\rho_\mathbb{R} (\mathbf{z}_n) \right)=
\mathbf{x} \in 
\Gamma \otimes \mathbb{R}.
\end{equation*}
By (\ref{228-2015}), it is enough to show that for any 
$f\in C^{\infty}_{0}(\Gamma \otimes \mathbb R)$,
\begin{equation*}
\lim_{n\rightarrow \infty }
\mathcal{L}_{\gamma_p}^{[nt]} 
\mathcal{P}_{n^{-1/2}} f (x_n, \mathbf{z}_n)= e^{-t\Delta/2} f( \mathbf{x}),
\end{equation*}
that is,
\begin{align*}
\lim_{n\rightarrow \infty}
\sum_{y \in V} p([nt], x_n, y) &
f\left( n^{-1/2} \left( \Phi_0(y)-\rho_{\mathbb{R}} (\mathbf{z}_n) 
-[nt] \rho_{\mathbb{R}} (\gamma_{p} ) \right) \right) \\
&=\int_{\Gamma \otimes \mathbb{R}} \frac{1}{(2\pi t)^{d/2}} 
\exp \left( -\frac{|\mathbf{x}-\mathbf{y}|^2}{2t} \right) f(\mathbf{y}) d \mathbf{y}.
\end{align*}

We set
$$ {\bf{x}}_{n}:=n^{-1/2}
\big( \Phi_{0}(x_{n})-\rho_{\mathbb R}({\bf{z}}_{n}) \big),
\quad {\bf{y}}_{n}(y):=
n^{-1/2} \big( \Phi_{0}(y)-\rho_{\mathbb R}({\bf{z}}_{n}) -[nt]\rho_{\mathbb R}(\gamma_{p})
\big)
\quad (y\in V).
$$
By the LCLT (\ref{Sunada-LCLT}), for an arbitrary small $\varepsilon>0$,
there is a positive integer $n_{0}$ such that for any $n \geq n_0$,
\begin{equation*}
\Big \vert (2\pi [nt])^{d/2} p([nt], x_{n},y) m(y)^{-1}
-{\rm{vol}}({\rm{Alb}}^{\Gamma})
\exp \Big( -\frac{ n\vert  
\mathbf{x}_n-\mathbf{y}_n(y)
\vert_{g_{0} } ^{2}} {2[nt]}  \Big) 
\Big \vert
<\varepsilon
\end{equation*}
uniformly for $x_n, y \in V$. 
Dividing both sides by $(2\pi [nt])^{d/2} m(y)^{-1}$, for any $\varepsilon^\prime>0$, 
there exists $n_1 \in \mathbb{N}$ such that for any $n \geq n_1$, 
\begin{equation*}
\Big \vert  p([nt], x_{n},y) 
-\frac{{\rm{vol}}({\rm{Alb}}^{\Gamma}) m(y)}{ (2\pi [nt])^{d/2}}
\exp \Big( -\frac{ n\vert  
\mathbf{x}_n-\mathbf{y}_n(y)
\vert_{g_{0} } ^{2}} {2[nt]}  \Big) 
\Big \vert
<\varepsilon^\prime n^{-d/2}.
\end{equation*}
Then we obtain
\begin{align*}
&\Big \vert \sum_{y \in V} 
\Big\{ p([nt], x_n , y) f ( \mathbf{y}_n (y)) 
-\frac{{\rm{vol}}({\rm{Alb}}^{\Gamma}) m(y)}{ (2\pi [nt])^{d/2}}
\exp \Big( -\frac{ n\vert  
\mathbf{x}_n-\mathbf{y}_n(y)
\vert_{g_{0} } ^{2}} {2[nt]}  \Big) f( \mathbf{y}_n (y)) \Big\}  \Big\vert \\
&\leq \sum_{y \in V}
\Big \vert  p([nt], x_{n},y) 
-\frac{{\rm{vol}}({\rm{Alb}}^{\Gamma}) m(y)}{ (2\pi [nt])^{d/2}}
\exp \Big( -\frac{ n\vert  
\mathbf{x}_n-\mathbf{y}_n(y)
\vert_{g_{0} } ^{2}} {2[nt]}  \Big) 
\Big \vert | f( \mathbf{y}_n (y)) \vert \\
&\leq \varepsilon^\prime n^{-d/2}
\sum_{y \in V} \vert f(\mathbf{y}_n( y) ) \vert.
\end{align*}
Since the support of $f$ is compact, 
\begin{eqnarray}
n^{-d/2}\sum_{y \in V} \vert f(\mathbf{y}_n( y)) \vert
&=& n^{-d/2} \sum_{ y_0 \in \mathcal{F}} \sum_{\sigma \in \Gamma} 
\left \vert f\Big( \mathbf{y}_n(y_0) +\frac{\sigma}{\sqrt{n}} \Big) \right\vert 
\nonumber \\
&\rightarrow&
\frac{\vert V_0 \vert}{\mathrm{vol}(\mathrm{Alb}^\Gamma )} 
\int_{\Gamma \otimes \mathbb{R}} \vert f(\mathbf{y}) \vert d \mathbf{y}
\qquad \mbox{as }n\to \infty.
\nonumber 
\end{eqnarray}
Then we have
\begin{eqnarray}
\lefteqn{
\lim_{n\rightarrow \infty}
\mathcal{L}_{\gamma_p}^{[nt]} \mathcal{P}_{n^{-1/2}}  f ( x_n, \mathbf{z}_n) }
\nonumber \\
&=& \lim_{n \rightarrow \infty } 
\frac{1}{(2\pi [nt])^{d/2}} \sum_{y \in V} m(y) \mathrm{vol} ( \mathrm{Alb}^\Gamma) 
\exp\Big( - \frac{n}{2[nt]} \vert \mathbf{y}_n(y) - \mathbf{x}_n \vert_{g_0}^2 \Big)
f( \mathbf{y}_n(y))
\nonumber \\
&=&\lim_{n \rightarrow \infty} 
\frac{1}{(2 \pi n t )^{d/2}}
\sum_{y \in V} m(y) \mathrm{vol} ( \mathrm{Alb}^\Gamma) 
\exp\Big( - \frac{1}{2t} \vert \mathbf{y}_n(y) - \mathbf{x}_n \vert_{g_0}^2 \Big)
f( \mathbf{y}_n(y)).
\nonumber
\end{eqnarray}
By \cite[pp. 655]{KS00}, for any $\varepsilon >0$, there exists $n_2 \in \mathbb{N}$ 
such that for any $n \geq n_2$, 
\begin{equation*}
\Big\vert \exp\Big( - \frac{1}{2t} \vert \mathbf{y}_n(y) - \mathbf{x}_n \vert_{g_0}^2 \Big)
-\exp\Big( - \frac{1}{2t} \vert \mathbf{y}_n(y) - \mathbf{x} \vert_{g_0}^2 \Big)
\Big\vert < \varepsilon
\end{equation*}
uniformly for $y \in V$. 

Thus we conclude
\begin{eqnarray}
\lefteqn{
\lim_{n\rightarrow \infty}  
\mathcal{L}_{\gamma_p}^{[nt]} \mathcal{P}_{n^{-1/2}}  f ( x_n, \mathbf{z}_n) }
\nonumber \\
&= &
\lim_{n \rightarrow \infty } 
\frac{1}{(2\pi nt)^{d/2}} \sum_{y \in V} m(y) \mathrm{vol} ( \mathrm{Alb}^\Gamma) 
\exp\Big( - \frac{1}{2t} \vert \mathbf{y}_n(y) - \mathbf{x} \vert_{g_0}^2 \Big)
f( \mathbf{y}_n(y)) 
\nonumber \\
&=&
\frac{1}{(2 \pi t)^{d/2}} \lim_{n \rightarrow \infty}
\sum_{ y_0 \in \mathcal{F}} m(y_0) \sum_{\sigma \in \Gamma} 
\frac{ \mathrm{vol} (\mathrm{Alb}^\Gamma) }{n^{d/2}} 
\nonumber \\
&\mbox{ }&
\quad
\times
f\Big( \mathbf{y}_n(y_0) + \frac{\sigma}{\sqrt{n}} \Big) 
\exp \Big ( - \frac{\vert \mathbf{x}- ( \mathbf{y}_n (y_0) + \frac{\sigma}{\sqrt{n}} ) 
\vert_{g_0}^2 }{2n} \Big)
\nonumber \\
&=& \frac{1}{(2\pi t)^{d/2}} 
\int_{\Gamma \otimes \mathbb{R}} f(\mathbf{y}) 
\exp \Big( -\frac{ \vert \mathbf{x}-\mathbf{y} \vert_{g_0}^2 }{2t} \Big) d \mathbf{y},
\nonumber
\end{eqnarray}
whence (\ref{Trotter-1-2}) follows.
\section{Examples of the modified 
standard realization}\label{section Example}
In this final section, we give several examples of the modified standard realization 
of crystal lattices associated with non-symmetric random walks.
See \cite{KSS, KS00, KS00-2} for the symmetric case. Here we 
write
${\hat{\alpha}}:=\alpha+\alpha^\prime$ and $\check{\alpha}:=\alpha-\alpha^\prime$ for two numbers 
$\alpha, \alpha^\prime$, and denote 
$\sigma \otimes 1 \in \Gamma \otimes \mathbb R$ by
the same symbol $\sigma$ for simplicity of notation.
\subsection{The $2$-dimensional square lattice}
\begin{figure}[tbph]
\begin{center}
\scalebox{1.12}{\input{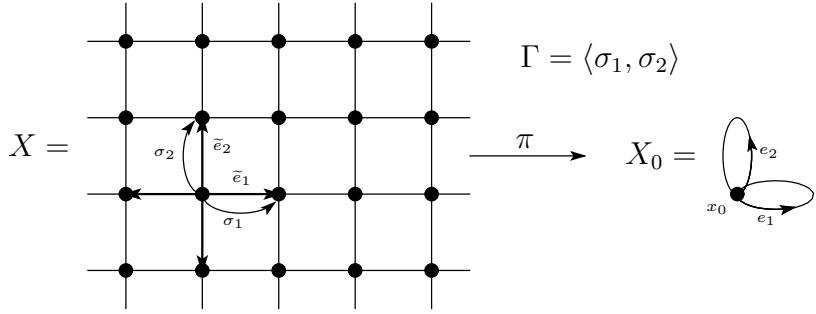}}
\end{center}
\caption{Square lattice graph and its quotient}
\label{figure:square1}
\end{figure}

Let $X=(V,E)$ be the $2$-dimensional square lattice graph. Namely, $V=\mathbb{Z}^2$ and 
\begin{equation*}
E=\{ ({\mathbf{x}},{\mathbf{y}}) \in V\times V \, | \,  {\bf{y}}-{\bf{x}}\in \{ \pm (1,0), \pm (0, 1) \} \}.
\end{equation*}
We consider a random walk on $X$ whose 
transition probability is given by
\begin{align*}
& p(({\bf{x}}, {\bf{x}}+(1,0))=\alpha, \qquad p(({\bf{x}}, {\bf{x}}-(1,0))=\alpha^\prime, 
\nonumber \\
& p(({\bf{x}}, {\bf{x}}+(0,1))=\beta, \qquad p(({\bf{x}}, {\bf{x}}-(0,1))=\beta^\prime
\end{align*}
for every ${\bf{x}}$. 
Here we assume 
$$\alpha, \alpha^\prime, \beta, \beta^\prime\geq 0 \quad \mbox{ and  } \quad
\alpha+\alpha^\prime +\beta+\beta^\prime=1.$$ 
It is easy to see that $X$ is invariant under the $\mathbb{Z}^2$-action 
generated by $\sigma_1({\bf{x}})={\bf{x}}+(1,0)$ and 
$\sigma_2({\bf{x}})={\bf{x}}+(0,1)$. 
Its quotient $X_0$ is a $2$-bouquet graph 
$X_0=(V_0, E_0)=(\{x_0\}, \{ e_1, e_2 \})$ 
(see Figure \ref{figure:square1}). 

By definition, we have
$\gamma_p=\check{\alpha} [e_1]+\check{\beta} [e_2] \in {\rm H}^{1}(X_{0}, \mathbb R)$ and $m(x_{0})=1$.
Since $X_0$ is a bouquet graph, the exterior derivative 
$d:C^0(X_0, \mathbb{R}) \rightarrow C^1(X_0,\mathbb{R})$ 
is the $0$-map. Then
\begin{equation*}
\mathrm{H}^1(X_0, \mathbb{R})\cong \mathcal{H}^1(X_0)=C^1(X_0,\mathbb{R}).
\end{equation*}
We define the canonical surjective linear map $\rho_{\mathbb R}: {\rm{H}}_{1}(X_{0}, \mathbb R)
\to \Gamma \otimes \mathbb R$ by $$\rho_{\mathbb R}([e_{i}]):=\sigma_{i} \qquad
(i=1,2).$$ 

We are going to determine the modified standard realization $\Phi_{0}: X\to 
(\Gamma \otimes \mathbb R,g_{0})$. 
Let ${\widetilde e_{i}}$ ($i=1,2$) be a lift of $e_{i}$ to $X$, and we
set $\Phi_{0}(o({\widetilde e_{1}}))=\Phi_{0}(o({\widetilde e_{2}}))={\bf 0}$. Noting that
the asymptotic direction is given by
$\rho_{\mathbb R}(\gamma_{p})=\check{\alpha} \sigma_{1} 
+\check{\beta} \sigma_{2}$, we easily see that 
\begin{equation}
\Phi_{0}(t({\widetilde e_{1}}))={\sigma_{1}}, \quad 
\Phi_{0}(t({\widetilde e_{2}}))={\sigma_{2}}
\label{7-1-0}
\nonumber
\end{equation}
is the modified harmonic realization.
We next take the dual basis $\{ \omega_1, \omega_2 \} (\subset {\rm Hom}(\Gamma, \mathbb R))$ of 
$\{ \sigma_1, \sigma_2 \}$. Namely,
$\omega_{i}[ \sigma_{j} ]_{\Gamma \otimes \mathbb R}=\delta_{ij}$ ($i,j=1,2$).
Because 
$\omega \in {\rm Hom}(\Gamma, \mathbb R)$ is identified with
$\omega=\hspace{-1mm}\mbox{ }^{t}\rho_{\mathbb R}(\omega) \in {\rm H}^{1}(X_{0}, \mathbb R)$, we have
\begin{align}
 \omega_{i}(e_{j})
=
\subscripts
{{\rm H}^{1}(X_{0},\mathbb R)}
{\big \langle 
\hspace{-1mm}\mbox{ }^{t}
\rho_{\mathbb R}(\omega_{i}), [e_{j}] \big \rangle}
{{\rm H}_{1}(X_{0}, \mathbb R)}
=
\subscripts
{{\rm Hom}(\Gamma, \mathbb R)}
{\big \langle \omega_{i}, \rho_{\mathbb R}([e_{j}]) \big \rangle}
{\Gamma \otimes \mathbb R}
=\delta_{ij}.
\label{7-1-0.5}
\end{align}
By direct computation, we have
\begin{equation}
\langle  \! \langle \omega_1, \omega_1 \rangle \! \rangle = 
\hat{\alpha}-\check{\alpha}^2, \quad
\langle \! \langle \omega_2, \omega_2 \rangle \! \rangle =
\hat{\beta} -\check{\beta}^2, \quad
\langle \! \langle \omega_1, \omega_2 \rangle \! \rangle =
-\check{\alpha}\check{\beta},
\label{7-1-1}
\end{equation}
and thus the volume of the $\Gamma$-Jacobian torus $\mathrm{Jac}^\Gamma$ is given by
\begin{equation*}
\mathrm{vol}(\mathrm{Jac}^\Gamma )=
\mathrm{vol}(\mathrm{Alb}^\Gamma )^{-1}=
\sqrt{ \det  (\langle \! \langle \omega_i, \omega_j \rangle \! \rangle)_{i,j=1}^{2} }=
\sqrt{\hat{\alpha}\hat{\beta}
-\hat{\alpha}\check{\beta}^2
-\check{\alpha}^2\hat{\beta}}.
\end{equation*}
Since $g_{0}=\big ( \langle \sigma_i, \sigma_j \rangle_{g_{0}} \big )_{i,j=1}^{2}$ is the inverse matrix of
$\big (\langle \! \langle \omega_i, \omega_j \rangle \! \rangle \big)_{i,j=1}^{2}$, 
we then obtain the Albanese metric $g_{0}$ on 
$\Gamma \otimes \mathbb R$ as
follows:
\begin{align}
\langle \sigma_1, \sigma_1 \rangle_{g_{0}} &= 
(\hat{\beta} -\check{\beta}^2)
\mathrm{vol}(\mathrm{Alb}^\Gamma )^{2},
\quad
\nonumber \\
\langle \sigma_1, \sigma_2 \rangle_{g_{0}} &= 
\check{\alpha}\check{\beta}\mathrm{vol}(\mathrm{Alb}^\Gamma )^{2},
\nonumber \\
\langle \sigma_2, \sigma_2 \rangle_{g_{0}} 
&= 
(\hat{\alpha} -\check{\alpha}^2)
\mathrm{vol}(\mathrm{Alb}^\Gamma )^{2}.
\label{7-1-2}
\nonumber
\end{align}

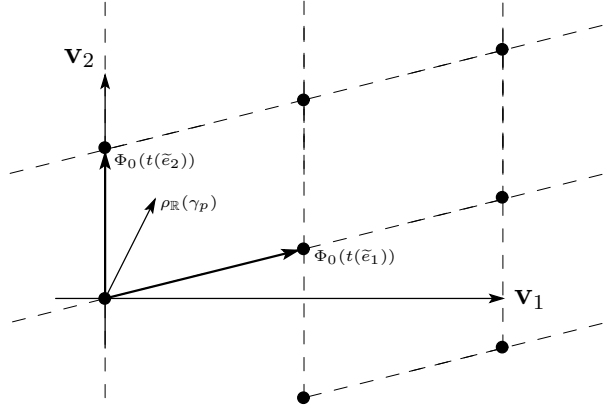
\begin{figure}[tbph]
\begin{center}
\scalebox{1.05}{
{\unitlength 0.1in%
\begin{picture}( 32.4300, 21.1000)(  0.6000,-27.1000)%
%
\special{pn 8}%
\special{pa 443 2160}%
\special{pa 2784 2160}%
\special{fp}%
\special{sh 1}%
\special{pa 2784 2160}%
\special{pa 2717 2140}%
\special{pa 2731 2160}%
\special{pa 2717 2180}%
\special{pa 2784 2160}%
\special{fp}%
%
\special{pn 4}%
\special{pa 3303 2291}%
\special{pa 2264 2550}%
\special{da 0.070}%
%
\special{pn 4}%
\special{pa 3303 1510}%
\special{pa 2264 1770}%
\special{da 0.070}%
%
\special{pn 4}%
\special{pa 3303 730}%
\special{pa 2264 990}%
\special{da 0.070}%
%
\special{pn 4}%
\special{pa 1743 1510}%
\special{pa 1743 730}%
\special{da 0.070}%
%
\special{pn 4}%
\special{pa 2784 2420}%
\special{pa 1743 2680}%
\special{da 0.070}%
%
\special{pn 4}%
\special{pa 1743 1380}%
\special{pa 1743 600}%
\special{da 0.070}%
%
\special{pn 4}%
\special{pa 2784 1380}%
\special{pa 2784 600}%
\special{da 0.070}%
%
\special{pn 4}%
\special{pa 2784 1510}%
\special{pa 2784 730}%
\special{da 0.070}%
%
\special{pn 4}%
\special{pa 2784 2420}%
\special{pa 2784 1641}%
\special{da 0.070}%
%
\special{pn 4}%
\special{pa 702 1380}%
\special{pa 702 600}%
\special{da 0.070}%
%
\special{pn 4}%
\special{pa 1223 1251}%
\special{pa 184 1510}%
\special{da 0.070}%
%
\special{pn 4}%
\special{pa 1223 2029}%
\special{pa 184 2291}%
\special{da 0.070}%
%
\special{pn 4}%
\special{pa 702 2680}%
\special{pa 702 1901}%
\special{da 0.070}%
%
\special{pn 4}%
\special{pa 2784 1641}%
\special{pa 1743 1901}%
\special{da 0.070}%
%
\special{pn 4}%
\special{pa 1743 2680}%
\special{pa 1743 1901}%
\special{da 0.070}%
%
\special{pn 4}%
\special{pa 2784 1641}%
\special{pa 2784 860}%
\special{da 0.070}%
%
\special{pn 4}%
\special{pa 2784 860}%
\special{pa 1743 1120}%
\special{da 0.070}%
%
\special{pn 4}%
\special{pa 1743 1120}%
\special{pa 702 1380}%
\special{da 0.070}%
%
\special{pn 4}%
\special{pa 1743 1901}%
\special{pa 1743 1120}%
\special{da 0.070}%
%
\special{pn 8}%
\special{pa 702 2420}%
\special{pa 702 990}%
\special{fp}%
\special{sh 1}%
\special{pa 702 990}%
\special{pa 682 1057}%
\special{pa 702 1043}%
\special{pa 722 1057}%
\special{pa 702 990}%
\special{fp}%
%
\special{pn 13}%
\special{pa 702 2160}%
\special{pa 1705 1910}%
\special{fp}%
\special{sh 1}%
\special{pa 1705 1910}%
\special{pa 1635 1907}%
\special{pa 1653 1923}%
\special{pa 1645 1946}%
\special{pa 1705 1910}%
\special{fp}%
\put(17.9500,-19.0000){\makebox(0,0)[lt]{{\tiny $\Phi_0(t(\widetilde{e}_1))$}}}%
%
\special{pn 13}%
\special{pa 702 2160}%
\special{pa 702 1405}%
\special{fp}%
\special{sh 1}%
\special{pa 702 1405}%
\special{pa 682 1472}%
\special{pa 702 1458}%
\special{pa 722 1472}%
\special{pa 702 1405}%
\special{fp}%
\put(7.4800,-13.9700){\makebox(0,0)[lt]{{\tiny $\Phi_0(t(\widetilde{e}_2))$}}}%
%
\special{pn 8}%
\special{pa 702 2160}%
\special{pa 963 1641}%
\special{fp}%
\special{sh 1}%
\special{pa 963 1641}%
\special{pa 915 1692}%
\special{pa 939 1689}%
\special{pa 951 1710}%
\special{pa 963 1641}%
\special{fp}%
\put(9.8900,-17.0500){\makebox(0,0)[lb]{{\tiny $\rho_\mathbb{R}(\gamma_p)$}}}%
\put(29.2500,-21.6500){\makebox(0,0){${\bf v}_1$}}%
\put(6.5000,-9.5000){\makebox(0,0)[rb]{${\bf v}_2$}}%
%
\special{sh 1.000}%
\special{ia 700 2160 30 30  0.0000000  6.2831853}%
\special{pn 8}%
\special{ar 700 2160 30 30  0.0000000  6.2831853}%
%
\special{sh 1.000}%
\special{ia 1740 1900 30 30  0.0000000  6.2831853}%
\special{pn 8}%
\special{ar 1740 1900 30 30  0.0000000  6.2831853}%
%
\special{sh 1.000}%
\special{ia 1740 2680 30 30  0.0000000  6.2831853}%
\special{pn 8}%
\special{ar 1740 2680 30 30  0.0000000  6.2831853}%
%
\special{sh 1.000}%
\special{ia 1740 1120 30 30  0.0000000  6.2831853}%
\special{pn 8}%
\special{ar 1740 1120 30 30  0.0000000  6.2831853}%
%
\special{sh 1.000}%
\special{ia 700 1370 30 30  0.0000000  6.2831853}%
\special{pn 8}%
\special{ar 700 1370 30 30  0.0000000  6.2831853}%
%
\special{sh 1.000}%
\special{ia 2780 2415 30 30  0.0000000  6.2831853}%
\special{pn 8}%
\special{ar 2780 2415 30 30  0.0000000  6.2831853}%
%
\special{sh 1.000}%
\special{ia 2780 1630 30 30  0.0000000  6.2831853}%
\special{pn 8}%
\special{ar 2780 1630 30 30  0.0000000  6.2831853}%
%
\special{sh 1.000}%
\special{ia 2780 855 30 30  0.0000000  6.2831853}%
\special{pn 8}%
\special{ar 2780 855 30 30  0.0000000  6.2831853}%
\end{picture}}%
}
\end{center}
\caption{Modified standard realization of the square lattice graph}
\label{figure:square-mharm}
\end{figure}

Let $\{ v_1, v_2 \}$ be the orthonormal basis of $\mathrm{Hom}(\Gamma, \mathbb{R})
(\subset \mathcal{H}^1(X_0))$ given by the Gram--Schmidt orthonormalization of $\{\omega_1, \omega_2\}$.
It follows from (\ref{7-1-1}) that
$$ v_{1}=\frac{1}{\sqrt{ \hat{\alpha}-\check{\alpha}^2}}\omega_{1}, \quad v_{2}=
{\mathrm{vol}(\mathrm{Alb}^\Gamma )} \Big(
\frac{\check{\alpha}\check{\beta}}
{\sqrt{\hat{\alpha}-\check{\alpha}^2}}\omega_{1}+
{\sqrt{\hat{\alpha}-\check{\alpha}^2}} \omega_{2} \Big),
$$
and thus we obtain
\begin{align*}
v_1(e_1)=&\frac{1}{\sqrt{ \hat{\alpha}-\check{\alpha}^2}},& 
v_1(e_2)&=0,\\
v_2(e_1)=&\frac{\check{\alpha}\check{\beta}
\hspace{0.8mm}
\mathrm{vol}(\mathrm{Alb}^\Gamma )}
{\sqrt{\hat{\alpha}-\check{\alpha}^2}}, &
v_2(e_2)&={\sqrt{\hat{\alpha}-\check{\alpha}^2}}
\hspace{0.8mm}
\mathrm{vol}(\mathrm{Alb}^\Gamma ).
\end{align*}
Let $\{ {\bf v}_{1}, {\bf v}_{2} \}$ denote the dual basis of 
$\{ v_1, v_2 \}$. Needless to say,  $\{ {\bf v}_{1}, {\bf v}_{2} \}$ is an orthonormal basis of
$\Gamma \otimes \mathbb R$. 
As (\ref{7-1-0.5}), we have
\begin{equation}
\sigma_{i}=
\sum_{j=1}^{2}
\subscripts
{{\rm Hom}(\Gamma, \mathbb R)}
{\big \langle v_{j}, \sigma_{i}  \big \rangle}
{\Gamma \otimes \mathbb R}
{\bf v}_{j}=\sum_{j=1}^{2} v_{j}(e_{i}) {\bf v}_{j} \qquad (i,j=1,2).
\label{7-1-eq}
\nonumber
\end{equation}
It means that $\sigma_{i} \in \Gamma \otimes \mathbb R$ may be identified with 
$(v_{1}(e_{i}), v_{2}(e_{i})) \in \mathbb R^{2}$. In this sense, 
the standard realization $\Phi_0:X\rightarrow (\Gamma\otimes \mathbb{R}, g_{0}) \cong 
(\mathbb R^{2}, \{ {\bf v}_{1}, {\bf v}_{2} \})$ is also
given by $\Phi_{0}(o( \widetilde{e}_1 ))=\Phi_{0}(o( \widetilde{e}_2 ))=(0,0)$ and
\begin{align*}
\Phi_0 (t( \widetilde{e}_1 )) =
\frac{1}{\sqrt{ \hat{\alpha}-\check{\alpha}^2}}
\Big (1 ,  
\check{\alpha}\check{\beta} \hspace{0.8mm}
\mathrm{vol}(\mathrm{Alb}^\Gamma ) \Big), 
\quad
\Phi_0 (t( \widetilde{e}_2 )) =& \Big( 0, 
{\sqrt{\hat{\alpha}-\check{\alpha}^2}}
\hspace{0.8mm}
\mathrm{vol}(\mathrm{Alb}^\Gamma )
\Big)
\end{align*}
(see Figure \ref{figure:square-mharm}). 

In particular, when the random walk is simple (i.e., $\alpha=\alpha'=\beta=\beta'=1/4$), 
$$\mathrm{vol}(\mathrm{Alb}^\Gamma )=2, \quad
\Phi_0 (t( \widetilde{e}_1 )) =
\big (\sqrt{2} , 0 \big), \quad  
\Phi_0 (t( \widetilde{e}_2 )) =\big( 0, \sqrt{2} \big).$$
As mentioned before, this realization is different from the 
standard realization in \cite[page 685]{KSS} due to the difference of both
the flat metric (\ref{flatmetric-norm}) and the invariant measure $m$.
\subsection{The triangular lattice}
\begin{figure}[tbph]
\begin{center}
\scalebox{1.0}{\input{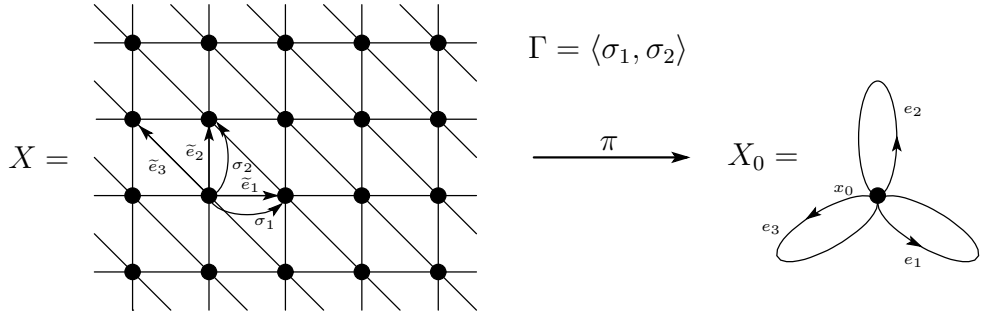}}
\end{center}
\caption{Triangular lattice and its quotient}
\label{figure:triangular}
\end{figure}

We consider a class of non-symmetric random walks on the triangular lattice discussed in
Ishiwata--Kawabi--Teruya \cite{IKT}.
Let $X=(V,E)$ be a triangular lattice, where $V=\mathbb{Z}^2$ and 
\begin{equation*}
E=\big \{ (\mathbf{x}, \mathbf{y} ) \in V \times V \, | \, 
\mathbf{y}-\mathbf{x}\in \{ \pm (1, 0), \pm (0, 1), \pm (-1, 1) \} 
\big \}.
\end{equation*}
The transition probability of the random walk is given by
\begin{align*}
& p((\mathbf{x},\mathbf{x}+(1,0)))=\alpha,  \qquad  p((\mathbf{x},\mathbf{x}-(1,0)))=\alpha^\prime, \\
& p((\mathbf{x},\mathbf{x}+(0,1)))=\beta^\prime, \qquad  p((\mathbf{x},\mathbf{x}-(0,1)))=\beta, \\
& p((\mathbf{x},\mathbf{x}+(-1,1)))=\gamma,  \qquad  p((\mathbf{x},\mathbf{x}-(-1,1)))=\gamma^\prime.
\end{align*}
Here we assume 
$$\alpha, \alpha^\prime, \beta, \beta^\prime, \gamma, \gamma^\prime>0, \quad 
{\hat \alpha}+{\hat \beta}+{\hat \gamma}=1 \quad \mbox{and} 
\quad {\check \alpha}={\check \beta}={\check \gamma}=\kappa \geq 0.$$
The constant $\kappa$ can be regarded the intensity of non-symmetry of the random walk, and 
under this assumption, the random walk on $X$ is aperiodic (i.e., $K=1$). 
We also mention that the random walk is of period $K=3$ when $\alpha'=\beta'=\gamma'=0$.
It is easy to see that $X$ is invariant under the $\mathbb Z^{2}$-action 
generated by
$\sigma_1(\mathbf{x})=\mathbf{x}+(1,0)$ and $\sigma_2(\mathbf{x})=\mathbf{x}+(0,1)$.
Its quotint graph $X_0=(V_0,E_0)$ is a $3$-bouquet graph consisting of $V_0=\{ x_0 \}$, 
$E_0=\{e_1, e_2, e_3\}$ (see Figure \ref{figure:triangular}). 

The first homology group $\mathrm{H}_1(X_0, \mathbb{R})$ is spanned by $[e_1], [e_2], [e_3]$. 
By definition, 
we have
$$\gamma_p=\check{\alpha} [e_1] -\check{\beta} [e_2] +\check{\gamma} [e_3]=\kappa
([e_1] - [e_2] +[e_3]).$$ 
Since the quotient graph $X_0$ is a bouquet graph, we see $\mathrm{H}^1(X_0,\mathbb{R})\cong
\mathcal{H}^1(X_0) =C^1(X_0,\mathbb{R})$ and $m(x_{0})=1$.
We define the canonical surjective linear map 
$\rho_{\mathbb{R}} : \mathrm{H}_1(X_0, \mathbb{R}) \rightarrow \Gamma\otimes \mathbb{R}$ by
\begin{equation*}
\rho_{\mathbb{R}} ([e_1]):= \sigma_1, \quad \, \rho_{\mathbb{R}}([e_2])
:=\sigma_2, \quad \rho_{\mathbb{R}}([ e_3 ]):= \sigma_2 - \sigma_1.
\end{equation*}
Then the asymptotic direction is given by
\begin{equation*}
\rho_{\mathbb{R}}( \gamma_p)=\left( \check{\alpha}-\check{\gamma}\right) \sigma_1 
+\left( \check{\gamma} -\check{\beta}\right) \sigma_2 ={\bf 0},
\end{equation*}
and we find that $\rho_{\mathbb R}(\gamma_{p})={\bf 0}$ is equivalent to
${\check \alpha}={\check \beta}={\check \gamma}$.

We introduce the basis $\{u_{1}, u_{2} \}$ in ${\rm Hom}(\Gamma, \mathbb R)$ by 
$$u_{1} \big( (m,n) \big)=m, \quad u_{2} \big( (m,n) \big )=n  \quad \mbox{for }(m,n)
\in \mathbb Z^{2}\cong \Gamma.$$
Note that $\{ u_{1}, u_{2} \}$ is the dual basis of $\{ \sigma_{1}, \sigma_{2} \}$. 
Let $\{ \omega_1, \omega_2, \omega_3 \}(\subset {\rm H}^{1}(X_{0}, \mathbb R))$ be the dual basis of 
$\{ [e_1], [e_2], [e_3] \}(\subset {\rm H}_{1}(X_{0}, \mathbb R))$. 
By direct computation, we have
\begin{equation}
\begin{cases} 
~\langle \! \langle \omega_1, \omega_1 \rangle \! \rangle = 
\hat{\alpha} -\check{\alpha}^2=\hat{\alpha}-\kappa^{2},
\qquad
\langle \! \langle \omega_1, \omega_2 \rangle \! \rangle =\check{\alpha}\check{\beta}=\kappa^{2}, 
\\
~\langle \! \langle \omega_2, \omega_2 \rangle \! \rangle =\hat{\beta}-\check{\beta}^2=\hat{\beta}-\kappa^{2},  
\qquad
\langle \! \langle \omega_2, \omega_3 \rangle \! \rangle =\check{\beta}\check{\gamma}=\kappa^{2}, 
\\
~\langle \! \langle \omega_3, \omega_3 \rangle \! \rangle =\hat{\gamma}-\check{\gamma}^2
={\hat \gamma}-\kappa^{2} ,
\qquad
\langle \! \langle \omega_1, \omega_3 \rangle \! \rangle = -\check{\alpha}\check{\gamma}=-\kappa^{2}.
\end{cases}
\label{7-2-1}
\end{equation}
Since $X$ is a non-maximal abelian covering graph of $X_{0}$ with the covering transformation group $\Gamma \cong
\mathbb Z^{2}$,
we need to find a $\mathbb Z$-basis of the lattice
\begin{equation*}
L=\left\{ \omega \in {\rm{H}}^1(X_0, \mathbb Z) \, | \, \omega ( [c])=0 \mbox{ for every cycle $\widetilde {c}$ 
on $X$} \right\}.
\end{equation*}
It is easy to find that $u_1=\hspace{-1mm}\mbox{ }^{t}\rho_{\mathbb R}(u_{1})
=\omega_1-\omega_3$ and 
$u_2=\hspace{-1mm}\mbox{ }^{t}\rho_{\mathbb R}(u_{2})=
\omega_2+\omega_3$ form a $\mathbb{Z}$-basis of 
the lattice $L$.
It follows from (\ref{7-2-1}) that
\begin{align}
\langle \! \langle u_1, u_1 \rangle \! \rangle &={\hat \alpha}+{\hat \gamma}, \quad
\langle \! \langle u_1, u_2 \rangle \! \rangle =-{\hat \gamma}, \quad
\langle \! \langle u_2, u_2 \rangle \! \rangle =
\hat{\beta}
+\hat{\gamma},
\label{7-2-3}
\end{align}
and thus we obtain
\begin{equation*}
 \mathrm{vol}(\mathrm{Jac}^\Gamma)=\mathrm{vol}(\mathrm{Alb}^\Gamma)^{-1}=
 \sqrt{ \det \left( \langle \! \langle u_i, u_j \rangle \! \rangle \right)_{i,j=1}^{2}}= 
\sqrt{ \hat{\alpha}\hat{\beta}+\hat{\beta}\hat{\gamma}+\hat{\gamma}\hat{\alpha}},
\end{equation*}
and the Albanese metric $g_{0}$ on $\Gamma \otimes \mathbb R$ as follows:
\begin{align}
\langle \sigma_1, \sigma_1 \rangle_{g_{0}} &= 
\frac{ \hat{\beta} +\hat{\gamma} }{\hat{\alpha}\hat{\beta}+\hat{\beta}\hat{\gamma}+\hat{\gamma}\hat{\alpha}},
\quad
\nonumber \\
\langle \sigma_1, \sigma_2 \rangle_{g_{0}} 
&= 
\frac{ \hat{\gamma}}{\hat{\alpha}\hat{\beta}+\hat{\beta}\hat{\gamma}+\hat{\gamma}\hat{\alpha}},
\quad
\nonumber \\
\langle \sigma_2, \sigma_2 \rangle_{g_{0}} &= 
\frac{ \hat{\alpha} +\hat{\gamma}}{\hat{\alpha}\hat{\beta}+\hat{\beta}\hat{\gamma}+\hat{\gamma}\hat{\alpha}}.
\label{7-2-2}
\nonumber
\end{align}
 
\begin{figure}[htbp]
\begin{center}
\scalebox{1.0}{
{\unitlength 0.1in%
\begin{picture}( 30.0500, 22.1400)( 11.8000,-23.3700)%
%
\special{pn 4}%
\special{pa 2673 2167}%
\special{pa 3098 2337}%
\special{da 0.070}%
%
\special{pn 4}%
\special{pa 3355 1485}%
\special{pa 4038 1757}%
\special{da 0.070}%
%
\special{pn 4}%
\special{pa 3355 528}%
\special{pa 4038 802}%
\special{da 0.070}%
%
\special{pn 4}%
\special{pa 2673 255}%
\special{pa 3355 528}%
\special{da 0.070}%
%
\special{pn 4}%
\special{pa 1988 1894}%
\special{pa 2671 2167}%
\special{da 0.070}%
%
\special{pn 4}%
\special{pa 2673 1210}%
\special{pa 3355 1485}%
\special{da 0.070}%
%
\special{pn 4}%
\special{pa 3355 1485}%
\special{pa 2673 2167}%
\special{da 0.070}%
%
\special{pn 4}%
\special{pa 3355 528}%
\special{pa 2673 1210}%
\special{da 0.070}%
%
\special{pn 4}%
\special{pa 1307 664}%
\special{pa 1989 939}%
\special{da 0.070}%
%
\special{pn 4}%
\special{pa 2262 1621}%
\special{pa 1580 2303}%
\special{da 0.070}%
%
\special{pn 4}%
\special{pa 1989 939}%
\special{pa 1307 1621}%
\special{da 0.070}%
%
\special{pn 4}%
\special{pa 1989 939}%
\special{pa 2673 1210}%
\special{da 0.070}%
%
\special{pn 4}%
\special{pa 1307 1621}%
\special{pa 1989 1894}%
\special{da 0.070}%
%
\special{pn 13}%
\special{pa 1989 1894}%
\special{pa 2642 1241}%
\special{fp}%
\special{sh 1}%
\special{pa 2642 1241}%
\special{pa 2581 1274}%
\special{pa 2604 1279}%
\special{pa 2609 1302}%
\special{pa 2642 1241}%
\special{fp}%
\put(29.8800,-11.8200){\makebox(0,0){{\tiny $\Phi_0(t(\widetilde{e}_1))$}}}%
\put(22.8800,-9.0000){\makebox(0,0){{\tiny $\Phi_0(t(\widetilde{e}_2))$}}}%
%
\special{pn 4}%
\special{pa 1779 188}%
\special{pa 1306 661}%
\special{da 0.070}%
%
\special{pn 4}%
\special{pa 4042 798}%
\special{pa 3358 1482}%
\special{da 0.070}%
%
\special{pn 4}%
\special{pa 2671 2163}%
\special{pa 2497 2337}%
\special{da 0.070}%
%
\special{pn 4}%
\special{pa 4042 1762}%
\special{pa 3479 2325}%
\special{da 0.070}%
%
\special{pn 8}%
\special{pa 1200 1893}%
\special{pa 4074 1893}%
\special{fp}%
\special{sh 1}%
\special{pa 4074 1893}%
\special{pa 4007 1873}%
\special{pa 4021 1893}%
\special{pa 4007 1913}%
\special{pa 4074 1893}%
\special{fp}%
%
\special{pn 8}%
\special{pa 1982 2331}%
\special{pa 1982 284}%
\special{fp}%
\special{sh 1}%
\special{pa 1982 284}%
\special{pa 1962 351}%
\special{pa 1982 337}%
\special{pa 2002 351}%
\special{pa 1982 284}%
\special{fp}%
%
\special{pn 13}%
\special{pa 1982 1888}%
\special{pa 1982 972}%
\special{fp}%
\special{sh 1}%
\special{pa 1982 972}%
\special{pa 1962 1039}%
\special{pa 1982 1025}%
\special{pa 2002 1039}%
\special{pa 1982 972}%
\special{fp}%
%
\special{pn 4}%
\special{pa 2359 128}%
\special{pa 2671 252}%
\special{da 0.070}%
%
\special{pn 4}%
\special{pa 3717 166}%
\special{pa 3035 849}%
\special{da 0.070}%
%
\special{pn 4}%
\special{pa 2765 166}%
\special{pa 2083 849}%
\special{da 0.070}%
\put(41.5500,-18.8800){\makebox(0,0){{\tiny ${\bf v}_1$}}}%
\put(19.7600,-1.8800){\makebox(0,0){{\tiny ${\bf v}_2$}}}%
%
\special{pn 4}%
\special{pa 4042 2253}%
\special{pa 4042 194}%
\special{da 0.070}%
%
\special{pn 4}%
\special{pa 1306 2253}%
\special{pa 1306 194}%
\special{da 0.070}%
%
\special{pn 4}%
\special{pa 2671 2253}%
\special{pa 2671 194}%
\special{da 0.070}%
%
\special{pn 4}%
\special{pa 3353 2253}%
\special{pa 3353 194}%
\special{da 0.070}%
%
\special{sh 1.000}%
\special{ia 1982 1894 36 36  0.0000000  6.2831853}%
\special{pn 8}%
\special{ar 1982 1894 36 36  0.0000000  6.2831853}%
%
\special{sh 1.000}%
\special{ia 2671 2163 36 36  0.0000000  6.2831853}%
\special{pn 8}%
\special{ar 2671 2163 36 36  0.0000000  6.2831853}%
%
\special{sh 1.000}%
\special{ia 1306 1619 36 36  0.0000000  6.2831853}%
\special{pn 8}%
\special{ar 1306 1619 36 36  0.0000000  6.2831853}%
%
\special{sh 1.000}%
\special{ia 4042 1762 35 35  0.0000000  6.2831853}%
\special{pn 8}%
\special{ar 4042 1762 35 35  0.0000000  6.2831853}%
%
\special{sh 1.000}%
\special{ia 3353 1487 36 36  0.0000000  6.2831853}%
\special{pn 8}%
\special{ar 3353 1487 36 36  0.0000000  6.2831853}%
%
\special{sh 1.000}%
\special{ia 2671 1200 36 36  0.0000000  6.2831853}%
\special{pn 8}%
\special{ar 2671 1200 36 36  0.0000000  6.2831853}%
%
\special{sh 1.000}%
\special{ia 1982 936 36 36  0.0000000  6.2831853}%
\special{pn 8}%
\special{ar 1982 936 36 36  0.0000000  6.2831853}%
%
\special{sh 1.000}%
\special{ia 1306 661 36 36  0.0000000  6.2831853}%
\special{pn 8}%
\special{ar 1306 661 36 36  0.0000000  6.2831853}%
%
\special{sh 1.000}%
\special{ia 2671 254 36 36  0.0000000  6.2831853}%
\special{pn 8}%
\special{ar 2671 254 36 36  0.0000000  6.2831853}%
%
\special{sh 1.000}%
\special{ia 3353 523 36 36  0.0000000  6.2831853}%
\special{pn 8}%
\special{ar 3353 523 36 36  0.0000000  6.2831853}%
%
\special{sh 1.000}%
\special{ia 4042 804 35 35  0.0000000  6.2831853}%
\special{pn 8}%
\special{ar 4042 804 35 35  0.0000000  6.2831853}%
%
\special{pn 4}%
\special{pa 1192 1571}%
\special{pa 1300 1614}%
\special{da 0.070}%
%
\special{pn 4}%
\special{pa 1307 1621}%
\special{pa 1180 1748}%
\special{da 0.070}%
%
\special{pn 4}%
\special{pa 1307 663}%
\special{pa 1180 790}%
\special{da 0.070}%
%
\special{pn 4}%
\special{pa 1192 613}%
\special{pa 1300 656}%
\special{da 0.070}%
%
\special{pn 4}%
\special{pa 4044 1760}%
\special{pa 4185 1816}%
\special{da 0.070}%
%
\special{pn 4}%
\special{pa 4143 1661}%
\special{pa 4042 1762}%
\special{da 0.070}%
%
\special{pn 4}%
\special{pa 4143 703}%
\special{pa 4042 804}%
\special{da 0.070}%
%
\special{pn 4}%
\special{pa 4044 802}%
\special{pa 4185 858}%
\special{da 0.070}%
\end{picture}}%
}
\end{center}
\caption{Modified standard realization of the triangular lattice in the case 
$\rho_\mathbb{R} (\gamma_p)={\bf 0}$.}
\label{figure:triangular-mharm}
\end{figure}
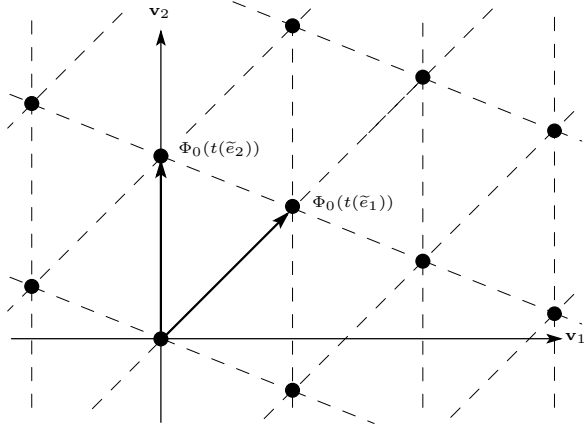

We now determine the modified standard realization $\Phi_{0}: X \to (\Gamma\otimes \mathbb{R},g_{0})$.
Let ${\widetilde e_{i}}$ ($i=1,2$) be a lift of $e_{i}$ to $X$, and we
set $\Phi_{0}(o({\widetilde e_{1}}))=\Phi_{0}(o({\widetilde e_{2}}))={\bf 0}$. As in Example 7.1,
we also see that 
\begin{equation}
\Phi_{0}(t({\widetilde e_{1}}))={\sigma_{1}}, \quad 
\Phi_{0}(t({\widetilde e_{2}}))={\sigma_{2}}
\label{7-2-4}
\nonumber
\end{equation}
is the modified harmonic realization.
Let $\{ v_1, v_2\}$ be the Gram--Schmidt orthogonalization of $\{u_1,u_2\}$,
and let $\{ {\bf v}_{1}, {\bf v}_{2} \}$
be the dual basis of
$\{ v_1, v_2\}$ in $\Gamma \otimes \mathbb R$.
By (\ref{7-2-3}), we have
$$ v_{1}=\frac{1}{\sqrt{ \hat{\alpha}+\hat{\gamma}}}u_{1}, \quad v_{2}=
{\sqrt{ \frac{\hat \alpha+\hat \gamma}{
\hat{\alpha}\hat{\beta}+\hat{\beta}\hat{\gamma}+\hat{\gamma}\hat{\alpha}}}}
\Big(
\frac{\hat \gamma}{\hat \alpha+\hat \gamma}u_{1}+u_{2}\Big),
$$
and this implies
\begin{align*}
&v_1(e_1)= \frac{1}{\sqrt{ \hat{\alpha}+\hat{\gamma} }},  &
v_1(e_2)=&0 ,\\
&v_2(e_1)=\frac{\hat \gamma \hspace{0.8mm} {\rm vol}({\rm Alb}^{\Gamma})}{\sqrt{
\hat{\alpha}+\hat{\gamma}
}}, &
v_2(e_2)=&
{\sqrt{ \hat{\alpha}+\hat{\gamma} }}
\hspace{0.8mm}
{\rm vol}({\rm Alb}^{\Gamma}).
\end{align*}
Therefore the modified standard realization 
$\Phi_{0}: X \rightarrow (\Gamma\otimes \mathbb{R},g_{0}) \cong (\mathbb R^{2}, \{ {\bf v}_{1}, {\bf v}_{2} \})$ 
is given by $\Phi_{0}(o( \widetilde{e}_1 ))=\Phi_{0}(o( \widetilde{e}_2 ))=(0,0)$ and
\begin{align*}
\Phi_0 (t( \widetilde{e}_1 )) =
\frac{1}{\sqrt{ \hat{\alpha}+\hat{\gamma}}}
\Big(1,  \hat \gamma \hspace{0.8mm} {\rm vol}({\rm Alb}^{\Gamma}) \Big), 
\quad
\Phi_0 (t( \widetilde{e}_2 )) =& \Big( 0, 
{\sqrt{ \hat{\alpha}+\hat{\gamma} }}
\hspace{0.8mm}
{\rm vol}({\rm Alb}^{\Gamma}) \Big)
\end{align*}
(see Figure \ref{figure:triangular-mharm}).
This realization is same as in \cite[page 133]{IKT}, 
and if we consider the simple random walk (i.e., $\alpha=\alpha'=\beta=\beta'=\gamma=\gamma'=1/6$),
$$\mathrm{vol}(\mathrm{Alb}^\Gamma )={\sqrt 3}, \quad
\Phi_0 (t( \widetilde{e}_1 )) =
\big (\frac{\sqrt 6}{2} , \frac{\sqrt 2}{2} \big), \quad  
\Phi_0 (t( \widetilde{e}_2 )) =\big( 0, \sqrt{2} \big).$$

Furthermore, Theorems \ref{AE-LCLT} and \ref{a_1} (see also \cite[Theorem 2.2]{IKT}) also yield
the following precise asymptotics of the 
$n$-step transition probability $p(n,x,y)$:
\begin{align*}
2\pi n \cdot p(n,x,y) &={\rm{vol}}({\rm{Alb}}^{\Gamma}) \exp 
\Big (-\frac{\vert \Phi_{0}(y)-\Phi_{0}(x) \vert_{g_{0}}^{2}}{2n} \Big)
\nonumber \\
&\quad \times \big (1+a_{1}(\kappa; \Phi_{0}(y)-\Phi_{0}(x))n^{-1} \big)+O(n^{-3/2})
\end{align*}
as $n \to \infty$ uniformly for $x,y\in V$ with $\vert \Phi_{0}(y)-\Phi_{0}(x) \vert_{g_{0}}\leq Cn^{1/6}$. Here
the coefficient $a_{1}(\kappa; \Phi_{0}(y)-\Phi_{0}(x))$ is explicitly given by
\begin{align*}
a_{1}(\kappa; \Phi_{0}(y)-\Phi_{0}(x))&=
-1+\frac{{\rm vol}({\rm Alb}^{\Gamma})^{4}}{8} \big \{ {\hat \alpha}({\hat \beta}+{\hat \gamma})^{2}
+{\hat \beta}({\hat \gamma}+{\hat \alpha})^{2}+{\hat \gamma}({\hat \alpha}+{\hat \beta})^{2} \big \}
\nonumber \\
&\quad 
+ {\rm vol}({\rm Alb}^{\Gamma})^{4} \Big \{ ({\hat \alpha}{\hat \beta}-2{\hat \beta}{\hat \gamma}
+{\hat \gamma}{\hat \alpha})\big(\Phi_{0}(y)-\Phi_{0}(x) \big)_{1}
\nonumber \\
&\qquad \qquad \qquad \qquad +
(-{\hat \alpha}{\hat \beta}-{\hat \beta}{\hat \gamma}
+2{\hat \gamma}{\hat \alpha})\big(\Phi_{0}(y)-\Phi_{0}(x) \big)_{2} \Big \} \kappa
\nonumber \\
&\quad
+\frac{3}{8}{\rm vol}({\rm Alb}^{\Gamma})^{4} \Big (-1+5
{\hat \alpha}{\hat \beta}{\hat \gamma}\hspace{0.8mm}{\rm vol}({\rm Alb}^{\Gamma})^{2} \Big) \kappa^{2}.
\end{align*}
\subsection{The hexagonal lattice}
\begin{figure}[htbp]
\begin{center}
\scalebox{1.0}{
{\unitlength 0.1in%
\begin{picture}( 38.5900, 12.0900)(  3.9000,-14.0900)%
%
\special{pn 8}%
\special{pa 1488 804}%
\special{pa 1929 363}%
\special{da 0.070}%
\special{sh 1}%
\special{pa 1929 363}%
\special{pa 1868 396}%
\special{pa 1891 401}%
\special{pa 1896 424}%
\special{pa 1929 363}%
\special{fp}%
\special{pa 1488 804}%
\special{pa 1035 345}%
\special{da 0.070}%
\special{sh 1}%
\special{pa 1035 345}%
\special{pa 1068 406}%
\special{pa 1072 383}%
\special{pa 1096 378}%
\special{pa 1035 345}%
\special{fp}%
%
\special{pn 13}%
\special{pa 521 321}%
\special{pa 521 200}%
\special{fp}%
\special{pa 1488 200}%
\special{pa 1488 321}%
\special{fp}%
\special{pa 2455 321}%
\special{pa 2455 200}%
\special{fp}%
\special{pa 1971 1288}%
\special{pa 1971 1409}%
\special{fp}%
\special{pa 1004 1409}%
\special{pa 1004 1288}%
\special{fp}%
%
\special{pn 13}%
\special{pa 521 804}%
\special{pa 521 1288}%
\special{fp}%
%
\special{pn 13}%
\special{pa 1004 321}%
\special{pa 1004 804}%
\special{fp}%
%
\special{pn 13}%
\special{pa 2455 804}%
\special{pa 2455 1288}%
\special{fp}%
%
\special{pn 13}%
\special{pa 1971 321}%
\special{pa 1971 804}%
\special{fp}%
%
\special{pn 13}%
\special{pa 1488 804}%
\special{pa 1488 1288}%
\special{fp}%
%
\special{pn 13}%
\special{pa 400 321}%
\special{pa 2576 321}%
\special{fp}%
\special{pa 2576 804}%
\special{pa 400 804}%
\special{fp}%
\special{pa 400 1288}%
\special{pa 2576 1288}%
\special{fp}%
%
\special{sh 0}%
\special{ia 1488 804 42 42  0.0000000  6.2831853}%
\special{pn 8}%
\special{ar 1488 804 42 42  0.0000000  6.2831853}%
%
\special{sh 0}%
\special{ia 2455 804 42 42  0.0000000  6.2831853}%
\special{pn 8}%
\special{ar 2455 804 42 42  0.0000000  6.2831853}%
%
\special{sh 0}%
\special{ia 521 804 42 42  0.0000000  6.2831853}%
\special{pn 8}%
\special{ar 521 804 42 42  0.0000000  6.2831853}%
%
\special{sh 0}%
\special{ia 1971 1288 43 43  0.0000000  6.2831853}%
\special{pn 8}%
\special{ar 1971 1288 43 43  0.0000000  6.2831853}%
%
\special{sh 0}%
\special{ia 1004 1288 43 43  0.0000000  6.2831853}%
\special{pn 8}%
\special{ar 1004 1288 43 43  0.0000000  6.2831853}%
%
\special{sh 0}%
\special{ia 1004 321 43 42  0.0000000  6.2831853}%
\special{pn 8}%
\special{ar 1004 321 43 42  0.0000000  6.2831853}%
%
\special{sh 0}%
\special{ia 1971 321 43 42  0.0000000  6.2831853}%
\special{pn 8}%
\special{ar 1971 321 43 42  0.0000000  6.2831853}%
%
\special{sh 1.000}%
\special{ia 1488 321 42 42  0.0000000  6.2831853}%
\special{pn 8}%
\special{ar 1488 321 42 42  0.0000000  6.2831853}%
%
\special{sh 1.000}%
\special{ia 1971 804 43 42  0.0000000  6.2831853}%
\special{pn 8}%
\special{ar 1971 804 43 42  0.0000000  6.2831853}%
%
\special{sh 1.000}%
\special{ia 2455 1288 42 43  0.0000000  6.2831853}%
\special{pn 8}%
\special{ar 2455 1288 42 43  0.0000000  6.2831853}%
%
\special{sh 1.000}%
\special{ia 2455 321 42 42  0.0000000  6.2831853}%
\special{pn 8}%
\special{ar 2455 321 42 42  0.0000000  6.2831853}%
%
\special{sh 1.000}%
\special{ia 521 321 42 42  0.0000000  6.2831853}%
\special{pn 8}%
\special{ar 521 321 42 42  0.0000000  6.2831853}%
%
\special{sh 1.000}%
\special{ia 521 1288 42 43  0.0000000  6.2831853}%
\special{pn 8}%
\special{ar 521 1288 42 43  0.0000000  6.2831853}%
%
\special{sh 1.000}%
\special{ia 1488 1288 42 43  0.0000000  6.2831853}%
\special{pn 8}%
\special{ar 1488 1288 42 43  0.0000000  6.2831853}%
%
\special{sh 1.000}%
\special{ia 1004 804 43 42  0.0000000  6.2831853}%
\special{pn 8}%
\special{ar 1004 804 43 42  0.0000000  6.2831853}%
\put(14.8200,-7.1400){\makebox(0,0){$\widetilde{\bf{x}}_1$}}%
\put(14.9400,-14.2100){\makebox(0,0){$\widetilde{\bf{x}}_2$}}%
\put(17.3600,-5.0800){\makebox(0,0)[rb]{$\sigma_1$}}%
\put(12.2800,-5.0800){\makebox(0,0)[lb]{$\sigma_2$}}%
%
\special{pn 13}%
\special{pa 2696 804}%
\special{pa 3542 804}%
\special{fp}%
\special{sh 1}%
\special{pa 3542 804}%
\special{pa 3475 784}%
\special{pa 3489 804}%
\special{pa 3475 824}%
\special{pa 3542 804}%
\special{fp}%
\put(31.0100,-7.5600){\makebox(0,0){$\pi$}}%
\put(30.8900,-4.4800){\makebox(0,0){$\Gamma=\langle \sigma_1, \sigma_2 \rangle$}}%
%
\special{pn 8}%
\special{ar 4026 804 181 363  0.0000000  6.2831853}%
%
\special{pn 8}%
\special{pa 4026 1167}%
\special{pa 4026 442}%
\special{fp}%
\put(40.9200,-12.3300){\makebox(0,0)[lb]{${\bf{x}}_1$}}%
\put(40.8600,-3.8700){\makebox(0,0)[lb]{${\bf{x}}_2$}}%
%
\special{sh 0}%
\special{ia 4026 1167 42 42  0.0000000  6.2831853}%
\special{pn 8}%
\special{ar 4026 1167 42 42  0.0000000  6.2831853}%
%
\special{sh 1.000}%
\special{ia 4026 442 42 43  0.0000000  6.2831853}%
\special{pn 8}%
\special{ar 4026 442 42 43  0.0000000  6.2831853}%
%
\special{pn 8}%
\special{pa 4189 647}%
\special{pa 4165 575}%
\special{fp}%
\special{sh 1}%
\special{pa 4165 575}%
\special{pa 4167 645}%
\special{pa 4182 626}%
\special{pa 4205 632}%
\special{pa 4165 575}%
\special{fp}%
%
\special{pn 8}%
\special{pa 4026 629}%
\special{pa 4026 587}%
\special{fp}%
\special{sh 1}%
\special{pa 4026 587}%
\special{pa 4006 654}%
\special{pa 4026 640}%
\special{pa 4046 654}%
\special{pa 4026 587}%
\special{fp}%
%
\special{pn 8}%
\special{pa 3875 611}%
\special{pa 3887 581}%
\special{fp}%
\special{sh 1}%
\special{pa 3887 581}%
\special{pa 3844 635}%
\special{pa 3867 631}%
\special{pa 3881 650}%
\special{pa 3887 581}%
\special{fp}%
\put(42.4900,-7.6200){\makebox(0,0)[lb]{$e_1$}}%
\put(40.4400,-7.6200){\makebox(0,0)[lb]{$e_2$}}%
\put(38.1400,-7.6200){\makebox(0,0)[rb]{$e_3$}}%
%
\special{pn 8}%
\special{pa 1790 804}%
\special{pa 1862 804}%
\special{fp}%
\special{sh 1}%
\special{pa 1862 804}%
\special{pa 1795 784}%
\special{pa 1809 804}%
\special{pa 1795 824}%
\special{pa 1862 804}%
\special{fp}%
\special{pa 1488 1076}%
\special{pa 1488 1143}%
\special{fp}%
\special{sh 1}%
\special{pa 1488 1143}%
\special{pa 1508 1076}%
\special{pa 1488 1090}%
\special{pa 1468 1076}%
\special{pa 1488 1143}%
\special{fp}%
\special{pa 1192 804}%
\special{pa 1125 804}%
\special{fp}%
\special{sh 1}%
\special{pa 1125 804}%
\special{pa 1192 824}%
\special{pa 1178 804}%
\special{pa 1192 784}%
\special{pa 1125 804}%
\special{fp}%
\put(17.6000,-8.4700){\makebox(0,0)[rt]{$\widetilde{e}_1$}}%
\put(14.4500,-10.6400){\makebox(0,0)[rb]{$\widetilde{e}_2$}}%
\put(12.7000,-7.6800){\makebox(0,0)[rb]{$\widetilde{e}_3$}}%
\end{picture}}%
}
\end{center}
\caption{Hexagonal lattice and its quotient}
\label{figure:hex1}
\end{figure}
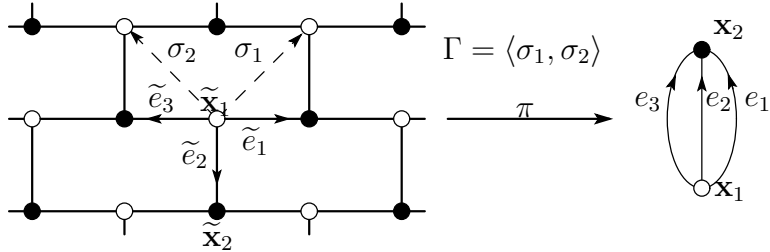

We discuss the modified standard realization 
of the hexagonal lattice $X=(V,E)$, where 
$V=\mathbb{Z}^2=\{ {\bf{x}}=(x_1, x_2) \, | \, x_1, x_2 \in \mathbb{Z} \}$ and 
\begin{equation*}
E=\left\{ 
({\bf{x}}, {\bf{y}}) \in V\times V \,  |  \, {\bf{x}}-{\bf{y}} =\pm (1,0), 
{\bf{x}}-{\bf{y}}=(0, (-1)^{x_1+x_2})
\right\}
\end{equation*}
(see Figure \ref{figure:hex1}).

We consider a random walk on $X$ given in the following manner:
If ${\bf{x}}=(x_1, x_2) \in V$ is a vertex so that $x_1+x_2$ is even, then we set
\begin{align*}
p({\bf{x}}, {\bf{x}}+(1,0))=\alpha,   \quad p({\bf{x}}, {\bf{x}}-(0,1))=\beta, \quad
p({\bf{x}}, {\bf{x}}-(1,0))=\gamma 
\end{align*}
with $\alpha,\beta,\gamma\geq 0$ and $\alpha+\beta+\gamma=1$. If $x_1+x_2$ is odd,
then we set
\begin{align*}
p({\bf{x}}, {\bf{x}}-(1,0))=\alpha^\prime,  \quad p({\bf{x}}, {\bf{x}}+(0,1))=\beta^\prime,
\quad p({\bf{x}}, {\bf{x}}+(1,0))=\gamma^\prime
\end{align*}
with $\alpha^\prime, \beta^\prime, \gamma^\prime\geq 0$ and
 $\alpha^\prime+\beta^\prime+\gamma^\prime=1$.
Here we see that $X$ is invariant under the action $\Gamma=\langle \sigma_1, \sigma_2 
\rangle(\cong \mathbb Z^{2}) $ generated by
$\sigma_1({\bf{x}})={\bf{x}}+(1,1)$ and $\sigma_2({\bf{x}})={\bf{x}}+(-1,1)$.
Then the quotient of $X$ by the action $\Gamma$ is a finite graph $X_0=(V_0, E_0)$ 
consisting of two vertices $V_0=\{ x_1, x_2 \}$ with three multiple edges 
$E=E_{x_1} \cup E_{x_2}=\{e_1, e_2, e_3\} \cup 
\{ \overline{e}_1, \overline{e}_2, \overline{e}_3 \}$.
The transition probability of the corresponding random walk on $X_0$ is given by
$$ p(e_{1})=\alpha, \quad p(\overline{e}_1)=\alpha^\prime,  \quad p(e_2)=\beta, 
\quad p(\overline{e}_2)=\beta^\prime, \quad 
p(e_3)=\gamma, \quad p(\overline{e}_3)=\gamma^\prime.$$

The first homology group $\mathrm{H}_1(X_0, \mathbb{R})$ is spanned by the two cycles
$[c_1]:=[e_1* \overline{e}_2]$ and $[c_2]:=[e_3 * \overline{e}_2 ]$. 
Solving (\ref{m-PF-intro}), we obtain $m(x_{1})=m(x_{2})=1/2$. Thus the homological direction $\gamma_{p}$
is given by
\begin{equation*}
\gamma_p=
\frac{1}{2}
\left\{ (\alpha-\alpha^\prime)e_1 +(\beta-\beta^\prime)e_2 
+(\gamma-\gamma^\prime) e_3 \right\}
=\frac{\check{\alpha}}{2}[c_1] +\frac{\check{\gamma}}{2} [c_2].
\end{equation*}
We define the canonical surjective linear map 
$\rho_\mathbb{R}: \mathrm{H}_1(X_0,\mathbb{R}) \rightarrow \Gamma \otimes \mathbb{R}$ by
$$
\rho_\mathbb{R} ( [c_1] ) := \sigma_1, \quad \rho_\mathbb{R} ([c_2]) := \sigma_2.
$$
Then the asymptotic direction is given by
\begin{equation}
\rho_\mathbb{R}( \gamma_p)= \frac{\check{\alpha}}{2} \sigma_1
+\frac{\check{\gamma}}{2} \sigma_2. 
\label{hex asymptotic}
\end{equation}

We are going to determine the modified standard realization $\Phi_0 :X \rightarrow 
(\Gamma\otimes \mathbb{R},g_{0})$.
Set $\widetilde{{\bf{x}}}_1=(0,0)$ and $\widetilde{{\bf{x}}}_2=(0,1)$ in $V$. 
Without loss of generality, we may put $\Phi_0(\widetilde{{\bf{x}}}_1)={\bf 0} \in \Gamma \otimes \mathbb{R}$.
It follows from equation (\ref{modified harmonic realization}) that
\begin{equation}
\Phi_0 ( \widetilde{{\bf{x}}}_1)-\alpha^\prime \sigma_1 -\gamma^\prime 
\sigma_2 =\Phi_0(\widetilde{{\bf{x}}}_2)+\rho_{\mathbb{R}} ( \gamma_p).
\label{7-3-1}
\end{equation} 
Combining (\ref{hex asymptotic}) with (\ref{7-3-1}), we obtain 
\begin{equation}
{\Phi_0}(\widetilde{{\bf{x}}}_1)={\bf 0}, \quad
{\Phi_0}(\widetilde{{\bf{x}}}_2)= -\frac{ \hat{\alpha}}{2} \sigma_1 -
\frac{ \hat{\gamma}}{2} \sigma_2.
\label{hex x-2}
\end{equation}
We define the basis $\{\omega_{1}, \omega_{2} \}$ in ${\rm Hom}(\Gamma, \mathbb R)$ by 
$$\omega_{1} \big( (m-n, m+n) \big)=m, \quad \omega_{2} \big( (m-n,m+n) \big )=n  \quad \mbox{for }(m,n)
\in \mathbb Z^{2}\cong \Gamma.$$
It is the dual basis of $\{ \sigma_{1}, \sigma_{2} \}$. Identifying 
$\omega_{i}\in {\rm Hom}(\Gamma, \mathbb R)$ ($i=1,2$) with 
$\hspace{-1mm}\mbox{ }^{t}\rho_{\mathbb R}(\omega_{i}) \in {\rm H}^{1}(X_{0}, \mathbb R)$
as before, we can also see that $\{ \omega_1, \omega_2\}$ is the dual basis
$\{ [c_1], [c_2] \} (\subset \mathrm{H}_1(X_0, \mathbb{R}))$. Recalling
$\mathrm{H}^1(X_0,\mathbb{R}) 
\cong \mathcal{H}^1(X_0) \subset C^1(X_0,\mathbb{R})$, we have
\begin{equation}
\begin{array}{l}
\omega_1(e_1)- \omega_1(e_2)=1, \qquad
\omega_1(e_3)- \omega_1(e_2)=0,  
\\
\omega_2(e_1)-\omega_2(e_2)=0, \qquad
\omega_2(e_3)-\omega_2(e_2)=1.
\label{alg-1}
\end{array}
\end{equation}
By definition of the modified harmonicity 
(\ref{M-harmonic}), we observe that $\omega_1$ and $\omega_2$ also satisfy 
\begin{align}
\hat{\alpha} \omega_1 (e_1) +\hat{\beta} \omega_1(e_2) +\hat{\gamma}\omega_1(e_3)= 0 \quad \mbox{and} \quad
\hat{\alpha} \omega_2 (e_1) +\hat{\beta} \omega_2(e_2) +\hat{\gamma}\omega_2(e_3)= 0,
\label{alg-2}
\end{align}
respectively.
Solving the algebraic equations (\ref{alg-1}) and (\ref{alg-2}), we obtain
\begin{align*}
\omega_1(e_1)&=1-\frac{\hat{\alpha}}{2}= \frac{\hat{\beta}+\hat{\gamma}}{2}, \quad 
\omega_1(e_2)=\omega_1(e_3)=-\frac{\hat{\alpha}}{2}, \\
\omega_2(e_1)&=\omega_2(e_2)= -\frac{\hat{\gamma}}{2}, \quad 
\quad \omega_2(e_3)=
1-\frac{\hat{\gamma}}{2}=
\frac{\hat{\alpha}+\hat{\beta}}{2}.
\end{align*}
Then by direct computation, we have 
\begin{align*}
\langle \! \langle \omega_1, \omega_1 \rangle \! \rangle=
\frac{ \hat{\alpha}(\hat{\beta}+\hat{\gamma}) - \check{\alpha}^2 }{4}, \quad
\langle \! \langle \omega_1, \omega_2 \rangle \! \rangle=
-\frac{\hat{\alpha}\hat{\gamma}+\check{\alpha}\check{\gamma}}{4}, \quad
\langle \! \langle \omega_2, \omega_2 \rangle \! \rangle=
\frac{(\hat{\alpha}+\hat{\beta})\hat{\gamma} - \check{\gamma}^2  }{4}.
\end{align*}
Thus the volume of the $\Gamma$-Jacobian torus $\mathrm{Jac}^\Gamma$ and
the Albanese metric $g_{0}$ on $\Gamma \otimes \mathbb R$ are given by
\begin{align*}
\mathrm{vol}(\mathrm{Jac}^\Gamma)&=\mathrm{vol}(\mathrm{Alb}^\Gamma)^{-1}
=\sqrt{ \det 
(\langle \! \langle \omega_i, \omega_j \rangle \! \rangle )_{i,j=1}^{2} }
\nonumber \\
&=
\frac{\sqrt{ 2\hat{\alpha} \hat{\beta} \hat{\gamma} +
\hat{\alpha}\check{\alpha}(\hat{\beta}\check{\gamma}+\check{\beta}\hat{\gamma})
+\hat{\beta}\check{\beta} (\hat{\alpha}\check{\gamma}+ \check{\alpha}\hat{\gamma})
+\hat{\gamma}\check{\gamma}( \hat{\alpha}\check{\beta} +\check{\alpha}\hat{\beta})}}{4},
\end{align*}
and
\begin{align}
\langle \sigma_1, \sigma_1 \rangle_{g_{0}} &= 
\frac{ \mathrm{vol}(\mathrm{Alb}^\Gamma)^{2}}{4} \big\{ 
( \hat{\alpha}+\hat{\beta})\hat{\gamma} - \check{\gamma}^2 \big \},
\nonumber \\
\langle \sigma_1, \sigma_2 \rangle_{g_{0}} 
&= 
\frac{ \mathrm{vol}(\mathrm{Alb}^\Gamma)^{2}}{4} \big (
 \hat{\alpha}\hat{\gamma} +\check{\alpha}\check{\gamma} \big),
\nonumber \\
\langle \sigma_2, \sigma_2 \rangle_{g_{0}} &= 
\frac{ \mathrm{vol}(\mathrm{Alb}^\Gamma)^{2}}{4} \big\{ 
 \hat{\alpha}(\hat{\beta}+\hat{\gamma}) - \check{\alpha}^2 \big \},
\label{7-3-2}
\nonumber
\end{align}
respectively. 

Let $\{ v_1, v_2 \}$ be the Gram--Schmidt orthogonalization of $\{ \omega_1, \omega_2 \}$, that is,
\begin{equation*}
v_1=\frac{ \omega_1}{\| \omega_1 \|}, \quad 
v_{2}={\rm vol}({\rm Alb}^{\Gamma}) \Vert \omega_{1} \Vert
\Big( \omega_2- \frac{ \langle \! \langle \omega_2, \omega_1 \rangle \! \rangle}
{\Vert \omega_1 \Vert^{2}} \omega_1 \Big).
\end{equation*}
Recalling that
$\{ \omega_{1}, \omega_{2} \}$ is the dual basis of $\{ [c_{1}], [c_{2}] \}$,
we obtain
\begin{align}
v_{1}([c_{i}])=\frac{\delta_{i1}}{ \| \omega_{1} \|},
\quad
v_{2}([c_{i}])=
{\rm vol}({\rm Alb}^{\Gamma}) \Vert \omega_{1} \Vert
\Big( \delta_{i2}- \frac{ \langle \! \langle \omega_2, \omega_1 \rangle \! \rangle}
{\Vert \omega_1 \Vert^{2}} \delta_{i1} \Big) \qquad (i=1,2),
\nonumber
\end{align}
and hence
\begin{equation}
\begin{array}{l}
\sigma_{1} = v_{1}([c_{1}]) {\bf v}_{1}+v_{2}([c_{1}]) {\bf v}_{2}
=\frac{1}{\| \omega_1 \|} {\bf v}_{1}
-\frac{ \langle \! \langle \omega_1, \omega_2 \rangle \! \rangle 
\mathrm{vol}(\mathrm{Alb}^\Gamma)} {\| \omega_1 \| } {\bf v}_{2},
\\ \\
\sigma_{2} =v_{2}([c_{1}]) {\bf v}_{1}+v_{2}([c_{2}]) {\bf v}_{2}=
\| \omega_1 \|  \mathrm{vol}( \mathrm{Alb}^\Gamma)  {\bf v}_{2}.
\label{7-3-5}
\end{array}
\end{equation}

\begin{figure}[htbp]
\begin{center}
\scalebox{1.2}{\input{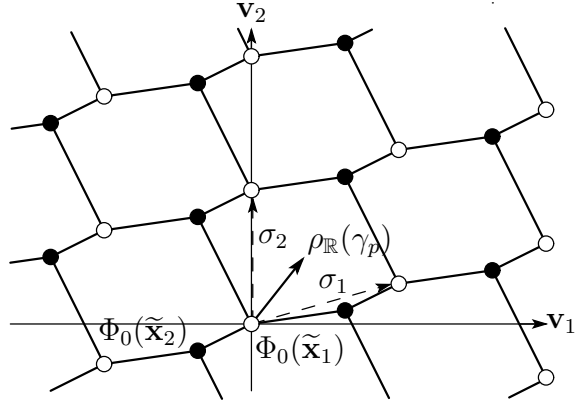}}
\end{center}
\caption{Modified standard realization of the hexagonal lattice}
\label{figure:hex3}
\end{figure}

Combining (\ref{hex x-2}) with (\ref{7-3-5}), we finally find that
the modified standard realization 
$\Phi_{0}: X \rightarrow (\Gamma\otimes \mathbb{R},g_{0}) \cong (\mathbb R^{2}, \{ {\bf v}_{1}, {\bf v}_{2} \})$ 
is given by $\Phi_{0}(\widetilde{{\bf{x}}}_1)=(0,0)$ and
\begin{align}
{\Phi}_{0}(\widetilde{{\bf{x}}}_2)&= \bigg(-
\frac{ \hat \alpha}{\sqrt{
\hat{\alpha}(\hat{\beta}+\hat{\gamma}) - \check{\alpha}^2
}}
\hspace{0.8mm}
,
\hspace{0.8mm}
%
%
\frac{ 
-2{\hat \alpha}{\hat \gamma}-{\hat \alpha}{\check \alpha}{\check \gamma}+{\check \alpha}^{2}{\hat \alpha}
}{4\sqrt{
\hat{\alpha}(\hat{\beta}+\hat{\gamma}) - \check{\alpha}^2
}} {\rm vol}({\rm Alb}^{\Gamma})
\bigg)
\nonumber
\end{align}
(see Figure \ref{figure:hex3}). 
In particular, when the random walk is simple (i.e., $\alpha=\alpha'=\beta=\beta'=\gamma=\gamma'=1/3$),
we have
\begin{align*}
& \mathrm{vol}(\mathrm{Alb}^\Gamma )=3{\sqrt 3}, \quad
\Phi_0 (\widetilde{{\bf{x}}}_2) = \Big( -\frac{\sqrt 2}{2}, -\frac{\sqrt 6}{2} \Big),
\quad
\sigma_{1}=\Big(\frac{3{\sqrt 2}}{2}, \frac{\sqrt 6}{2} \Big), \quad 
\sigma_{2}=\Big(0, {\sqrt 6}\Big).
\end{align*}
It means that the shape of the fundamental pattern of $\Phi_{0}(X)$ is
the equilateral hexagonal lattice with the common length of the sides 
${\sqrt 2}$.
\vspace{3mm} \\
{\bf{Acknowledgement.}}~The authors would like to thank Professor Toshikazu Sunada 
for sending the lecture slide \cite{Sunada-Lecture} 
during the preparation of the present paper. 
Also they are grateful to Professors Atsushi Katsuda, Kazumasa Kuwada, Stefan Neukamm, 
Tomoyuki Shirai and Ryokichi Tanaka 
for useful conversations and 
constant encouragement. Finally, the second author
thanks his graduate student Ryuya Namba for 
helpful discussion on Subsection 7.3.

\end{document}